\documentclass[11pt]{article}
\usepackage{epic,latexsym,amssymb}
\usepackage{color}
\usepackage{tikz}

\newtheorem{thm}{Theorem}
\newtheorem{lem}[thm]{Lemma}
\newtheorem{quest}{Question}
\newtheorem{fact}{Fact}
\newtheorem{claim}{Claim}

\newtheorem{conj}{Conjecture}

\newcommand{\qed}{$\Box$}

\newcommand{\gt}{\gamma_t}

\newcommand{\dbar}{{\overline{d}}}

\newcommand{\smallqed}{{\tiny ($\Box$)}}

\newcommand{\TR}[1]{\mbox{$\tau(#1)$}}

\newtheorem{definition}[fact]{Definition}

\newcommand{\cB}{{\cal B}}

\newcommand{\cH}{{\cal H}}

\newenvironment{unnumbered}[1]{\trivlist \item [\hskip \labelsep {\bf
#1}]\ignorespaces\it}{\endtrivlist}

\newcommand{\hyperedgetwo}[6]{
	\pgfmathsetmacro\Done{sqrt((#4-#1)^2+(#5-#2)^2)}
	\pgfmathsetmacro\angleone{(#2>#5)*(180+asin((#4-#1)/ \Done)-asin((#1-#4)/ \Done))+asin((#1-#4)/ \Done)+asin((#3-#6)/\Done)}
	\pgfmathsetmacro\angleone{\angleone-360*(\angleone>0)-360*(\angleone>360)}
	\draw ([shift=(\angleone:#3)] #1,#2)--([shift=(\angleone:#6)]#4,#5);
	\pgfmathsetmacro\Dtwo{sqrt((#1-#4)^2+(#2-#5)^2)}
	\pgfmathsetmacro\angletwo{(#5>#2)*(180+asin((#1-#4)/ \Dtwo)-asin((#4-#1)/ \Dtwo))+asin((#4-#1)/ \Dtwo)+asin((#6-#3)/\Dtwo)}
	\pgfmathsetmacro\angletwo{\angletwo-360*(\angletwo>0)-360*(\angletwo>360)}
	\draw ([shift=(\angletwo:#6)] #4,#5)--([shift=(\angletwo:#3)]#1,#2);
	\draw (#1,#2)+(\angletwo:#3) arc(\angletwo:\angleone+360*(\angleone<\angletwo):#3);
	\draw (#4,#5)+(\angleone:#6) arc(\angleone:\angletwo+360*(\angletwo<\angleone):#6);
}

\def\vertex(#1){\put(#1){\circle*{2}}}
\def\vertexo(#1){\put(#1){\circle{2}}}
\def\vert(#1){\put(#1){\circle*{1.5}}}
\def\verto(#1){\put(#1){\circle{1.5}}}
\def\lab(#1)#2{\put(#1){\makebox(0,0)[c]{#2}}}
\newcommand{\proof}{\noindent\textbf{Proof. }}
\newcommand{\2}{ \vspace{0.2cm} }
\newcommand{\1}{ \vspace{0.1cm} }

\newcommand{\PicTxt}[1]{{\small #1}}
\newcommand{\PicTxtII}[1]{{\normalsize #1}}
\newcommand{\Vx}[1]{V(#1)}

\newcommand{\NEW}[1]{{\color{red} #1}}

\begin{document}


\title{Transversals in $4$-Uniform Hypergraphs}




\author{Michael A. Henning${}^{1,}$\thanks{Research
supported in part by the South African
National Research Foundation and the University of Johannesburg} \, and Anders Yeo${}^{1,2}$\\
\\
${}^1$ Department of Pure and Applied Mathematics\\
University of Johannesburg \\
Auckland Park, 2006 South Africa \\
Email:  mahenning@uj.ac.za \\
\\
${}^2$ Engineering Systems and Design \\
Singapore University of Technology and Design \\
20 Dover Drive Singapore, 138682, Singapore \\
Email: andersyeo@gmail.com
}

\date{}
\maketitle

\begin{abstract}
Let $H$ be a $3$-regular $4$-uniform hypergraph on $n$ vertices. The transversal number $\tau(H)$ of $H$ is the minimum number of vertices that intersect every edge. Lai and Chang [J. Combin. Theory Ser. B 50 (1990), 129--133] proved that $\tau(H) \le 7n/18$. Thomass\'{e} and Yeo [Combinatorica 27 (2007), 473--487] improved this bound and showed that $\tau(H) \le 8n/21$. We provide a further improvement and prove that $\tau(H) \le 3n/8$, which is best possible due to a hypergraph of order eight. More generally, we show that if $H$ is a $4$-uniform hypergraph on $n$ vertices and $m$ edges with maximum degree $\Delta(H) \le 3$, then $\tau(H) \le n/4 + m/6$, which proves a known conjecture. We show that an easy corollary of our main result is that the total domination number of a graph on $n$ vertices with minimum degree at least~4 is at most $3n/7$,  which was the main result of the Thomass\'{e}-Yeo paper [Combinatorica 27 (2007), 473--487].
\end{abstract}


\newpage
\section{Notation and Definitions}

In this paper we continue the study of transversals in hypergraphs.
Hypergraphs are systems of sets which are conceived as natural
extensions of graphs.  A \emph{hypergraph} $H = (V,E)$ is a finite
set $V = V(H)$ of elements, called \emph{vertices}, together with a
finite multiset $E = E(H)$ of subsets of $V$, called
\emph{hyperedges} or simply \emph{edges}. The \emph{order} of $H$ is
$n(H) = |V|$ and the \emph{size} of $H$ is $m(H) = |E|$.

A $k$-\emph{edge} in $H$ is an edge of size~$k$. The hypergraph $H$
is said to be $k$-\emph{uniform} if every edge of $H$ is a $k$-edge.
Every (simple) graph is a $2$-uniform hypergraph. Thus graphs are
special hypergraphs.  For $i \ge 2$, we denote the number
of edges in $H$ of size $i$ by $e_i(H)$. The \emph{degree} of a
vertex $v$ in $H$, denoted by $d_H(v)$ or simply by $d(v)$ if $H$ is
clear from the context, is the number of edges of $H$ which contain
$v$. The maximum degree among the vertices of $H$ is denoted by
$\Delta(H)$. We say that two edges in $H$ \emph{overlap} if they intersect in at least two vertices.

Two vertices $x$ and $y$ of $H$ are \emph{adjacent} if there is an
edge $e$ of $H$ such that $\{x,y\}\subseteq e$. The
\emph{neighborhood} of a vertex $v$ in $H$, denoted $N_H(v)$ or
simply $N(v)$ if $H$ is clear from the context, is the set of all
vertices different from $v$ that are adjacent to $v$.
Two vertices $x$ and $y$ of $H$ are \emph{connected} if there is a
sequence $x=v_0,v_1,v_2\ldots,v_k=y$ of vertices of $H$ in which
$v_{i-1}$ is adjacent to $v_i$ for $i=1,2,\ldots,k$. A
\emph{connected hypergraph} is a hypergraph in which every pair of
vertices are connected. A maximal connected subhypergraph of $H$ is a
\emph{component} of $H$. Thus, no edge in $H$ contains vertices from
different components.

A subset $T$ of vertices in a hypergraph $H$ is a \emph{transversal}
(also called \emph{vertex cover} or \emph{hitting set} in many
papers) if $T$ has a nonempty intersection with every edge of $H$.
The \emph{transversal number} $\TR{H}$ of $H$ is the minimum size of
a transversal in $H$. A transversal of size~$\TR{H}$ is called a
$\TR{H}$-set.  Transversals in hypergraphs are well studied in the
literature (see, for
example,~\cite{BuHeTu12,ChMc,CoHeSl79,HeLo12,HeYe08,HeYe10,HeYe12a,HeYe13,LaCh90,ThYe07,Tu90}).

Given a hypergraph $H$ and subsets $X,Y \subseteq V(H)$ of vertices,
we let $H(X,Y)$ denote the hypergraph obtained by deleting all
vertices in $X \cup Y$ from $H$ and removing all edges containing
vertices from $X$ and removing the vertices in $Y$ from any remaining
edges. When we use the definition $H(X,Y)$ we furthermore assume that
no edges of size zero are created. That is, there is no edge $e \in
E(H)$ such that $\Vx{e} \subseteq Y \setminus X$. In this case we
note that if add $X$ to any $\tau(H(X,Y))$-set, then we get a
transversal of $H$, implying that $\tau(H) \le |X| + \tau(H(X,Y))$.
We will often use this fact throughout the paper.


A \emph{total dominating set}, also called a TD-set, of a graph $G$ with no isolated vertex is a set $S$ of vertices of $G$ such that every vertex is adjacent to a vertex in $S$. The \emph{total domination number} of $G$, denoted by $\gt(G)$, is the minimum cardinality of a TD-set of $G$. Total domination in graphs is now well studied in graph theory. The literature on the subject has been surveyed and detailed in the recent book~\cite{MHAYbookTD}. A recent paper on the topic can be found in~\cite{HeKlRa15}.

\section{The Family, $\cB$, of Hypergraphs}

In this section, we define a family, $\cB$, of ``bad" hypergraphs as
follows.

\begin{definition} \label{create_B}
Let $\cB$ be the class  of \textbf{bad hypergraphs} defined as
exactly those that can be generated using the operations (A)-(D)
below.

\begin{description}
  \item[(A):] Let $H_2$ be the hypergraph with two vertices
      $\{x,y\}$ and one edge $\{x,y\}$ and let $H_2$ belong to
      $\cB$.

  \item[(B):] Given any $B' \in \cB$ containing a $2$-edge
      $\{u,v\}$, define $B$ as follows. Let $V(B)=V(B') \cup
      \{x,y\}$ and let $E(B)=E(B') \cup \{\{u,v,x\}, \{u,v,y\},
      \{x,y\}\} \setminus \{u,v\}$.  Now add $B$ to $\cB$.

  \item[(C):] Given any $B' \in \cB$ containing a $3$-edge
      $\{u,v,w\}$, define $B$ as follows. Let $V(B)=V(B') \cup
      \{x,y\}$ and let $E(B)=E(B') \cup \{\{u,v,w,x\},
      \{u,v,w,y\}, \{x,y\}\} \setminus \{u,v,w\}$.  Now add $B$
      to $\cB$.

  \item[(D):] Given any $B_1,B_2 \in \cB$, such that $B_i$
      contains a $2$-edge $\{u_i,v_i\}$, for $i=1,2$, define $B$
      as follows. Let $V(B)=V(B_1) \cup V(B_2) \cup \{x\}$ and
      let $E(B)=E(B_1) \cup E(B_2) \cup \{\{u_1,v_1,x\},
      \{u_2,v_2,x\}, \{u_1,v_1,u_2,v_2\}\} \setminus
      \{\{u_1,v_1\},\{u_2,v_2\}\}$. Now add $B$ to $\cB$.
\end{description}

\begin{figure}[tb]
\begin{center}
\unitlength 0.40mm \linethickness{0.4pt}
\begin{picture}(80,100)
\put(50,40){\circle{11}} \put(50,40){\makebox(0,0){{\PicTxt{$u$}}}}
\put(50,80){\circle{11}} \put(50,80){\makebox(0,0){{\PicTxt{$v$}}}}
\drawline(40,40)(40.2185,37.9209)(40.8645,35.9326)(41.9098,34.1221)
\drawline(41.9098,34.1221)(43.3087,32.5686)(45,31.3397)(46.9098,30.4894)
\drawline(46.9098,30.4894)(48.9547,30.0548)(51.0453,30.0548)(53.0902,30.4894)
\drawline(53.0902,30.4894)(55,31.3397)(56.6913,32.5685)(58.0902,34.1221)
\drawline(58.0902,34.1221)(59.1354,35.9326)(59.7815,37.9209)(60,40)
\drawline(40,40)(40,80)
\drawline(60,80)(59.7815,82.0791)(59.1355,84.0674)(58.0902,85.8779)
\drawline(58.0902,85.8779)(56.6913,87.4314)(55,88.6602)(53.0902,89.5106)
\drawline(53.0902,89.5106)(51.0453,89.9452)(48.9547,89.9452)(46.9098,89.5106)
\drawline(46.9098,89.5106)(45,88.6603)(43.3087,87.4315)(41.9098,85.8779)
\drawline(41.9098,85.8779)(40.8646,84.0674)(40.2185,82.0791)(40,80)
\drawline(60,80)(60,40)
\put(50,15){\makebox(0,0){{\PicTxtII{$H_2$}}}}
\put(90,65){\makebox(0,0){{\huge $\Rightarrow$ }}}
\end{picture}
\unitlength 0.4mm
\begin{picture}(100,100)
\put(80,40){\circle{11}} \put(80,40){\makebox(0,0){{\PicTxt{$x$}}}}
\put(80,80){\circle{11}} \put(80,80){\makebox(0,0){{\PicTxt{$y$}}}}
\put(40,40){\circle{11}} \put(40,40){\makebox(0,0){{\PicTxt{$u$}}}}
\put(40,80){\circle{11}} \put(40,80){\makebox(0,0){{\PicTxt{$v$}}}}
\drawline(70,40)(70.2185,37.9209)(70.8645,35.9326)(71.9098,34.1221)
\drawline(71.9098,34.1221)(73.3087,32.5686)(75,31.3397)(76.9098,30.4894)
\drawline(76.9098,30.4894)(78.9547,30.0548)(81.0453,30.0548)(83.0902,30.4894)
\drawline(83.0902,30.4894)(85,31.3397)(86.6913,32.5685)(88.0902,34.1221)
\drawline(88.0902,34.1221)(89.1354,35.9326)(89.7815,37.9209)(90,40)
\drawline(70,40)(70,80)
\drawline(90,80)(89.7815,82.0791)(89.1355,84.0674)(88.0902,85.8779)
\drawline(88.0902,85.8779)(86.6913,87.4314)(85,88.6602)(83.0902,89.5106)
\drawline(83.0902,89.5106)(81.0453,89.9452)(78.9547,89.9452)(76.9098,89.5106)
\drawline(76.9098,89.5106)(75,88.6603)(73.3087,87.4315)(71.9098,85.8779)
\drawline(71.9098,85.8779)(70.8646,84.0674)(70.2185,82.0791)(70,80)
\drawline(90,80)(90,40)
\drawline(88.4853,71.5147)(89.9776,73.3332)(91.0866,75.4078)(91.7694,77.6589)
\drawline(91.7694,77.6589)(92,80)(91.7694,82.3411)(91.0865,84.5922)
\drawline(91.0865,84.5922)(89.9776,86.6669)(88.4853,88.4853)(86.6668,89.9776)
\drawline(86.6668,89.9776)(84.5922,91.0866)(82.3411,91.7694)(80,92)
\drawline(88.4853,71.5147)(48.4853,31.5147)
\drawline(28,40)(28.1966,37.8369)(28.7798,35.7447)(29.7306,33.7919)
\drawline(29.7306,33.7919)(31.0179,32.0425)(32.5994,30.5538)(34.4233,29.3745)
\drawline(34.4233,29.3745)(36.43,28.5433)(38.5536,28.0875)(40.7245,28.0219)
\drawline(40.7245,28.0219)(42.8718,28.3487)(44.925,29.0572)(46.8168,30.1242)(48.4853,31.5147)
\drawline(28,40)(28,80)
\drawline(40,92)(37.9162,91.8177)(35.8958,91.2763)(34,90.3923)
\drawline(34,90.3923)(32.2865,89.1925)(30.8075,87.7135)(29.6077,86)
\drawline(29.6077,86)(28.7237,84.1042)(28.1823,82.0838)(28,80)
\drawline(40,92)(80,92)
\drawline(80,26)(82.5236,26.2293)(84.9645,26.9098)(87.2427,28.0191)
\drawline(87.2427,28.0191)(89.2837,29.5209)(91.0206,31.3659)(92.3964,33.4939)
\drawline(92.3964,33.4939)(93.3661,35.835)(93.8979,38.3125)(93.9745,40.8453)
\drawline(93.9745,40.8453)(93.5932,43.3504)(92.7666,45.7458)(91.5218,47.9529)(89.8995,49.8995)
\drawline(80,26)(40,26)
\drawline(26,40)(26.2127,37.5689)(26.8443,35.2117)(27.8756,33)
\drawline(27.8756,33)(29.2754,31.001)(31.001,29.2754)(33,27.8756)
\drawline(33,27.8756)(35.2117,26.8443)(37.5689,26.2127)(40,26)
\drawline(26,40)(26,80)
\drawline(49.8995,89.8995)(47.778,91.6406)(45.3576,92.9343)(42.7313,93.731)
\drawline(42.7313,93.731)(40,94)(37.2687,93.731)(34.6424,92.9343)
\drawline(34.6424,92.9343)(32.222,91.6406)(30.1005,89.8995)(28.3594,87.778)
\drawline(28.3594,87.778)(27.0657,85.3576)(26.269,82.7313)(26,80)
\drawline(49.8995,89.8995)(89.8995,49.8995)
\put(60,15){\makebox(0,0){{\PicTxtII{$B$ }}}}
\end{picture}
\unitlength 0.25mm

\end{center}
\vskip -1 cm
\caption{An illustration of Step~(B) in  Definition~\ref{create_B}.}
\label{H2B4}
\end{figure}
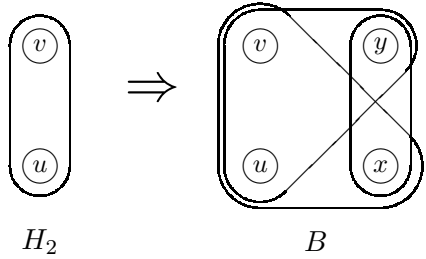


\begin{figure}[tb]
\begin{center}
 \unitlength 0.33mm
\linethickness{0.4pt}
\begin{picture}(320,160)
\put(80,40){\circle{11}}
\put(80,40){\makebox(0,0){{\PicTxt{$\bullet$}}}}
\put(80,80){\circle{11}}
\put(80,80){\makebox(0,0){{\PicTxt{$u$}}}}
\put(40,40){\circle{11}}
\put(40,40){\makebox(0,0){{\PicTxt{$w$}}}}
\put(40,80){\circle{11}}
\put(40,80){\makebox(0,0){{\PicTxt{$v$}}}}
\drawline(70,40)(70.2185,37.9209)(70.8645,35.9326)(71.9098,34.1221)
\drawline(71.9098,34.1221)(73.3087,32.5686)(75,31.3397)(76.9098,30.4894)
\drawline(76.9098,30.4894)(78.9547,30.0548)(81.0453,30.0548)(83.0902,30.4894)
\drawline(83.0902,30.4894)(85,31.3397)(86.6913,32.5685)(88.0902,34.1221)
\drawline(88.0902,34.1221)(89.1354,35.9326)(89.7815,37.9209)(90,40)
\drawline(70,40)(70,80)
\drawline(90,80)(89.7815,82.0791)(89.1355,84.0674)(88.0902,85.8779)
\drawline(88.0902,85.8779)(86.6913,87.4314)(85,88.6602)(83.0902,89.5106)
\drawline(83.0902,89.5106)(81.0453,89.9452)(78.9547,89.9452)(76.9098,89.5106)
\drawline(76.9098,89.5106)(75,88.6603)(73.3087,87.4315)(71.9098,85.8779)
\drawline(71.9098,85.8779)(70.8646,84.0674)(70.2185,82.0791)(70,80)
\drawline(90,80)(90,40)
\drawline(88.4853,71.5147)(89.9776,73.3332)(91.0866,75.4078)(91.7694,77.6589)
\drawline(91.7694,77.6589)(92,80)(91.7694,82.3411)(91.0865,84.5922)
\drawline(91.0865,84.5922)(89.9776,86.6669)(88.4853,88.4853)(86.6668,89.9776)
\drawline(86.6668,89.9776)(84.5922,91.0866)(82.3411,91.7694)(80,92)
\drawline(88.4853,71.5147)(48.4853,31.5147)
\drawline(28,40)(28.1966,37.8369)(28.7798,35.7447)(29.7306,33.7919)
\drawline(29.7306,33.7919)(31.0179,32.0425)(32.5994,30.5538)(34.4233,29.3745)
\drawline(34.4233,29.3745)(36.43,28.5433)(38.5536,28.0875)(40.7245,28.0219)
\drawline(40.7245,28.0219)(42.8718,28.3487)(44.925,29.0572)(46.8168,30.1242)(48.4853,31.5147)
\drawline(28,40)(28,80)
\drawline(40,92)(37.9162,91.8177)(35.8958,91.2763)(34,90.3923)
\drawline(34,90.3923)(32.2865,89.1925)(30.8075,87.7135)(29.6077,86)
\drawline(29.6077,86)(28.7237,84.1042)(28.1823,82.0838)(28,80)
\drawline(40,92)(80,92)
\drawline(80,26)(82.5236,26.2293)(84.9645,26.9098)(87.2427,28.0191)
\drawline(87.2427,28.0191)(89.2837,29.5209)(91.0206,31.3659)(92.3964,33.4939)
\drawline(92.3964,33.4939)(93.3661,35.835)(93.8979,38.3125)(93.9745,40.8453)
\drawline(93.9745,40.8453)(93.5932,43.3504)(92.7666,45.7458)(91.5218,47.9529)(89.8995,49.8995)
\drawline(80,26)(40,26)
\drawline(26,40)(26.2127,37.5689)(26.8443,35.2117)(27.8756,33)
\drawline(27.8756,33)(29.2754,31.001)(31.001,29.2754)(33,27.8756)
\drawline(33,27.8756)(35.2117,26.8443)(37.5689,26.2127)(40,26)
\drawline(26,40)(26,80)
\drawline(49.8995,89.8995)(47.778,91.6406)(45.3576,92.9343)(42.7313,93.731)
\drawline(42.7313,93.731)(40,94)(37.2687,93.731)(34.6424,92.9343)
\drawline(34.6424,92.9343)(32.222,91.6406)(30.1005,89.8995)(28.3594,87.778)
\drawline(28.3594,87.778)(27.0657,85.3576)(26.269,82.7313)(26,80)
\drawline(49.8995,89.8995)(89.8995,49.8995)
\put(60,15){\makebox(0,0){{\PicTxtII{$B'$ }}}}
\put(145,65){\makebox(0,0){{\huge $\Rightarrow$ }}}
\put(280,40){\circle{11}}
\put(280,40){\makebox(0,0){{\PicTxt{$\bullet$}}}}
\put(280,80){\circle{11}}
\put(280,80){\makebox(0,0){{\PicTxt{$u$}}}}
\put(240,40){\circle{11}}
\put(240,40){\makebox(0,0){{\PicTxt{$w$}}}}
\put(240,80){\circle{11}}
\put(240,80){\makebox(0,0){{\PicTxt{$v$}}}}
\put(240,120){\circle{11}}
\put(240,120){\makebox(0,0){{\PicTxt{$x$}}}}
\put(200,80){\circle{11}}
\put(200,80){\makebox(0,0){{\PicTxt{$y$}}}}
\drawline(246.364,113.636)(247.548,115.098)(248.402,116.775)(248.889,118.592)
\drawline(248.889,118.592)(248.988,120.471)(248.693,122.329)(248.019,124.086)
\drawline(248.019,124.086)(246.994,125.664)(245.664,126.994)(244.086,128.019)
\drawline(244.086,128.019)(242.329,128.693)(240.471,128.988)(238.592,128.889)
\drawline(238.592,128.889)(236.775,128.402)(235.098,127.548)(233.636,126.364)
\drawline(246.364,113.636)(206.364,73.636)
\drawline(193.636,86.364)(192.517,85.0001)(191.685,83.4442)(191.173,81.7558)
\drawline(191.173,81.7558)(191,80)(191.173,78.2442)(191.685,76.5558)
\drawline(191.685,76.5558)(192.517,74.9999)(193.636,73.636)(195,72.5168)
\drawline(195,72.5168)(196.556,71.6851)(198.244,71.1729)(200,71)
\drawline(200,71)(201.756,71.1729)(203.444,71.6851)(205,72.5168)(206.364,73.636)
\drawline(193.636,86.364)(233.636,126.364)
\drawline(270,40)(270.219,37.9209)(270.865,35.9326)(271.91,34.1221)
\drawline(271.91,34.1221)(273.309,32.5686)(275,31.3397)(276.91,30.4894)
\drawline(276.91,30.4894)(278.955,30.0548)(281.045,30.0548)(283.09,30.4894)
\drawline(283.09,30.4894)(285,31.3397)(286.691,32.5685)(288.09,34.1221)
\drawline(288.09,34.1221)(289.135,35.9326)(289.781,37.9209)(290,40)
\drawline(270,40)(270,80)
\drawline(290,80)(289.781,82.0791)(289.135,84.0674)(288.09,85.8779)
\drawline(288.09,85.8779)(286.691,87.4314)(285,88.6602)(283.09,89.5106)
\drawline(283.09,89.5106)(281.045,89.9452)(278.955,89.9452)(276.91,89.5106)
\drawline(276.91,89.5106)(275,88.6603)(273.309,87.4315)(271.91,85.8779)
\drawline(271.91,85.8779)(270.865,84.0674)(270.219,82.0791)(270,80)
\drawline(290,80)(290,40)
\drawline(280,29)(281.983,29.1802)(283.901,29.7148)(285.691,30.5864)
\drawline(285.691,30.5864)(287.294,31.7664)(288.659,33.2161)(289.74,34.8881)
\drawline(289.74,34.8881)(290.502,36.7275)(290.92,38.6741)(290.98,40.6642)
\drawline(290.98,40.6642)(290.68,42.6325)(290.031,44.5145)(289.053,46.2487)(287.778,47.7782)
\drawline(280,29)(240,29)
\drawline(229,40)(229.167,38.0899)(229.663,36.2378)(230.474,34.5)
\drawline(230.474,34.5)(231.574,32.9293)(232.929,31.5735)(234.5,30.4737)
\drawline(234.5,30.4737)(236.238,29.6634)(238.09,29.1671)(240,29)
\drawline(229,40)(229,80)
\drawline(247.778,87.7782)(246.111,89.1462)(244.21,90.1627)(242.146,90.7886)
\drawline(242.146,90.7886)(240,91)(237.854,90.7886)(235.79,90.1627)
\drawline(235.79,90.1627)(233.889,89.1462)(232.222,87.7782)(230.854,86.1113)
\drawline(230.854,86.1113)(229.837,84.2095)(229.211,82.146)(229,80)
\drawline(247.778,87.7782)(287.778,47.7782)
\drawline(200,93)(197.464,92.7502)(195.025,92.0104)(192.778,90.8091)
\drawline(192.778,90.8091)(190.808,89.1924)(189.191,87.2224)(187.99,84.9749)
\drawline(187.99,84.9749)(187.25,82.5362)(187,80)(187.25,77.4638)
\drawline(187.25,77.4638)(187.99,75.0251)(189.191,72.7776)(190.808,70.8076)
\drawline(200,93)(280,93)
\drawline(289.192,70.8076)(290.809,72.7776)(292.01,75.0251)(292.75,77.4638)
\drawline(292.75,77.4638)(293,80)(292.75,82.5362)(292.01,84.9749)
\drawline(292.01,84.9749)(290.809,87.2224)(289.192,89.1924)(287.222,90.8091)
\drawline(287.222,90.8091)(284.975,92.0104)(282.536,92.7502)(280,93)
\drawline(289.192,70.8076)(249.192,30.8076)
\drawline(230.808,30.8076)(232.544,29.351)(234.506,28.218)(236.635,27.443)
\drawline(236.635,27.443)(238.867,27.0495)(241.133,27.0495)(243.365,27.443)
\drawline(243.365,27.443)(245.494,28.218)(247.456,29.351)(249.192,30.8076)
\drawline(230.808,30.8076)(190.808,70.8076)
\drawline(250.607,130.607)(248.334,132.472)(245.74,133.858)(242.926,134.712)
\drawline(242.926,134.712)(240,135)(237.074,134.712)(234.26,133.858)
\drawline(234.26,133.858)(231.666,132.472)(229.393,130.607)(227.528,128.334)
\drawline(227.528,128.334)(226.142,125.74)(225.288,122.926)(225,120)
\drawline(250.607,130.607)(290.607,90.6066)
\drawline(290.607,69.3934)(292.287,71.3964)(293.595,73.6607)(294.489,76.1177)
\drawline(294.489,76.1177)(294.943,78.6927)(294.943,81.3073)(294.489,83.8823)
\drawline(294.489,83.8823)(293.595,86.3393)(292.287,88.6036)(290.607,90.6066)
\drawline(290.607,69.3934)(250.607,29.3934)
\drawline(225,40)(225.246,37.2962)(225.975,34.6809)(227.163,32.2399)
\drawline(227.163,32.2399)(228.772,30.0532)(230.749,28.1922)(233.029,26.7182)
\drawline(233.029,26.7182)(235.537,25.6792)(238.192,25.1094)(240.906,25.0274)
\drawline(240.906,25.0274)(243.59,25.4359)(246.156,26.3215)(248.521,27.6552)(250.607,29.3934)
\drawline(225,40)(225,120)
\put(240,15){\makebox(0,0){{\PicTxtII{$B$ }}}}
\end{picture}
 \unitlength 0.33mm
\linethickness{0.4pt}
\begin{picture}(360,140)
\put(80,40){\circle{13}}
\put(80,40){\makebox(0,0){{\PicTxt{$u_1$}}}}
\put(80,80){\circle{13}}
\put(80,80){\makebox(0,0){{\PicTxt{$v_1$}}}}
\put(40,40){\circle{13}}
\put(40,40){\makebox(0,0){{\PicTxt{$\bullet$}}}}
\put(40,80){\circle{13}}
\put(40,80){\makebox(0,0){{\PicTxt{$\bullet$}}}}
\drawline(70,40)(70.2185,37.9209)(70.8645,35.9326)(71.9098,34.1221)
\drawline(71.9098,34.1221)(73.3087,32.5686)(75,31.3397)(76.9098,30.4894)
\drawline(76.9098,30.4894)(78.9547,30.0548)(81.0453,30.0548)(83.0902,30.4894)
\drawline(83.0902,30.4894)(85,31.3397)(86.6913,32.5685)(88.0902,34.1221)
\drawline(88.0902,34.1221)(89.1354,35.9326)(89.7815,37.9209)(90,40)
\drawline(70,40)(70,80)
\drawline(90,80)(89.7815,82.0791)(89.1355,84.0674)(88.0902,85.8779)
\drawline(88.0902,85.8779)(86.6913,87.4314)(85,88.6602)(83.0902,89.5106)
\drawline(83.0902,89.5106)(81.0453,89.9452)(78.9547,89.9452)(76.9098,89.5106)
\drawline(76.9098,89.5106)(75,88.6603)(73.3087,87.4315)(71.9098,85.8779)
\drawline(71.9098,85.8779)(70.8646,84.0674)(70.2185,82.0791)(70,80)
\drawline(90,80)(90,40)
\drawline(88.4853,71.5147)(89.9776,73.3332)(91.0866,75.4078)(91.7694,77.6589)
\drawline(91.7694,77.6589)(92,80)(91.7694,82.3411)(91.0865,84.5922)
\drawline(91.0865,84.5922)(89.9776,86.6669)(88.4853,88.4853)(86.6668,89.9776)
\drawline(86.6668,89.9776)(84.5922,91.0866)(82.3411,91.7694)(80,92)
\drawline(88.4853,71.5147)(48.4853,31.5147)
\drawline(28,40)(28.1966,37.8369)(28.7798,35.7447)(29.7306,33.7919)
\drawline(29.7306,33.7919)(31.0179,32.0425)(32.5994,30.5538)(34.4233,29.3745)
\drawline(34.4233,29.3745)(36.43,28.5433)(38.5536,28.0875)(40.7245,28.0219)
\drawline(40.7245,28.0219)(42.8718,28.3487)(44.925,29.0572)(46.8168,30.1242)(48.4853,31.5147)
\drawline(28,40)(28,80)
\drawline(40,92)(37.9162,91.8177)(35.8958,91.2763)(34,90.3923)
\drawline(34,90.3923)(32.2865,89.1925)(30.8075,87.7135)(29.6077,86)
\drawline(29.6077,86)(28.7237,84.1042)(28.1823,82.0838)(28,80)
\drawline(40,92)(80,92)
\drawline(80,26)(82.5236,26.2293)(84.9645,26.9098)(87.2427,28.0191)
\drawline(87.2427,28.0191)(89.2837,29.5209)(91.0206,31.3659)(92.3964,33.4939)
\drawline(92.3964,33.4939)(93.3661,35.835)(93.8979,38.3125)(93.9745,40.8453)
\drawline(93.9745,40.8453)(93.5932,43.3504)(92.7666,45.7458)(91.5218,47.9529)(89.8995,49.8995)
\drawline(80,26)(40,26)
\drawline(26,40)(26.2127,37.5689)(26.8443,35.2117)(27.8756,33)
\drawline(27.8756,33)(29.2754,31.001)(31.001,29.2754)(33,27.8756)
\drawline(33,27.8756)(35.2117,26.8443)(37.5689,26.2127)(40,26)
\drawline(26,40)(26,80)
\drawline(49.8995,89.8995)(47.778,91.6406)(45.3576,92.9343)(42.7313,93.731)
\drawline(42.7313,93.731)(40,94)(37.2687,93.731)(34.6424,92.9343)
\drawline(34.6424,92.9343)(32.222,91.6406)(30.1005,89.8995)(28.3594,87.778)
\drawline(28.3594,87.778)(27.0657,85.3576)(26.269,82.7313)(26,80)
\drawline(49.8995,89.8995)(89.8995,49.8995)
\put(60,15){\makebox(0,0){{\PicTxtII{$B_1$ }}}}
\put(120,40){\circle{13}}
\put(120,40){\makebox(0,0){{\PicTxt{$u_2$}}}}
\put(120,80){\circle{13}}
\put(120,80){\makebox(0,0){{\PicTxt{$v_2$}}}}
\drawline(110,40)(110.219,37.9209)(110.865,35.9326)(111.91,34.1221)
\drawline(111.91,34.1221)(113.309,32.5686)(115,31.3397)(116.91,30.4894)
\drawline(116.91,30.4894)(118.955,30.0548)(121.045,30.0548)(123.09,30.4894)
\drawline(123.09,30.4894)(125,31.3397)(126.691,32.5685)(128.09,34.1221)
\drawline(128.09,34.1221)(129.135,35.9326)(129.781,37.9209)(130,40)
\drawline(110,40)(110,80)
\drawline(130,80)(129.781,82.0791)(129.135,84.0674)(128.09,85.8779)
\drawline(128.09,85.8779)(126.691,87.4314)(125,88.6602)(123.09,89.5106)
\drawline(123.09,89.5106)(121.045,89.9452)(118.955,89.9452)(116.91,89.5106)
\drawline(116.91,89.5106)(115,88.6603)(113.309,87.4315)(111.91,85.8779)
\drawline(111.91,85.8779)(110.865,84.0674)(110.219,82.0791)(110,80)
\drawline(130,80)(130,40)
\put(120,15){\makebox(0,0){{\PicTxtII{$B_2$ }}}}
\put(182,60){\makebox(0,0){{\huge $\Rightarrow$ }}}
\put(280,40){\circle{13}}
\put(280,40){\makebox(0,0){{\PicTxt{$u_1$}}}}
\put(280,80){\circle{13}}
\put(280,80){\makebox(0,0){{\PicTxt{$v_1$}}}}
\put(240,40){\circle{13}}
\put(240,40){\makebox(0,0){{\PicTxt{$\bullet$}}}}
\put(240,80){\circle{13}}
\put(240,80){\makebox(0,0){{\PicTxt{$\bullet$}}}}
\put(320,40){\circle{13}}
\put(320,40){\makebox(0,0){{\PicTxt{$u_2$}}}}
\put(320,80){\circle{13}}
\put(320,80){\makebox(0,0){{\PicTxt{$v_2$}}}}
\put(300,120){\circle{13}}
\put(300,120){\makebox(0,0){{\PicTxt{$x$}}}}
\drawline(270,40)(270.152,38.2635)(270.603,36.5798)(271.34,35)
\drawline(271.34,35)(272.34,33.5721)(273.572,32.3396)(275,31.3397)
\drawline(275,31.3397)(276.58,30.6031)(278.264,30.1519)(280,30)
\drawline(270,40)(270,80)
\drawline(280,90)(278.264,89.8481)(276.58,89.3969)(275,88.6603)
\drawline(275,88.6603)(273.572,87.6604)(272.34,86.4279)(271.34,85)
\drawline(271.34,85)(270.603,83.4202)(270.152,81.7365)(270,80)
\drawline(280,90)(320,90)
\drawline(330,80)(329.848,81.7365)(329.397,83.4202)(328.66,85)
\drawline(328.66,85)(327.66,86.4279)(326.428,87.6604)(325,88.6603)
\drawline(325,88.6603)(323.42,89.3969)(321.736,89.8481)(320,90)
\drawline(330,80)(330,40)
\drawline(320,30)(321.736,30.1519)(323.42,30.6031)(325,31.3397)
\drawline(325,31.3397)(326.428,32.3396)(327.66,33.5721)(328.66,35)
\drawline(328.66,35)(329.397,36.5798)(329.848,38.2635)(330,40)
\drawline(320,30)(280,30)
\drawline(288.485,71.5147)(289.978,73.3332)(291.087,75.4078)(291.769,77.6589)
\drawline(291.769,77.6589)(292,80)(291.769,82.3411)(291.087,84.5922)
\drawline(291.087,84.5922)(289.978,86.6669)(288.485,88.4853)(286.667,89.9776)
\drawline(286.667,89.9776)(284.592,91.0866)(282.341,91.7694)(280,92)
\drawline(288.485,71.5147)(248.485,31.5147)
\drawline(228,40)(228.197,37.8369)(228.78,35.7447)(229.731,33.7919)
\drawline(229.731,33.7919)(231.018,32.0425)(232.599,30.5538)(234.423,29.3745)
\drawline(234.423,29.3745)(236.43,28.5433)(238.554,28.0875)(240.725,28.0219)
\drawline(240.725,28.0219)(242.872,28.3487)(244.925,29.0572)(246.817,30.1242)(248.485,31.5147)
\drawline(228,40)(228,80)
\drawline(240,92)(237.916,91.8177)(235.896,91.2763)(234,90.3923)
\drawline(234,90.3923)(232.287,89.1925)(230.807,87.7135)(229.608,86)
\drawline(229.608,86)(228.724,84.1042)(228.182,82.0838)(228,80)
\drawline(240,92)(280,92)
\drawline(280,26)(282.524,26.2293)(284.964,26.9098)(287.243,28.0191)
\drawline(287.243,28.0191)(289.284,29.5209)(291.021,31.3659)(292.396,33.4939)
\drawline(292.396,33.4939)(293.366,35.835)(293.898,38.3125)(293.974,40.8453)
\drawline(293.974,40.8453)(293.593,43.3504)(292.767,45.7458)(291.522,47.9529)(289.899,49.8995)
\drawline(280,26)(240,26)
\drawline(226,40)(226.213,37.5689)(226.844,35.2117)(227.876,33)
\drawline(227.876,33)(229.275,31.001)(231.001,29.2754)(233,27.8756)
\drawline(233,27.8756)(235.212,26.8443)(237.569,26.2127)(240,26)
\drawline(226,40)(226,80)
\drawline(249.899,89.8995)(247.778,91.6406)(245.358,92.9343)(242.731,93.731)
\drawline(242.731,93.731)(240,94)(237.269,93.731)(234.642,92.9343)
\drawline(234.642,92.9343)(232.222,91.6406)(230.101,89.8995)(228.359,87.778)
\drawline(228.359,87.778)(227.066,85.3576)(226.269,82.7313)(226,80)
\drawline(249.899,89.8995)(289.9,49.8995)
\drawline(312.522,126.261)(311.073,128.567)(309.204,130.549)(306.987,132.132)
\drawline(306.987,132.132)(304.506,133.255)(301.854,133.877)(299.132,133.973)
\drawline(299.132,133.973)(296.443,133.541)(293.889,132.596)(291.565,131.174)
\drawline(291.565,131.174)(289.561,129.329)(287.952,127.131)(286.8,124.664)
\drawline(286.8,124.664)(286.146,122.019)(286.018,119.299)(286.418,116.604)
\drawline(312.522,126.261)(332.522,86.261)
\drawline(334,80)(333.94,81.2964)(333.76,82.5816)(333.462,83.8446)
\drawline(333.462,83.8446)(333.048,85.0746)(332.522,86.261)
\drawline(334,80)(334,40)
\drawline(306.418,36.6045)(307.322,34.0611)(308.697,31.7385)(310.493,29.723)
\drawline(310.493,29.723)(312.642,28.0895)(315.064,26.8988)(317.67,26.1952)
\drawline(317.67,26.1952)(320.363,26.0047)(323.042,26.3345)(325.608,27.1723)
\drawline(325.608,27.1723)(327.966,28.487)(330.027,30.2297)(331.716,32.3356)
\drawline(331.716,32.3356)(332.969,34.7264)(333.74,37.3133)(334,40)
\drawline(306.418,36.6045)(286.418,116.604)
\drawline(315.522,116.119)(315.98,119.198)(315.833,122.308)(315.086,125.33)
\drawline(315.086,125.33)(313.769,128.15)(311.93,130.662)(309.64,132.77)
\drawline(309.64,132.77)(306.985,134.395)(304.065,135.475)(300.992,135.969)
\drawline(300.992,135.969)(297.881,135.859)(294.85,135.149)(292.014,133.865)
\drawline(292.014,133.865)(289.481,132.056)(287.346,129.791)(285.689,127.155)
\drawline(315.522,116.119)(295.522,36.1194)
\drawline(264,40)(264.297,36.9294)(265.179,33.973)(266.611,31.2407)
\drawline(266.611,31.2407)(268.541,28.834)(270.896,26.8423)(273.591,25.3398)
\drawline(273.591,25.3398)(276.523,24.3823)(279.585,24.0054)(282.662,24.2231)
\drawline(282.662,24.2231)(285.641,25.0272)(288.409,26.388)(290.865,28.2548)
\drawline(290.865,28.2548)(292.917,30.5583)(294.489,33.2127)(295.522,36.1194)
\drawline(264,40)(264,80)
\drawline(265.689,87.1554)(265.088,85.7995)(264.615,84.3938)(264.274,82.9504)
\drawline(264.274,82.9504)(264.069,81.4815)(264,80)
\drawline(265.689,87.1554)(285.689,127.155)
\put(280,15){\makebox(0,0){{\PicTxtII{$B$ }}}}
\end{picture}

\end{center}
\vskip -0.4 cm   \caption{An illustration of Steps~(C) and~(D) in  Definition~\ref{create_B}.}
\label{B6B7}
\end{figure}
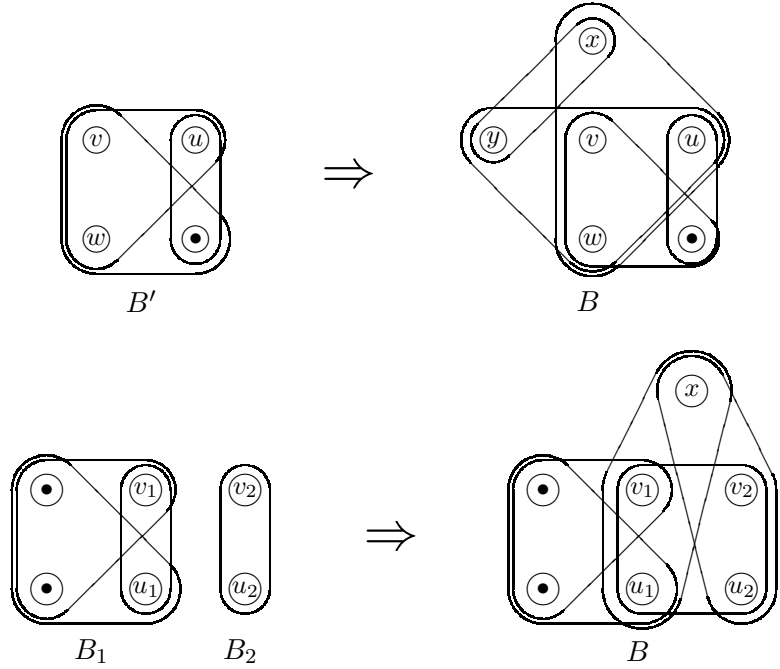

\vskip -0.3 cm We call the two vertices, $\{x,y\}$, added in step (A)
above an $\mathbf{(A)}$-\textbf{pair}.
Note that in operations (B) and (C), $\{a,b\}$ is an $(A)$-pair in $B$
if and only if it is an $(A)$-pair in $B'$. Analogously in operation (D),
$\{a,b\}$ is an $(A)$-pair in $B$ if and only if it is an $(A)$-pair in $B_1$ or $B_2$.

\end{definition}

The hypergraph $B \in \cB$ created by applying Step~(B) in
Definition~\ref{create_B} to the hypergraph $H_2$ is shown in
Figure~\ref{H2B4}, while Figure~\ref{B6B7} illustrates Step~(C) and
Step~(D) in Definition~\ref{create_B}.

We shall need the following definition.

\begin{definition}
If $H$ is a hypergraph, then let $b(H)$ denote the number of
components in $H$ that belong to $\cB$. Further for $i \ge 0$, let
$b^i(H)$ denote the maximum number of vertex disjoint subhypergraphs
in $H$ which are isomorphic to hypergraphs in $\cB$ and which are
intersected by exactly~$i$ other edges of $H$.
\end{definition}

\section{Main Results}

Let $\cH$ denote the class of all hypergraphs where all edges have
size at most four and at least two and with maximum degree at most
three. We shall prove the following result a proof of which is
presented in Section~\ref{S:proof}.

\begin{thm} \label{main_thm}
If $H \in \cH$, then
\[
24 \tau(H) \le 6n(H) + 4e_4(H) + 6e_3(H) + 10e_2(H) + 2b(H) + b^1(H).
\]
Furthermore if $b^1(H)$ is odd, then the above inequality is strict.
\end{thm}

Let $H$ be a $4$-uniform hypergraph with $\Delta(H) \le 3$. Since
every hypergraph in $\cB$ contains a $2$-edge or a $3$-edge, we note
that $b(H)=b^1(H)=0$. By the $4$-uniformity of $H$, we have that
$e_2(H)=e_3(H)=0$ and $e_4(H)=m(H)$. Therefore,
Theorem~\ref{main_thm} implies that $24\tau(H) \le 6n(H) + 4 m(H)$.
Hence as an immediate consequence of Theorem~\ref{main_thm} we have
our two main results.

\begin{thm} \label{main_thmA}
If $H$ is a $4$-uniform hypergraph with $\Delta(H) \le 3$, then
$\tau(H) \le \frac{n(H)}{4} + \frac{m(H)}{6}$.
\end{thm}


\begin{thm} \label{main_thmB}
If $H$ is a $3$-regular $4$-uniform hypergraph on $n$ vertices, then $\tau(H) \le \frac{3n(H)}{8}$.
\end{thm}

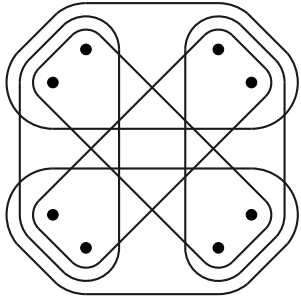
\begin{figure}[htb]
\begin{center}
\tikzstyle{vertexX}=[circle,draw, fill=black!100, minimum size=8pt, scale=0.5, inner sep=0.1pt]
\begin{tikzpicture}[scale=0.44]
 \draw (0,0) node {\mbox{ }};
\draw (4,3) node {\mbox{ }};
\node (a1) at (0.0,1.0) [vertexX] {};
\node (a2) at (1.0,0.0) [vertexX] {};
\node (b1) at (0.0,5.0) [vertexX] {};
\node (b2) at (1.0,6.0) [vertexX] {};
\node (c1) at (5.0,6.0) [vertexX] {};
\node (c2) at (6.0,5.0) [vertexX] {};
\node (d1) at (6.0,1.0) [vertexX] {};
\node (d2) at (5.0,0.0) [vertexX] {};
\draw[color=black!90, thick,rounded corners=4pt] (0.33828113841207375,-0.7497520578289725) arc (228.5689049390958:359.9147461611031:1.0); 
\draw[color=black!90, thick,rounded corners=4pt] (-0.99998850591691,0.9952054161738438) arc (180.27471047027626:221.4310950609044:1.0); 
\draw[color=black!90, thick,rounded corners=4pt] (-0.7497520578289725,5.661718861587926) arc (138.5689049390958:179.72528952972374:1.0); 
\draw[color=black!90, thick,rounded corners=4pt] (1.9999988929874426,6.001487959639652) arc (0.0852538388969144:131.4310950609044:1.0); 
\draw[color=black!90, thick,rounded corners=4pt] (0.33828113841207375,-0.7497520578289725) -- (-0.7497520578289704,0.3382811384120713);
\draw[color=black!90, thick,rounded corners=4pt] (-0.99998850591691,0.9952054161738438) -- (-0.99998850591691,5.004794583826157);
\draw[color=black!90, thick,rounded corners=4pt] (-0.7497520578289725,5.661718861587926) -- (0.3382811384120713,6.74975205782897);
\draw[color=black!90, thick,rounded corners=4pt] (1.9999988929874426,6.001487959639652) -- (1.9999988929874426,-0.0014879596396527753);
\draw[color=black!90, thick,rounded corners=4pt] (5.661718861587927,6.749752057828973) arc (48.56890493909578:179.91474616110307:1.0); 
\draw[color=black!90, thick,rounded corners=4pt] (6.99998850591691,5.004794583826157) arc (0.2747104702762737:41.4310950609044:1.0); 
\draw[color=black!90, thick,rounded corners=4pt] (6.749752057828972,0.3382811384120734) arc (318.5689049390958:359.72528952972374:1.0); 
\draw[color=black!90, thick,rounded corners=4pt] (4.000001107012557,-0.0014879596396528518) arc (180.08525383889693:311.43109506090445:1.0); 
\draw[color=black!90, thick,rounded corners=4pt] (5.661718861587927,6.749752057828973) -- (6.74975205782897,5.661718861587929);
\draw[color=black!90, thick,rounded corners=4pt] (6.99998850591691,5.004794583826157) -- (6.99998850591691,0.9952054161738435);
\draw[color=black!90, thick,rounded corners=4pt] (6.749752057828972,0.3382811384120734) -- (5.6617188615879295,-0.7497520578289695);
\draw[color=black!90, thick,rounded corners=4pt] (4.000001107012557,-0.0014879596396528518) -- (4.000001107012557,6.001487959639653);
\draw[color=black!90, thick,rounded corners=4pt] (-1.090271767160063,5.878241125052609) arc (141.14768868443804:269.77422548308033:1.4); 
\draw[color=black!90, thick,rounded corners=4pt] (0.9830127200795282,7.399896936320993) arc (90.69523094951481:128.85231130303504:1.4); 
\draw[color=black!90, thick,rounded corners=4pt] (5.878241125052609,7.090271767160063) arc (51.14768868443803:89.30476905048522:1.4); 
\draw[color=black!90, thick,rounded corners=4pt] (6.005516697885407,3.600010869312036) arc (-89.77422548308033:38.85231130303504:1.4); 
\draw[color=black!90, thick,rounded corners=4pt] (-1.090271767160063,5.878241125052609) -- (0.12175887518576367,7.090271767352078);
\draw[color=black!90, thick,rounded corners=4pt] (0.9830127200795282,7.399896936320993) -- (5.016987279920471,7.399896936320993);
\draw[color=black!90, thick,rounded corners=4pt] (5.878241125052609,7.090271767160063) -- (7.090271767352078,5.878241124814236);
\draw[color=black!90, thick,rounded corners=4pt] (6.005516697885407,3.600010869312036) -- (-0.005516697885407651,3.600010869312036);
\draw[color=black!90, thick,rounded corners=4pt] (0.1217588749473899,-1.0902717671600624) arc (231.14768868443804:269.3047690504852:1.4); 
\draw[color=black!90, thick,rounded corners=4pt] (-0.005516697885406376,2.399989130687964) arc (90.22577451691966:218.85231130303504:1.4); 
\draw[color=black!90, thick,rounded corners=4pt] (7.090271767160062,0.1217588749473898) arc (-38.85231131556196:89.77422548308027:1.4); 
\draw[color=black!90, thick,rounded corners=4pt] (5.016987279920471,-1.3998969363209934) arc (270.6952309495148:308.85231130303504:1.4); 
\draw[color=black!90, thick,rounded corners=4pt] (0.1217588749473899,-1.0902717671600624) -- (-1.0902717673520785,0.12175887518576367);
\draw[color=black!90, thick,rounded corners=4pt] (-0.005516697885406376,2.399989130687964) -- (6.005516697885408,2.399989130687964);
\draw[color=black!90, thick,rounded corners=4pt] (7.090271767160062,0.1217588749473898) -- (5.878241124814236,-1.0902717673520785);
\draw[color=black!90, thick,rounded corners=4pt] (5.016987279920471,-1.3998969363209934) -- (0.9830127200795284,-1.3998969363209934);
\draw[color=black!90, thick,rounded corners=4pt] (0.5843075616526987,-0.43266591811798066) arc (226.1461890060021:314.98844089708825:0.6); 
\draw[color=black!90, thick,rounded corners=4pt] (-0.42434965297924815,1.42417846717672) arc (135.01155910291172:223.8538109939979:0.6); 
\draw[color=black!90, thick,rounded corners=4pt] (5.4156924383473015,6.432665918117981) arc (46.14618900600209:134.98844089708828:0.6); 
\draw[color=black!90, thick,rounded corners=4pt] (6.424349652979248,4.57582153282328) arc (-44.98844089708825:43.85381099399791:0.6); 
\draw[color=black!90, thick,rounded corners=4pt] (0.5843075616526987,-0.43266591811798066) -- (-0.4326659181179808,0.5843075616526989);
\draw[color=black!90, thick,rounded corners=4pt] (-0.42434965297924815,1.42417846717672) -- (4.57582153282328,6.424349652979248);
\draw[color=black!90, thick,rounded corners=4pt] (5.4156924383473015,6.432665918117981) -- (6.432665918117981,5.4156924383473015);
\draw[color=black!90, thick,rounded corners=4pt] (6.424349652979248,4.57582153282328) -- (1.4241784671767195,-0.4243496529792485);
\draw[color=black!90, thick,rounded corners=4pt] (-0.43266591811798066,5.4156924383473015) arc (136.1461890060021:224.98844089708828:0.6); 
\draw[color=black!90, thick,rounded corners=4pt] (1.4241784671767197,6.424349652979249) arc (45.011559102911725:133.8538109939979:0.6); 
\draw[color=black!90, thick,rounded corners=4pt] (6.432665918117981,0.5843075616526987) arc (-43.85381099399791:44.98844089708825:0.6); 
\draw[color=black!90, thick,rounded corners=4pt] (4.57582153282328,-0.42434965297924815) arc (225.01155910291172:313.8538109939979:0.6); 
\draw[color=black!90, thick,rounded corners=4pt] (-0.43266591811798066,5.4156924383473015) -- (0.5843075616526989,6.432665918117981);
\draw[color=black!90, thick,rounded corners=4pt] (1.4241784671767197,6.424349652979249) -- (6.424349652979248,1.42417846717672);
\draw[color=black!90, thick,rounded corners=4pt] (6.432665918117981,0.5843075616526987) -- (5.4156924383473015,-0.4326659181179808);
\draw[color=black!90, thick,rounded corners=4pt] (4.57582153282328,-0.42434965297924815) -- (-0.4243496529792482,4.57582153282328);
 \end{tikzpicture}
\end{center}
\vskip -0.5cm
\caption{The $3$-regular $4$-uniform hypergraph, $H_8$.}
\label{hyper8}
\end{figure}

Theorem~\ref{main_thmA} and Theorem~\ref{main_thmB} are best possible due to the hypergraph, $H_8$, depicted in Figure~\ref{hyper8}, of order~$n = 8$, size~$m = 6$, satisfying $\tau(H_8) = 3 = \frac{3n}{8} = \frac{n}{4} + \frac{m}{6}$. 

As an application of our main result, Theorem~\ref{main_thmA}, we give very short proofs of the following three known results in Section~\ref{3n7}.

\begin{thm} \label{main_thmC_h1} {\rm (\cite{ChMc})}
If $H$ is a $4$-uniform hypergraph, then
$\tau(H) \le \frac{n(H)}{6} + \frac{m(H)}{3}$.
\end{thm}

\begin{thm} \label{main_thmC_h2} {\rm (\cite{ThYe07})}
If $H$ is a $4$-uniform hypergraph, then
$\tau(H) \le \frac{5n(H)}{21} + \frac{4m(H)}{21}$.
\end{thm}

Recall that $\delta(G)$ denotes the minimum degree of a graph $G$.

\begin{thm} \label{main_thmC} {\rm (\cite{ThYe07})}
If $G$ is a graph of order~$n$ with $\delta(G) \ge 4$, then $\gt(G) \le \frac{3n}{7}$.
\end{thm}

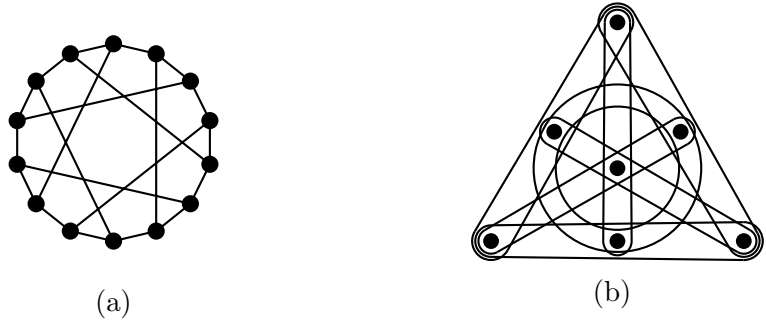
\begin{figure}[htb]
\tikzstyle{every node}=[circle, draw, fill=black!0, inner sep=0pt,minimum width=.2cm]
\begin{center}
\begin{tikzpicture}[thick,scale=.7]
  \draw(0,0) { 
    +(1.01,3.56) -- +(1.83,3.75)
    +(1.83,3.75) -- +(2.64,3.56)
    +(2.64,3.56) -- +(3.29,3.04)
    +(3.29,3.04) -- +(3.66,2.29)
    +(3.66,2.29) -- +(3.66,1.46)
    +(3.66,1.46) -- +(3.29,0.71)
    +(3.29,0.71) -- +(2.64,0.19)
    +(2.64,0.19) -- +(1.83,0.00)
    +(1.83,0.00) -- +(1.01,0.19)
    +(1.01,0.19) -- +(0.36,0.71)
    +(0.36,0.71) -- +(0.00,1.46)
    +(0.00,1.46) -- +(0.00,2.29)
    +(0.00,2.29) -- +(0.36,3.04)
    +(1.01,3.56) -- +(0.36,3.04)
    +(1.01,3.56) -- +(3.66,1.46)
    +(2.64,3.56) -- +(2.64,0.19)
    +(3.66,2.29) -- +(1.01,0.19)
    +(3.29,0.71) -- +(0.00,1.46)
    +(1.83,0.00) -- +(0.36,3.04)
    +(0.36,0.71) -- +(1.83,3.75)
    +(0.00,2.29) -- +(3.29,3.04)
    +(1.83,3.75) node[fill=black!100]{}
    +(2.64,3.56) node[fill=black!100]{}
    +(3.29,3.04) node[fill=black!100]{}
    +(3.66,2.29) node[fill=black!100]{}
    +(3.66,1.46) node[fill=black!100]{}
    +(3.29,0.71) node[fill=black!100]{}
    +(2.64,0.19) node[fill=black!100]{}
    +(1.83,0.00) node[fill=black!100]{}
    +(1.01,0.19) node[fill=black!100]{}
    +(0.36,0.71) node[fill=black!100]{}
    +(0.00,1.46) node[fill=black!100]{}
    +(0.00,2.29) node[fill=black!100]{}
    +(0.36,3.04) node[fill=black!100]{}
    +(1.01,3.56) node[fill=black!100]{}
    +(1.83,-1.2) node[rectangle, draw=white!0, fill=white!100]{(a)}   
    +(11.3,-1) node[rectangle, draw=white!0, fill=white!100]{(b)}   
  };
  { \begin{scope}[xshift=9cm,yshift=0cm,scale=.6]
\fill (0,0) circle (0.25cm); \fill (4,0) circle (0.25cm); \fill (8,0)
circle (0.25cm); \fill (2,3.464) circle (0.25cm); \fill (6,3.464)
circle (0.25cm); \fill (4,6.928) circle (0.25cm); \fill (4,2.3093)
circle (0.25cm); \hyperedgetwo{0}{0}{0.5}{8}{0}{0.6};
\hyperedgetwo{4}{6.928}{0.5}{0}{0}{0.6};
\hyperedgetwo{4}{6.928}{0.6}{8}{0}{0.5};
\hyperedgetwo{0}{0}{0.4}{6}{3.464}{0.45};
\hyperedgetwo{4}{6.928}{0.4}{4}{0}{0.45};
\hyperedgetwo{8}{0}{0.4}{2}{3.464}{0.45}; \draw (4,2.3093) circle
(1.9523cm); \draw (4,2.3093) circle (2.6523cm); 
{\large }; 
\end{scope}
};
\end{tikzpicture}
\end{center}
\vskip -0.6 cm \caption{The Heawood graph and the Fano Plane.}\index{Heawood graph} \label{Ch6:f:Heawood}
\end{figure}

Theorem~\ref{main_thmC_h2} and Theorem~\ref{main_thmC} were the main results in \cite{ThYe07}. 
Recall that the Heawood graph is the graph shown in Figure~\ref{Ch6:f:Heawood}(a) (which is the unique $6$-cage). The bipartite complement of the Heawood graph is the bipartite graph formed by taking the two partite sets of the Heawood graph and joining a vertex from one partite set to a vertex from the other partite set by an edge whenever they are not joined in the Heawood graph.  The bipartite complement of the Heawood graph can also be seen as the incidence bipartite graph of the complement of the Fano plane which is shown in Figure~\ref{Ch6:f:Heawood}(b).

\medskip


In fact, it is not difficult to prove the following improvement on Theorem~\ref{main_thmC} along the same lines as the proof of Theorem~\ref{main_thmC}. A proof of Theorem~\ref{main_thmD} is provided in Section~\ref{char:3n7}.

\begin{thm}{\rm (\cite{MHAYbookTD,AYimprove3n7})}
\label{main_thmD}
If $G$ is a connected graph  of order~$n$ with $\delta(G) \ge 4$,
then $\gt(G) \le 3n/7$. Furthermore we have equality if and only if
$G$ is the bipartite complement of the Heawood Graph.
\end{thm}

\subsection{Motivation}

There has been much interest in determining upper bounds on the
transversal number of a $3$-regular $4$-uniform hypergraph. In
particular, as a consequence of more general results we have the
Chv\'{a}tal-McDiarmid bound, the improved Lai-Chang bound and the
further improved Thomass\'{e}-Yeo bound. These bounds are summarized
in Theorem~\ref{t:knownbds}.

\begin{thm}
Let $H$ be a $3$-regular $4$-uniform hypergraph on $n$ vertices. Then the following bounds on $\tau(H)$ have been established.  \\
\indent {\rm (a)} $\tau(H) \le 5n/12$ {\rm (Chv\'{a}tal, McDiarmid~\cite{ChMc}).}  \\
\indent {\rm (b)} $\tau(H) \le 7n/18$ {\rm (Lai, Chang~\cite{LaCh90}).} \\
\indent {\rm (c)} $\tau(H) \le 8n/21$ {\rm (Thomass\'{e}, Yeo~\cite{ThYe07}).}
 \label{t:knownbds}
\end{thm}

In this paper, we provide a further improvement on the bounds in Theorem~\ref{t:knownbds} as shown in our main result, Theorem~\ref{main_thmB}, by proving that $\tau(H) \le 3n/8$. As mentioned above our bound is best possible, due to a hypergraph on eight vertices.

Motivated by comments and questions posed by Douglas
West~\cite{West}, the authors in~\cite{HeYe13} considered the
following slightly more general question.

\begin{quest}
For $k \ge 2$, let $H$ be a hypergraph on $n$ vertices with $m$ edges
and with every edge of size at least~$k$. Is it true that $\tau(H)
\le n/k + m/6$ holds for all $k$?
 \label{Q1}
\end{quest}

It is shown in~\cite{HeYe13} that Question~\ref{Q1} holds for $k =
2$ and a characterization of the extremal hypergraphs is given.
Chv\'{a}tal and McDiarmid~\cite{ChMc} proved that Question~\ref{Q1}
holds for $k = 3$ and the extremal hypergraphs are characterized
in~\cite{HeYe13}. Question~\ref{Q1} is not always true when $k \ge
4$ as shown in~\cite{HeYe13}. However the family of counterexamples
presented in~\cite{HeYe13} all satisfy $\Delta(H) \ge 4$. The
authors in~\cite{HeYe13} pose the following conjecture.

\begin{conj}{\rm (\cite{HeYe13})}
For all $k \ge 2$, if $H$ is a $k$-uniform hypergraph on $n$ vertices
with $m$ edges satisfying $\Delta(H) \le 3$, then $\tau(H) \le n/k +
m/6$.
 \label{conj1}
\end{conj}

As remarked earlier, Conjecture~\ref{conj1} always holds when $k \in
\{2,3\}$ (with no restriction on the maximum degree).
In~\cite{HeYe13} it is shown that Conjecture~\ref{conj1} is true
when $\Delta(H) \le 2$. However Conjecture~\ref{conj1} appears to be
a challenging conjecture for general $k \ge 4$ and for $\Delta(H) =
3$. We remark that if the conjecture is true, then this would imply
as a very special case a long standing conjecture due to Tuza and
Vestergaard~\cite{TuVe02} that if $H$ is a $3$-regular $6$-uniform
hypergraph, then $\tau(H) \le n/4$. In this paper, we prove that Conjecture~\ref{conj1} is true for $4$-uniform hypergraphs as shown in our main result, Theorem~\ref{main_thmA}.

\section{Preliminary Lemma}

We need the following lemma which proves a number of properties of
the hypergraphs that belong to the family $\cB$.

\begin{lem} \label{B_properties}
The following properties holds for all $B \in \cB$.
\begin{description}
\item[(i):] If $B$ was created from $B'$ in Step (B) or (C) in
    Definition~\ref{create_B}, then $\tau(B) = \tau(B') + 1$.
\item[(ii):] If $B$ was created from $B_1$ and $B_2$ in Step (D)
    in Definition~\ref{create_B}, then $\tau(B) = \tau(B_1) +
    \tau(B_2)$.
\item[(iii):] $\tau(B) = (6n(B) + 4e_4(B) + 6 e_3(B) + 10 e_2(B)
    +2)/24$.
\item[(iv):] All $(A)$-pairs are vertex disjoint (recall the definition of $(A)$-pairs, below the definition of Steps (A)-(D)).
\item[(v):] For all $e \in E(B)$ we have $\tau(B - e) =
    \tau(B)-1$.
\item[(vi):] For all $s \in V(B)$ there exists a $\tau(B)$-set
    containing $s$.
\item[(vii):] For all $s,t \in V(B)$ there exists a $\tau(B)$-set
    containing both $s$ and $t$ if and only if $\{s,t\}$ is not
    an $(A)$-pair.
\item[(viii):] Let $\{s_1,t_1\}$, $\{s_2,t_2\}$ and $\{s_3,t_3\}$
    be three subsets of $V(B)$. Then there exists a $\tau(B)$-set
    in $B$ intersecting all of these three sets.
\item[(ix):] There is no $4$-edge in $B$ intersecting three or
    more $2$-edges.
\item[(x):] If $B \ne H_2$, then $\delta(B) \ge 2$.
\item[(xi):] If $d_B(x)=2$, then $x$ is contained in a $3$-edge
    or a $2$-edge in $B$.
\item[(xii):] If $B\ne H_2$ and $e_2(B)>0$, then $B$ contains
    either two overlapping $3$-edges or two $4$-edges, $e_1$ and
    $e_2$, with $|V(e_1) \cap V(e_2)|=3$.
\item[(xiii):] If $B\ne H_2$ and $B$ does not contain two
    $4$-edges intersecting in three vertices, then every $2$-edge
    in $B$ intersects two overlapping $3$-edges.
\end{description}
\end{lem}
\proof \textbf{(i):} Suppose that $B$ was created from $B'$ in Step
(B) in Definition~\ref{create_B}. Name the vertices as in
Definition~\ref{create_B} and let $S$ be a $\tau(B)$-set. Since the
set $S$ intersects the $2$-edge $\{x,y\}$, we note that $|S \cap
\{x,y\}| \ge 1$. If $|S \cap \{x,y\}|=2$, then $(S \cup \{u\})
\setminus \{x\}$ is a $\tau(B)$-set. Hence we may choose the set $S$
so that $|S \cap \{x,y\}|=1$. This implies that $|S \cap \{u,v\}| \ge
1$ and that $S \setminus \{x,y\}$ is a transversal in $B'$ of size
$|S|-1$, and so $\tau(B') \le \tau(B) - 1$. Since every transversal
in $B'$ can be extended to a transversal in $B$ by adding to it the
vertex $x$, we have that $\tau(B) \le \tau(B') + 1$. Consequently,
$\tau(B) = \tau(B') + 1$, as desired. If $B$ was created from $B'$ in
Step~(C) in Definition~\ref{create_B}, then analogously to when $B$
was created in Step (B), we have that $\tau(B) = \tau(B') + 1$.

\2 \textbf{(ii):} Suppose that $B$ was created from $B_1$ and $B_2$
in Step~(D). Name the vertices as in Definition~\ref{create_B} and
let $S$ be a $\tau(B)$-set. Suppose $x \in S$. Since $S \cap
\{u_1,v_1,u_2,v_2\} \ne \emptyset$, we may assume, renaming vertices
if necessary, that $u_1 \in S$. Then, $(S \cup \{u_2\}) \setminus
\{x\}$ is a $\tau(B)$-set. Hence we may choose the set $S$ so that $x
\notin S$. In this case, $S \cap V(B_1)$ is a transversal in $B_1$
and $S \cap V(B_2)$ is a transversal in $B_2$, and so $\tau(B_1) +
\tau(B_2) \le |S \cap V(B_1)| + |S \cap V(B_2)| = |S| = \tau(B)$.
Furthermore, if $S_i$ is a transversal of $B_i$, for $i=1,2$, then
$S_1 \cup S_2$ is a transversal of $B$, and so $\tau(B) \le \tau(B_1)
+ \tau(B_2)$. Consequently, $\tau(B) = \tau(B_1) + \tau(B_2)$.

\2 \textbf{(iii):} We will show Part~(iii) by induction on the order,
$n(B)$, of the hypergraph $B$. If $n(G) = 2$, then $B=H_2$ was
created in step (A) in Definition~\ref{create_B}. In this case,
$\tau(B) = 1 = (12 + 10 + 2)/24 = (6n(B) + 4e_4(B) + 6 e_3(B) + 10
e_2(B) +2)/24$ and Part (iii) holds in this case. This establishes
the base case. Let $k \ge 3$ and assume that the formula holds for
all $B' \in \cB$ with $n(B') < k$ and let $B \in \cB$ have
order~$n(B) = k$.

Suppose that $B$ was created from $B'$ in Step (B) in
Definition~\ref{create_B}. By Part~(i), $\tau(B) = \tau(B') + 1$.
Applying the inductive hypothesis to $B'$, we therefore have that

\[
\begin{array}{lcl} \2
  \tau(B) &  = &  \tau(B')+1 \\   \2
   & = &  \frac{1}{24}( 6n(B') + 4e_4(B') + 6 e_3(B') + 10 e_2(B') + 2  ) + 1 \\  \2
   & = &  \frac{1}{24}( 6(n(B)-2) + 4e_4(B) + 6(e_3(B)-2) + 10 e_2(B) + 2  ) + 1 \\  \2
   & = &  \frac{1}{24}( 6n(B) + 4e_4(B) + 6 e_3(B) + 10 e_2(B) +2 ),
\end{array}
\]
and so Part (iii) holds in this case.
Suppose next that $B$ was created from $B'$ in Step~(C) in
Definition~\ref{create_B}. By Part~(i), $\tau(B) = \tau(B') + 1$.
Applying the inductive hypothesis to $B'$, we therefore have that

\[
\begin{array}{lcl} \2
  \tau(B) &  = &  \tau(B')+1 \\  \2
   & = &  \frac{1}{24}( 6n(B') + 4e_4(B') + 6 e_3(B') + 10 e_2(B') + 2  ) + 1 \\ \2
   & = &  \frac{1}{24}( 6(n(B)-2) + 4(e_4(B)-2) + 6(e_3(B)+1) + 10 (e_2(B) - 1) + 2  ) + 1 \\  \2
   & = &  \frac{1}{24}( 6n(B) + 4e_4(B) + 6 e_3(B) + 10 e_2(B) +2 ),
\end{array}
\]
and so Part (iii) holds in this case.
Suppose finally that $B$ was created from $B_1$ and $B_2$ in
Step~(D). By Part~(ii), $\tau(B) = \tau(B_1) + \tau(B_2)$. Applying
the inductive hypothesis to $B_1$ and $B_2$, we therefore have that
\[
\begin{array}{lcl} \2
  \tau(B) &  = &  \tau(B_1) + \tau(B_2) \\   \2
   & = &  \frac{1}{24}( 6(n(B)-1) + 4(e_4(B)-1) + 6(e_3(B)-2) + 10(e_2(B)+2) +2+2 ) \\  \2
   & = &  \frac{1}{24}( 6n(B) + 4e_4(B) + 6 e_3(B) + 10 e_2(B) +2 ),
\end{array}
\]
and so Part (iii) holds in this case. This completes the proof of
Part (iii).

\2 \textbf{(iv):} Part~(iv) follows easily by induction as no
operation can make $(A)$-pairs intersect.

\2 \textbf{(v):} We will prove Part~(v) by induction on the order,
$n(B)$, of the hypergraph $B$. Let $e \in E(B)$ be an arbitrary edge
in $B$. If $n(G) = 2$, then $B=H_2$ was created in step (A) in
Definition~\ref{create_B}. In this case, if $e$ denotes the edge of
$B$, then $\tau(B-e) = 0 = \tau(B)-1$ and Part (v) holds. This
establishes the base case. Let $k \ge 3$ and assume that the result
holds for all $B' \in \cB$ with $n(B') < k$ and let $B \in \cB$ have
order~$n(B) = k$. Let $e \in E(B)$ be an arbitrary edge in $B$.

Suppose that $B$ was created from $B'$ in Step (B) in
Definition~\ref{create_B} and name the vertices as in
Definition~\ref{create_B}. By Part~(i), $\tau(B') = \tau(B)-1$.
Suppose that $e=\{u,v,x\}$ or $e = \{u,v,y\}$. Renaming vertices, if
necessary, we may assume without loss of generality that
$e=\{u,v,x\}$. By induction there exists a $\tau(B' - \{u,v\})$-set,
$S'$, with $|S'|=\tau(B')-1$. Since $S' \cup \{y\}$ is a transversal
of $B-e$, we note that $\tau(B-e) \le |S'| + 1 = \tau(B') =
\tau(B)-1$, and so $\tau(B-e) \le \tau(B)-1$. Since deleting an edge
from a hypergraph can decrease the transversal number by at most one,
we have that $\tau(B-e) \ge \tau(B)-1$. Consequently, $\tau(B-e) =
\tau(B) - 1$, as desired. Suppose next that $e = \{x,y\}$. In this
case any transversal in $B'$ is a transversal in $B-e$, implying that
$\tau(B-e) \le \tau(B') = \tau(B)-1$. As observed earlier, $\tau(B-e)
\ge \tau(B) - 1$. Consequently, $\tau(B-e) = \tau(B) - 1$, as
desired. Suppose finally that $e \in E(B')$. By induction,
$\tau(B'-e) = \tau(B')-1$. Every $\tau(B'-e)$-set can be extended to
a transversal of $B - e$ by adding to it the vertex~$x$, implying
that $\tau(B) - 1 \le \tau(B-e) \le \tau(B'-e) + 1 = \tau(B') =
\tau(B)-1$. Consequently, $\tau(B-e) = \tau(B) - 1$, as desired.

If $B$ was created from $B'$ in Step~(C) in
Definition~\ref{create_B}, then the proof that Part~(v) holds is
analogous to when $B$ was created in Step~(B).

Suppose finally that $B$ was created from $B_1$ and $B_2$ in Step~(D)
and name the vertices as in Definition~\ref{create_B}. By Part~(ii),
$\tau(B) = \tau(B_1) + \tau(B_2)$. Suppose first that
$e=\{u_i,v_i,x\}$ for some $i \in \{1,2\}$. By induction there exists
a $\tau(B_i - \{u_i,v_i\})$-set, $S_i$, with $|S_i|=\tau(B_i)-1$. Let
$S_{3-i}$ be any $\tau(B_{3-i})$-set in $B_{3-i}$ and note that $S_1
\cup S_2$ is a transversal in $B-e$, and so $\tau(B) - 1 \le
\tau(B-e) \le |S_1| + |S_2| = \tau(B_1)+\tau(B_2)-1 = \tau(B) - 1$.
Consequently, $\tau(B-e) = \tau(B)-1$, as desired.
Suppose next that $e=\{u_1,v_1,u_2,v_2\}$. By induction there exists
a $\tau(B_i - \{u_i,v_i\})$-set, $T_i$, with $|T_i|=\tau(B_i)-1$.
Then, $T_1 \cup T_2 \cup \{x\}$ is a transversal in $B-e$, and so
$\tau(B) - 1 \le \tau(B-e) \le |T_1| + |T_2| + 1 = (\tau(B_1)-1) +
(\tau(B_2)-1) + 1 = \tau(B_1) + \tau(B_2) - 1 = \tau(B) - 1$.
Consequently, $\tau(B-e) = \tau(B)-1$, as desired.
Suppose finally that $e \in E(B_i)$ for some $i\in \{1,2\}$. By
induction there exists a $\tau(B_i - e)$-set, $D_i$, with
$|D_i|=\tau(B_i)-1$. Let $D_{3-i}$ be any $\tau(B_{3-i})$-set in
$B_{3-i}$ and note that $D_1 \cup D_2$ is a transversal in $B-e$, and
so $\tau(B) - 1 \le \tau(B-e) \le |D_1| + |D_2| =
\tau(B_1)+\tau(B_2)-1 = \tau(B) - 1$. Consequently, $\tau(B-e) =
\tau(B)-1$, as desired. This completes the proof of Part (v).

\2 \textbf{(vi):} We will prove Part~(vi) by induction on the order,
$n(B)$, of the hypergraph $B$. Let $s \in V(B)$ be an arbitrary
vertex in $B$. If $n(G) = 2$, then $B=H_2$ was created in step (A) in
Definition~\ref{create_B}. In this case, there clearly exists a
$\tau(B)$-set containing $s$ and Part (vi) holds. This establishes
the base case. Let $k \ge 3$ and assume that the result holds for all
$B' \in \cB$ with $n(B') < k$ and let $B \in \cB$ have order~$n(B) =
k$. Let $s \in V(B)$ be an arbitrary vertex in $B$.

Suppose that $B$ was created from $B'$ in Step (B) in
Definition~\ref{create_B} and name the vertices as in
Definition~\ref{create_B}. By Part~(i), $\tau(B') = \tau(B)-1$. On
the one hand, if $s \in \{x,y\}$, then adding the vertex $s$ to any
$\tau(B')$-set produces a transversal in $B$ of size $\tau(B') +1 =
\tau(B)$ containing $s$. On the other hand, if $s \notin \{x,y\}$,
then by induction let $S$ be any $\tau(B')$-set containing $s$ and
note that $S \cup \{x\}$ is a transversal of size $\tau(B')+1 =
\tau(B)$ in $B$ containing~$s$. In both cases, there exists a
$\tau(B)$-set containing~$s$.

If $B$ was created from $B'$ in Step~(C) in
Definition~\ref{create_B}, then the proof that Part~(vi) holds is
analogous to when $B$ was created in Step~(B).

Suppose finally that $B$ was created from $B_1$ and $B_2$ in Step~(D)
and name the vertices as in Definition~\ref{create_B}. By Part~(ii),
$\tau(B) = \tau(B_1) + \tau(B_2)$. Suppose first that $s=x$. In this
case, let $S_1$ be any $\tau(B_1)$-set and let $S_2$ be any $\tau(B_2
- \{u_2,v_2\})$-set. By Part~(v), $|S_2| = \tau(B_2)-1$. Thus the set
$S_1 \cup S_2 \cup \{x\}$ is a transversal in $B$ containing $s$ of
size $\tau(B_1) + (\tau(B_2)-1) +1 = \tau(B)$, as desired. Suppose
next that $s \ne x$. Renaming $B_1$ and $B_2$, if necessary, we may
assume that $s \in V(B_1)$. Applying the inductive hypothesis to
$B_1$, there exists a $\tau(B_1)$-set, $S_1$, containing $s$. Let
$S_2$ be a $\tau(B_2)$-set. Then, $S_1 \cup S_2$ is a transversal in
$B$ containing $s$ of size $\tau(B_1)+\tau(B_2) = \tau(B)$, which
completes the proof of Part~(vi).

\2 \textbf{(vii):} We will prove Part~(vii) by induction on the
order, $n(B)$, of the hypergraph $B$. Let $s,t \in V(B)$ be distinct
arbitrary vertices. If $n(G) = 2$, then $B=H_2$ was created in step~(A) in Definition~\ref{create_B}. In this case, $\{s,t\}$ is an
$(A)$-pair and there is no $\tau(B)$-set containing both $s$ and $t$.
This establishes the base case. Let $k \ge 3$ and assume that the
result holds for all $B' \in \cB$ with $n(B') < k$ and let $B \in
\cB$ have order~$n(B) = k$. Let $s,t \in V(B)$ be distinct arbitrary
vertices.

Suppose that $B$ was created from $B'$ in Step (B) in
Definition~\ref{create_B} and name the vertices as in
Definition~\ref{create_B}. By Part~(i), $\tau(B') = \tau(B)-1$.
Suppose first that $\{s,t\} = \{x,y\}$. Let $S'$ be any $\tau(B' -
\{u,v\})$-set. By Part~(v), $|S'| = \tau(B')-1$. The set $S' \cup
\{s,t\}$ is a transversal in $B$ of size $(\tau(B')-1) + 2 = \tau(B)$
containing $s$ and $t$, as desired. Suppose next that $|\{s,t\} \cap
\{x,y\}| = 1$. By Part~(vi) there exists a $\tau(B')$-set, $S''$, containing
the vertex in the set $\{s,t\}\setminus \{x,y\}$. Adding the vertex in
$\{s,t\} \cap \{x,y\}$ to $S''$ produces a transversal of size $\tau(B')+1=\tau(B)$ in $B$
containing $s$ and $t$, as desired. Finally consider the case when
$\{s,t\} \cap \{x,y\} = \emptyset$. If there exists a $\tau(B')$-set
containing both $s$ and $t$, then add $x$ to such a set in order to
obtain a $\tau(B)$-set containing $s$ and $t$. If there is no
$\tau(B')$-set containing both $s$ and $t$, then, by induction,
$\{s,t\}$ is an $(A)$-pair in $B'$ and therefore also an $(A)$-pair
in $B$.

We will now show that if $\{s,t\}$ is an $(A)$-pair in $B'$ (and therefore in $B$)
there is no $\tau(B)$-set containing $s$ and $t$.
For the sake of contradiction, assume that $S$ is a $\tau(B)$-set
containing $s$ and $t$.
If $S \cap V(B')$ is a transversal in $B'$, then since there is no
$\tau(B')$-set containing both $s$ and $t$ and $\{s,t\} \subseteq S
\cap V(B')$, we have that $\tau(B') < |S \cap V(B')|$. However since
$|S \cap \{x,y\}| \ge 1$, this implies that $|S| \ge |S \cap V(B')| +
1 > \tau(B') + 1 = \tau(B)$, a contradiction. Hence, the set $S \cap
V(B')$ is not a transversal in $B'$. The only edge of $B'$ that does
not intersect $S$ is the edge $\{u,v\}$, implying that $\{u,v\} \cap
S = \emptyset$ and $\{x,y\} \subseteq S$. In this case, $|S \cap
V(B')| = |S| - 2 = \tau(B) - 2 = \tau(B') - 1$. Hence adding the
vertex $v$ to the set $S \cap V(B')$ produces a transversal in $B'$
of size~$\tau(B')$ containing both $s$ and $t$, a contradiction.
Therefore if $\{s,t\}$ is an $(A)$-pair in $B'$, then there is no
$\tau(B)$-set containing $s$ and $t$.

If $B$ was created from $B'$ in Step~(C) in
Definition~\ref{create_B}, then the proof that Part~(vii) holds is
analogous to when $B$ was created in Step~(B).

Suppose finally that $B$ was created from $B_1$ and $B_2$ in Step~(D)
and name the vertices as in Definition~\ref{create_B}. By Part~(ii),
$\tau(B) = \tau(B_1) + \tau(B_2)$.  Suppose $x \in \{s,t\}$. Without
loss of generality we assume that $x=s$ and $t \in V(B_1)$. By
Part~(vi) there exists a $\tau(B_1)$-set, $S_1$, containing the
vertex $t$. Let $S_2$ be a $\tau(B_2 - \{u_2,v_2\})$-set. By
part~(v), $|S_2| = \tau(B_2) - 1$. Now the set $S_1 \cup S_2 \cup
\{x\}$ is a transversal in $B$ containing $s$ and $t$ of size $|S_1|
+ |S_2| + 1 = \tau(B_1) + (\tau(B_2)-1) +1 = \tau(B)$. Hence we may
assume that $x \notin \{s,t\}$, for otherwise the desired result
follows. Suppose $|\{s,t\} \cap V(B_1)|=1$. Renaming vertices if
necessary, we may assume that $s \in V(B_1)$ and $t \in V(B_2)$. By
Part~(vi) there exists a $\tau(B_1)$-set, $S_1$, containing the
vertex $s$ and a $\tau(B_2)$-set, $S_2$, containing the vertex $t$.
In this case, the set $S_1 \cup S_2$ is a transversal in $B$
containing $s$ and $t$ of size $|S_1| + |S_2| = \tau(B_1) + \tau(B_2)
= \tau(B)$. Hence without loss of generality we may assume that
$\{s,t\} \subseteq V(B_1)$.

If there exists a $\tau(B_1)$-set containing both $s$ and $t$, then
such a set can be extended to a $\tau(B)$-set containing $s$ and $t$
by adding to it a $\tau(B_2)$-set. Hence we may assume that there is
no $\tau(B_1)$-set containing both $s$ and $t$, for otherwise we are
done. By induction, the set $\{s,t\}$ is an $(A)$-pair in $B_1$ and
therefore also an $(A)$-pair in $B$. We will now show that in this
case there is no $\tau(B)$-set containing $s$ and $t$, which would
complete the proof of Part (vii). For the sake of contradiction,
assume that $S$ is a $\tau(B)$-set containing $s$ and $t$.

If $S \cap V(B_1)$ is a transversal in $B_1$, then since there is no
$\tau(B_1)$-set containing both $s$ and $t$ and $\{s,t\} \subseteq S
\cap V(B_1)$, we have that $\tau(B_1) < |S \cap V(B_1)|$. However $|S
\cap (V(B_2) \cup \{x\})| \ge \tau(B_2)$, implying that $|S| = |S
\cap V(B_1)| + |S \cap (V(B_2) \cup \{x\})| > \tau(B_1) + \tau(B_2) =
\tau(B)$, a contradiction. Hence, the set $S \cap V(B_1)$ is not a
transversal in $B'$.

The only edge in $B_1$ that is not intersected by the set $S$ is the
edge $\{u_1,v_1\}$, implying that $S \cap \{u_1,v_1\} = \emptyset$.
Since $|S \cap \{x,u_1,v_1\}| \ge 1$, this implies that $x \in S$.
Further since $|S \cap \{u_1,u_2,v_1,v_2\}| \ge 1$, this in turn
implies that $S \cap \{u_2,v_2\} \ne \emptyset$ and that the set $S
\cap V(B_2)$ is a transversal in $B_2$. Therefore, $|S \cap V(B_2)|
\ge \tau(B_2)$. Since the set $S \cap V(B_1)$ is a transversal in
$B_1 - \{u_1,v_1\}$, by Part~(v) we have that $|S \cap V(B_1)| \ge
\tau(B_1 - \{u_1,v_1\}) = \tau(B_1) - 1$. Hence, $\tau(B_1) +
\tau(B_2) = |S| = |S \cap V(B_1)| + |\{x\}| + |S \cap V(B_2)| \ge
(\tau(B_1) -1) + 1 + \tau(B_2) = \tau(B_1) + \tau(B_2)$. Thus we must
have equality throughout this inequality chain. In particular, we
have $|S \cap V(B_1)| = \tau(B_1) - 1$. But then the set $(S \cap
V(B_1)) \cup \{u_1\}$ is a transversal in $B_1$ of size~$\tau(B_1)$
containing both $s$ and $t$, a contradiction. Therefore if $\{s,t\}$
is an $(A)$-pair in $B_1$, then there is no $\tau(B)$-set containing
$s$ and $t$, which completes the proof of Part~(vii).

\2 \textbf{(viii):} We will prove Part~(viii) by induction on the
order, $n(B)$, of the hypergraph $B$. Let $Y_1=\{s_1,t_1\}$,
$Y_2=\{s_2,t_2\}$ and $Y_3=\{s_3,t_3\}$. If $n(G) = 2$, then $B=H_2$
was created in step (A) in Definition~\ref{create_B}. In this case,
$Y_1 = Y_2 = Y_3$ and the result holds trivially. This establishes
the base case. Let $k \ge 3$ and assume that the result holds for all
$B' \in \cB$ with $n(B') < k$ and let $B \in \cB$ have order~$n(B) =
k$.

Assume that $Y_1$, $Y_2$ and $Y_3$ are not vertex disjoint. Renaming
vertices, we may assume that $s_1=s_2$. If $s_1 \in Y_3$, then we are
done by part~(vi) since there exists a $\tau(B)$-set containing
$s_1$. Hence we may assume that $s_1 \notin Y_3$. However by
Part~(iv) either $\{s_1,s_3\}$ or $\{s_1,t_3\}$ is not an $(A)$-pair.
Renaming vertices in $Y_3$ if necessary, we may assume that
$\{s_1,s_3\}$ is not an $(A)$-pair. We are now done by Part~(vii)
since there exists a $\tau(B)$-set containing $s_1$ and $s_3$. Hence
we may assume that $Y_1$, $Y_2$ and $Y_3$ are vertex disjoint, for
otherwise the desired result follows. Let
$X=\{s_1,t_1,s_2,t_2,s_3,t_3\}$, and so $|V(B)| \ge |X| = 6$.

Suppose that $B$ was created from $B'$ in Step (B) in
Definition~\ref{create_B} and name the vertices as in
Definition~\ref{create_B}. By Part~(i), $\tau(B') = \tau(B)-1$.
Suppose $\{x,y\} \cap X = \emptyset$. Applying the inductive
hypothesis to $B'$, there exists a $\tau(B')$-set, $S'$, intersecting
$Y_1$, $Y_2$ and $Y_3$. But then the set $S' \cup \{x\}$ is a
$\tau(B)$-set intersecting $Y_1$, $Y_2$ and $Y_3$. Hence we may
assume, renaming vertices if necessary, that $s_1=x$. Since $Y_2$ and
$Y_3$ are vertex disjoint sets, the vertex $y$ belongs to at most one
of the sets, implying that there exists a vertex, $w_1$, in $Y_2
\setminus \{y\}$ and a vertex, $w_2$, in $Y_3 \setminus \{y\}$ that
together do not form an $(A)$-pair by Part~(iv). However, by
Part~(vii), this implies that there exists a $\tau(B')$-set, $S'$,
containing $w_1$ and $w_2$. Thus the set $S' \cup \{x\}$ is a
$\tau(B)$-set covering $Y_1$, $Y_2$ and $Y_3$.

If $B$ was created from $B'$ in Step~(C) in
Definition~\ref{create_B}, then the proof that Part~(viii) holds is
analogous to when $B$ was created in Step~(B).

Suppose finally that $B$ was created from $B_1$ and $B_2$ in Step~(D)
and name the vertices as in Definition~\ref{create_B}. By Part~(ii),
$\tau(B) = \tau(B_1) + \tau(B_2)$. For $i = 1,2$, let $X_i = X \cap
V(B_i)$. Then, $|X_1| \ge 3$ or $|X_2| \ge 3$. Renaming $B_1$ and
$B_2$ if necessary, we may assume without loss of generality that
$|X_1| \ge 3$.

If $|X_1| = 6$, then by induction there exists a $\tau(B_1)$-set,
$S_1$, covering all three sets, $Y_1$, $Y_2$ and $Y_3$. Let $S_2$ be
a $\tau(B_2)$-set. Then, $S_1 \cup S_2$ is a $\tau(B)$-set covering
$Y_1$, $Y_2$ and $Y_3$. Hence we may assume that $3 \le |X_1| \le 5$.
Further renaming $Y_1$, $Y_2$ and $Y_3$ if necessary, we may assume
by Part~(iv) that $\{s_1,s_2\} \subset V(B_1)$ and that $\{s_1,s_2\}$
is not an $(A)$-pair. Further since $|X_1| \le 5$, we may assume that
$|Y_3 \cap V(B_1)| \le 1$. By Part~(vii), there exists a
$\tau(B_1)$-set, $S_1$, containing $s_1$ and $s_2$. On the one hand
if $x \in Y_3$, then let $S_2'$ be a $\tau(B_2 - \{u_2,v_2\})$-set.
By part~(v), $|S_2'| = \tau(B_2) - 1$. In this case, the set $S_1
\cup S_2' \cup \{x\}$ is a transversal in $B$ of size~$|S_1| + |S_2'|
+ 1 = \tau(B_1) + (\tau(B_2) - 1) + 1 = \tau(B)$ covering $Y_1$,
$Y_2$ and $Y_3$. On the other hand, if $x \notin Y_3$, then $|Y_3
\cap V(B_2)| \ge 1$ and we may assume, renaming $s_3$ and $t_3$ if
necessary, that $s_3 \in V(B_2)$. By Part~(vi), there exists a
$\tau(B_2)$-set, $S_2$, containing $s_3$. In this case the set $S_1
\cup S_2$ is a $\tau(B)$-set covering $Y_1$, $Y_2$ and $Y_3$, which
completes the proof of Part~(viii).

\2 \textbf{(ix):} We will prove Part~(ix) by induction on the order,
$n(B)$, of the hypergraph $B$. Clearly, Part~(ix) is vacuously true
if $B = H_2$. This establishes the base case. Let $k \ge 3$ and
assume that the result holds for all $B' \in \cB$ with $n(B') < k$
and let $B \in \cB$ have order~$n(B) = k$. We first note that no
$2$-edges intersect in any hypergraph in $\cB$, as none of the steps
(A)-(D) in Definition~\ref{create_B} cause $2$-edges to intersect. In
particular, we note that in Step~(B) the $2$-edge $\{u,v\}$ in $B'$
does not intersect any other $2$-edge in $B'$. We now observe that no
$3$-edge in any $B \in \cB$ can intersect two $2$-edges in $B$, as
again none of the steps (A)-(D) in Definition~\ref{create_B} can
cause this to happen. In particular, we observe that in Step~(C) the
$3$-edge $\{u,v,w\}$ in $B'$ intersects at most one other $2$-edge in
$B'$. Finally we observe that no $4$-edge in $B \in \cB$ can
intersect three $2$-edges in $B$, as again none of the steps (A)-(D)
in Definition~\ref{create_B} can cause this to happen. Therefore,
Part~(ix) follows easily by induction.

\2 \textbf{(x):} Part~(x) follows easily by induction and the
observation that Steps (B)-(D) all increase the degrees of existing
vertices being operated on and introduce new vertices of degree two.

\2 \textbf{(xi):} We will prove Part~(xi) by induction on the order,
$n(B)$, of the hypergraph $B$. Clearly, Part~(xi) is vacuously true
if $B = H_2$. This establishes the base case. Let $k \ge 3$ and
assume that the result holds for all $B' \in \cB$ with $n(B') < k$
and let $B \in \cB$ have order~$n(B) = k$. Let $x \in V(B)$ be chosen
such that $d_B(x)=2$. As observed in the proof of Part~(x), Steps
(B)-(D) all increase the degrees of existing vertices being operated
on and introduce new vertices of degree two. If $x$ is a new vertex
of degree two added when constructing $B$, then by construction the
vertex $x$ belongs to a $2$-edge or a $3$-edge. If $x$ is not a new
vertex added when constructing $B$, then by considering Steps (A)-(D)
and Part~(x) above it is not difficult to see that Part~(xi) holds.
This completes the proof of Lemma~\ref{B_properties}.

\2 \textbf{(xii):} We will prove Part~(xii) by induction on the order, $n(B)$, of the hypergraph $B$.
It is not difficult to see that Part~(xii) holds if the order is at most four.
 Let $k \ge 5$ and assume that the result holds for all $B' \in \cB$ with $n(B') < k$
and let $B \in \cB$ have order~$n(B) = k$.
If $B$ was created using Step~(B) or (C), then clearly Part~(xii) holds. If $B$ was created using
Step~(D), then without loss of generality there is a $2$-edge in $B_1$ different from $\{u_1,v_1\}$ (otherwise there is a $2$-edge in $B_2$ different from
$\{u_2,v_2\}$) and Part~(xii) follows by induction on $B_1$.

\2 \textbf{(xiii):} As $B$ does not contain two $4$-edges intersecting in three vertices we note that Step~(C) was never performed in any step of constructing $B$ (as
no operation removes $4$-edges). As Step~(C) was never performed we note that no operation removes $3$-edges. As all $2$-edges in $B$ are created using Step~(B) (any $2$-edge
created in Step~(A) will be removed again by Step~(B) or Step~(D)) we note
that all $2$-edges in $B$ intersects two overlapping $3$-edges.~\qed

\newcommand{\DELETEthis}{
We will prove (XXX) by induction. Let $e \in E(B)$ and $s \in V(B)
\setminus \Vx{e}$ be arbitrary. If $B=H_2$ was created in step (A) in
Definition~\ref{create_B} then $s$ cannot exist and we are done
vacuously.

Now assume that $B$ was created from $B'$ in Step (B) and name the
vertices as in Definition~\ref{create_B}. First assume that
$e=\{u,v,x\}$ or $e = \{u,v,y\}$. If $s \in \{x,y\}$  then by (v)
there exists a $\tau(B' - \{u,v\})$-set, $S'$, with
$|S'|=\tau(B')-1$. As $S' \cup \{s\}$ is a transversal of $B-e$
containing $s$ and of size $|S'|+1 = \tau(B')-1+1 = \tau(B)-1$, we
may assume that $s \notin \{x,y\}$.  Now by induction there exists a
$\tau(B' - \{u,v\})$-set, $S'$, containing $s$ with
$|S'|=\tau(B')-1$. As $S' \cup (\{x,y\}\setminus \Vx{e})$ is a
transversal of $B-e$ containing $s$ and of size $|S'|+1 = \tau(B)-1$,
we are done in this case.  Now consider the case when $e=\{x,y\}$. In
this case any $\tau(B')$-set containing $s$ (which exits by (vi)) is
a $\tau(B-e)$-set containing $s$. Finally consider the case when $e
\in E(B')$ By induction or Part (v) we can let $S'$ be a
$\tau(B'-e)$-set of size $\tau(B')-1$ and containing $s$ if $s \in
V(B')$.  If $s \notin V(B')$ then note that $S' \cup \{s\}$ is a
$\tau(B-e)$-set of size $\tau(B)-1$ containing $s$.

If $B$ was created from in Step (C)  in Definition~\ref{create_B},
then the proof that (iii) holds is analogous to when $B$ was created
in Step (B).

Now assume that $B$ was created from $B_1$ and $B_2$ in Step (D) and
name the vertices as in Definition~\ref{create_B}. If
$e=\{u_1,v_1,x\}$, then by induction (or Part (v)) there exists a
$\tau(B_1 - \{u_1,v_1\})$-set, $S_1$, containing $s$ if $s \in
V(B_1)$. Let $S_2$ be any $\tau(B_2)$-set containing $s$ if $s \in
V(B_2)$ and note that $S_1 \cup S_2$ is a transversal in $B-e$ of
size $\tau(B_1)+\tau(B_2)-1 = \tau(B) - 1$ containing $s$. If $e =
\{u_2,v_2,x\}$ then analogously we are done. Now consider the case
when $e=\{u_1,v_1,u_2,v_2\}$ and let $S_i$ be a $\tau(B_i -
\{u_i,v_i\})$-set containing $s$ if $ s \in V(B_i)$ for $i=1,2$. Note
that $S_1 \cup S_2 \cup \{x\}$ is a transversal in $B-e$ of size
$\tau(B_1)-1 +\tau(B_2)-1 +1= \tau(B) - 1$ containing $s$. Finally if
$e \in E(B_i)$ for some $i\in \{1,2\}$, then without loss of
generality assume $i=1$ and note that by induction (or Part (v))
there exists a $\tau(B_1 - e)$-set, $S_1$, containing $s$ if $s \in
V(B_1)$. If $s=x$ then let $S^*$ contain any $\tau(B_2 -
\{u_2,v_2\})$-set together with $S_1 \cup \{x\}$, whereas if $s \ne
x$ then let $S^*$ be $S_1$ together with any $\tau(B_2)$-set
containing $s$ if $s \in V(B_2)$. Note that in both cases $S^*$ is a
transversal of $B-e$ containing $s$ and of size $\tau(B)-1$. This
completes the proof of Part (XXX). }


\section{Proof of Main Result}
\label{S:proof}

In this section, we present a proof of our main result, namely
Theorem~\ref{main_thm}. Recall its statement, where $\cH$ denotes the
class of hypergraphs where all edges have size at most four and at
least two and with maximum degree at most three.

\noindent \textbf{Theorem~\ref{main_thm}}. \emph{If $H \in \cH$, then
\[
24 \tau(H) \le 6n(H) + 4e_4(H) + 6e_3(H) + 10e_2(H) + 2b(H) + b^1(H).
\]
Furthermore if $b^1(H)$ is odd, then the above inequality is strict.
}

\noindent \textbf{Proof of Theorem~\ref{main_thm}.} Given any $H' \in
\cH$, let
\[
\phi(H') =  6n(H') + 4e_4(H') + 6e_3(H') + 10 e_2(H') +
2b(H') + b^1(H').
\]
We note that if $b^1(H')$ is odd, then $\phi(H')$ is odd. Hence if
$24 \tau(H) \le \phi(H')$ and $b^1(H')$ is odd, then $24 \tau(H) <
\phi(H')$.

If $e \in E(H')$, we let $\omega_{H'}(e)$, or simply $\omega(e)$ if
$H'$ is clear from the context, denote the contribution of the edge
$e$ to the expression $\phi(H')$; that is,

\[
\omega(e) = \left\{
\begin{array}{cl}
4 & \mbox{if $e$ is a $4$-edge} \\
6 & \mbox{if $e$ is a $3$-edge} \\
10 & \mbox{if $e$ is a $2$-edge}
\end{array}
\right.
\]

We refer to $\omega(e)$ as the \emph{weight} of the edge $e$. Suppose
to the contrary that the theorem is false. Among all counterexamples,
let $H$ be chosen so that $n(H)+m(H)$ is minimum. In particular,
$24 \tau(H) > \phi(H)$. We will often use the following fact.

\begin{description}
\item[Fact 1:] Let $H' \in \cH$ be a hypergraph with
    $n(H')+m(H') < n(H)+m(H)$. Then the following holds. \\
    \hspace*{0.5cm} {\rm (a)} $\phi(H) - \phi(H') < 24(\tau(H) - \tau(H'))$. \\
    \hspace*{0.5cm} {\rm (b)} If $H'=H(X,Y)$, then $\phi(H) -
    \phi(H') < 24|X|$.
\end{description}
\proof (a) Let $H' \in \cH$ satisfy $n(H')+m(H') < n(H)+m(H)$. If
$\phi(H) - \phi(H') \ge 24\tau(H) - 24\tau(H')$, then $24\tau(H') \ge
\phi(H') + (24\tau(H)-\phi(H)) > \phi(H')$, contradicting the
minimality of $H$. Hence, $\phi(H) - \phi(H') < 24\tau(H) -
24\tau(H')$.

(b) Further suppose $H'=H(X,Y)$. If $X'$ is a $\tau(H')$-set, then $X
\cup X'$ is a transversal in $H$, implying that $\phi(H) < 24\tau(H)
\le 24|X| + 24|X'| = 24\tau(H') + 24|X| \le \phi(H') + 24|X|$, or,
equivalently, $\phi(H) - \phi(H') < 24|X|$.~\smallqed

\medskip
In what follows we present a series of claims describing some
structural properties of $H$ which culminate in the implication of
its non-existence.

\begin{claim} \label{trivial}
No edge of $H$ is contained in another edge of $H$.
\end{claim}
\proof Let $e$ and $f$ be two distinct edges of $H$ and suppose to
the contrary that $V(e) \subseteq V(f)$. Let $H' = H - f$. By the
minimality of $H$, we have that $24 \tau(H') \le \phi(H')$. Since
every transversal of $H'$ is a transversal of $H$, we have that
$\tau(H) \le \tau(H')$. Hence, $24\tau(H) \le 24\tau(H') \le \phi(H')
= \phi(H) - \omega(f) \le \phi(H) - 4 < \phi(H)$, a
contradiction.~\smallqed

\begin{claim} \label{claim1}
The following hold in the hypergraph $H$. \\
\indent {\rm (a)} $H$ is connected. \\
\indent {\rm (b)} $b(H)=0$. \\
\indent {\rm (c)} $b^1(H)=0$.
\end{claim}
\proof (a) If $H$ is disconnected, then by the minimality of $H$ we
have that the theorem holds for all components of $H$ and therefore
also for $H$, a contradiction.

(b) If $b(H)>0$, then by Part (a), $H \in \cB$ and by Lemma
\ref{B_properties}(iii) we note that $H$ is not a counter-example to
the theorem, a contradiction.

(c) Suppose to the contrary that $b^1(H)>0$. Let $B \in \cB$ be a
subhypergraph in $H$ and let $e \in E(H)$ be the (unique) edge of
$E(H) \setminus E(B)$ intersecting $B$ in $H$. By
Lemma~\ref{B_properties}(vi) there exists a transversal $S$ of $B$
containing a vertex, $v$, in $e$. Let $H' = H(S,V(B) \setminus S)$.
If a vertex, $v'$, in $V(e) \setminus \{v\}$ contributes one to
$b(H')$, then necessarily $v'$ belongs to a component $B' \in \cB$ of
$H'$ and therefore contributes one to $b^1(H)$. In this case, $v'$
contributes one to $2b(H)+b^1(H)$ and two to $2b(H')+b^1(H')$.
If a vertex, $v'$, in $V(e) \setminus \{v\}$ contributes one to
$b^1(H')$, then $v'$ contributes one to $b^2(H)$. In this case, $v'$
contributes zero to $2b(H)+b^1(H)$ and one to $2b(H')+b^1(H')$.
In both cases, the vertex $v'$ increases $2b(H')+b^1(H')$ by one.
Since $|V(e)| \le 4$, we note that $|V(e) \setminus \{v\}| \le 3$.
Thus since each vertex in $V(e) \setminus \{v\}$ increases
$2b(H')+b^1(H')$ by at most one, and since the deletion of the
subhypergraph $B'$ from $H$ decreases $2b(H')+b^1(H')$ by one, we
have that $2b(H')+b^1(H')$ is at most two larger than $2b(H)+b^1(H)$;
that is,
\[
(2b(H)+b^1(H)) - (2b(H')+b^1(H')) \ge - 2.
\]
Further since $\omega_H(e) \ge 4$, and applying
Lemma~\ref{B_properties}(iii) to $B \in \cB$, we have that

\[
\begin{array}{lcl}
\phi(H)-\phi(H') & = & (6 n(B) + 4 e_4(B) + 6 e_3(B) + 10 e_2(B)) + \omega_H(e)  \\ \2
& &   \hspace*{0.5cm}  + (2b(H)+b^1(H)) - (2b(H')+b^1(H') \\   \2
& \ge & (24|S| - 2) + 4 - 2 \\
& = & 24|S|,
\end{array}
\]

\noindent contradicting Fact~1. Therefore, $b^1(H)=0$.~\smallqed

\begin{claim} \label{claim2}
$b^2(H)=0$.
\end{claim}
\proof Suppose to the contrary that $b^2(H)>0$. Let $B \in \cB$ be
any subhypergraph in $H$ contributing to $b^2(H)$ and let $f_1,f_2
\in E(H)$ be the two edges of $E(H) \setminus E(B)$ intersecting $B$
in $H$. We now show a number of subclaims.

\begin{unnumbered}{Subclaim~3(a)}
$|\Vx{f_i} \cap V(B)|=1$ for $i=1,2$. Further if $\Vx{f_i} \cap V(B)
= \{s_i\}$, then $s_1 \ne s_2$ and $\{s_1,s_2\}$ is an $(A)$-pair in $B$.
\end{unnumbered}
\textbf{Proof of Subclaim~3(a).} Suppose to the contrary that
$|\Vx{f_i} \cap V(B)| \ge 2$ for some $i=1,2$ or that $\Vx{f_i} \cap
V(B) = \{s_i\}$ but $\{s_1,s_2\}$ is not an $(A)$-pair in $B$. We now
choose a $\tau(H)$-set, $S$, as follows. If there exists a vertex $v
\in \Vx{f_1} \cap \Vx{f_2}$, then by Lemma~\ref{B_properties}(vi),
let $S$ be chosen to contain~$v$. If $f_i$ intersects $B$ in at least
two vertices for some $i \in \{1,2\}$, then by
Lemma~\ref{B_properties}(iv) we can find vertices $s_j \in \Vx{f_j}$
such that $\{s_1,s_2\}$ is not an $(A)$-pair in $B$. By
Lemma~\ref{B_properties}(vii), let $S$ be chosen to contain $s_1$ and
$s_2$. Finally if $\Vx{f_i} \cap V(B) = \{s_i\}$ where $s_1 \ne s_2$
but $\{s_1,s_2\}$ is not an $(A)$-pair in $B$, then by
Lemma~\ref{B_properties}(vii) let $S$ be chosen to contain $s_1$ and
$s_2$. In all three cases, we have that the $\tau(B)$-set, $S$,
covers $f_1$ and $f_2$. Let $H' = H(S,V(B)\setminus S)$. A similar
argument as in the proof of Claim~\ref{claim1}(c) shows that each
vertex in $\Vx{f_i} \setminus V(B)$ for $i \in \{1,2\}$, increases
$2b(H')+b^1(H')$ by at most one. Hence since $2b(H)+b^1(H)=0$ and
$|\Vx{f_1} \setminus S| + |\Vx{f_2} \setminus S| \le 6$, we have that
$2b(H')+b^1(H') \le 6$, and so
\[
(2b(H)+b^1(H)) - (2b(H')+b^1(H')) \ge - 6.
\]

Further since $\omega_H(f_i) \ge 4$ for $i \in \{1,2\}$, and applying
Lemma~\ref{B_properties}(iii) to $B \in \cB$, we have that

\[
\begin{array}{lcl}
\phi(H)-\phi(H') & = & (6 n(B) + 4 e_4(B) + 6 e_3(B) + 10 e_2(B)) + \omega_H(f_1) + \omega_H(f_2) \\ \2
& &   \hspace*{0.5cm}  + (2b(H)+b^1(H)) - (2b(H')+b^1(H') \\   \2
& \ge & (24|S| - 2) + 4 + 4 - 6 \\
& = & 24|S|,
\end{array}
\]

\noindent contradicting Fact~1 and proving Subclaim~3(a).~\smallqed

\begin{unnumbered}{Subclaim~3(b)}
$B = H_2$.
\end{unnumbered}
\textbf{Proof of Subclaim~3(b).} By Subclaim~3(a), we may assume relabeling vertices if necessary that $\Vx{f_1} \cap V(B) = \{s_1\}$ and $\Vx{f_2} \cap V(B) = \{s_2\}$ and that $\{s_1,s_2\}$ is an $(A)$-pair in $B$. Suppose to the contrary that $B \ne H_2$. Let $H'$ be obtained from $H$ by removing all edges in $B$ and all vertices $V(B) \setminus \{s_1,s_2\}$ and adding the $2$-edge $\{s_1,s_2\}$. We show that $\tau(H) \le \tau(H') + \tau(B) -1$. Let $S'$ be a $\tau(H')$-set such that $|S' \cap \{s_1,s_2\}|$ is a minimum.  Since $\{s_1,s_2\}$ is an edge in $H'$, we have that $|S' \cap \{s_1,s_2\}| \ge 1$. If $|S' \cap \{s_1,s_2\}|=2$, then by removing $s_2$ from $S'$ and replacing it with an arbitrary  vertex in $\Vx{f_2} \setminus \{s_2\}$ we get a contradiction to the
minimality of $|S' \cap \{s_1,s_2\}|$. Therefore, $|S' \cap \{s_1,s_2\}| = 1$. Renaming vertices if necessary, we may assume that $s_1 \in S'$.  By Lemma~\ref{B_properties}(vi) there exists a transversal, $S_B$, of $B$
containing the vertex $s_1$. Thus, $S' \cup S_B$ is a transversal in $H$ and $S' \cap S_B = \{s_1\}$, and so $\tau(H) \le |S'| + |S_B| - 1 = \tau(H') + \tau(B) -1$, as desired. Equivalently, $\tau(B) - 1 \ge \tau(H) - \tau(H')$.
By Lemma \ref{B_properties}(iii), we therefore have that
\[
\begin{array}{lcl} \1
\phi(H)-\phi(H') & = & \phi(B) - \phi(H_2) \\ \1
& = & 24 \tau(B) - 24 \tau(H_2) \\   \1
& = & 24 (\tau(B) - \tau(H_2)) \\ \1
& = & 24 (\tau(B) - 1) \\ \1
& \ge & 24 (\tau(H) - \tau(H')),
\end{array}
\]

\noindent
where $H_2$ is defined in Definition~\ref{create_B}(A),  contradicting Fact~1 and proving Subclaim~3(b).~\smallqed

\begin{unnumbered}{Subclaim~3(c)}
There is no edge $e \in E(H)$ with $\Vx{e} \subseteq (\Vx{f_1} \cup \Vx{f_2}) \setminus \{s_1,s_2\}$.
\end{unnumbered}
\textbf{Proof of Subclaim~3(c).}  Suppose to the contrary that there is an edge $e \in E(H)$ such that $\Vx{e} \subseteq (\Vx{f_1} \cup \Vx{f_2}) \setminus \{s_1,s_2\}$.   Let $H'$ be obtained from $H$ by deleting the vertices $s_1$ and $s_2$ and the edges $f_1,f_2,\{s_1,s_2\}$; that is, $H' = H(\{s_1,s_2\},\emptyset)$. Let $S'$ be a $\tau(H')$-set. Due to the existence of the edge $e$ we may assume without loss of generality that $f_1$ contains a vertex from $S'$. But then $S' \cup \{s_2\}$ is a transversal of $H$, implying that $\tau(H) \le |S'| + 1 = \tau(H')+1$.
Each vertex in $(\Vx{f_1} \cup \Vx{f_2}) \setminus \{s_1,s_2\}$ increases $2b(H')+b^1(H')$ by at most one. Thus since $2b(H)+b^1(H)=0$ and
$|(\Vx{f_1} \cup \Vx{f_2}) \setminus \{s_1,s_2\}| \le 6$, we have that $2b(H')+b^1(H') \le 6$. Further, $\omega(f_1) \ge 4$, $\omega(f_2) \ge 4$ and $\omega(\{s_1,s_2\}) = 10$. Therefore since the vertices $s_1$ and $s_2$ and the edges $f_1,f_2,\{s_1,s_2\}$ are removed from $H$ when constructing $H'$, we have that
\[
\begin{array}{lcl} \1
\phi(H)-\phi(H') & = & 6|\{s_1,s_2\}| +  \omega(f_1)  + \omega(f_2) + \omega(\{s_1,s_2\}) \\
& & \hspace*{0.5cm} - (2b(H')+b^1(H')) \\ \1
& = & 12 + 4 + 4 + 10 - 6 = 24  \\ \1
& \ge & 24 (\tau(H) - \tau(H')),
\end{array}
\]

\noindent
contradicting Fact~1 and proving Subclaim~3(c).~\smallqed

\begin{unnumbered}{Subclaim~3(d)}
$b(H-f_1-f_2) = 1$.
\end{unnumbered}
\textbf{Proof of Subclaim~3(d).}  Since $B$ is a component of $H -
f_1 - f_2$, we have that $b(H-f_1-f_2) \ge 1$. We show that
$b(H-f_1-f_2) = 1$. Suppose to the contrary that there exists a
component, $R \in \cB$, in $H-f_1-f_2$ which is different from $B$.
Since $b(H) = b^1(H) = 0$, the subhypergraph $R$ contributes to
$b^2(H)$, which by Subclaim~3(b) implies that $R=H_2$. By
Subclaim~3(a) we note that the $2$-edge in $R$ is a subset of
$(\Vx{f_1} \cup \Vx{f_2}) \setminus \{s_1,s_2\}$, a contradiction to
Subclaim~3(c).~\smallqed

\2
We now return to the proof of Claim~\ref{claim2}.  By Subclaim~3(a) and~3(b), we may assume that $B=H_2$, $V(B)=\{s_1,s_2\}$ and $\Vx{f_i} \cap V(B) = \{s_i\}$ for $i=1,2$. Let $X = (\Vx{f_1} \cup \Vx{f_2}) \setminus
\{s_1,s_2\}$ and assume without loss of generality that $|\Vx{f_1}| \le |\Vx{f_2}|$. Clearly, $1 \le |X| \le 6$. We now consider a number of different cases.

First consider the case when $|X|=1$. Assume that $X=\{x\}$, which
implies that $f_1=\{s_1,x\}$ and $f_2=\{s_2,x\}$. Let
$H'=H(\{x\},\emptyset)$. Suppose $d_H(x)=2$. Then, $H' = B$,
$b(H')=1$ and $b^1(H)=0$, implying that $\phi(H)-\phi(H') = 6|\{x\}|
+ \omega(f_1)  + \omega(f_2) - (2b(H')+b^1(H')) = 6 + (2\times 10) -
2  = 24 = 24|X|$, contradicting Fact~1. Hence, $d_H(x) \ge 3$.
Consequently since $\Delta(H) = 3$, we have that $d_H(x)=3$. Let $e$
be the edge of $H$ different from $f_1$ and $f_2$ containing $x$ and
note that $2b(H-e)+b^1(H-e) \le 3$, which implies that $2b(H') +
b^1(H') \le 5$. Therefore, $\phi(H)-\phi(H') = 6|\{x\}| + \omega(e) +
\omega(f_1)  + \omega(f_2) - (2b(H')+b^1(H')) \ge 6+4+ (2\times 10) -
5 > 24 = 24|X|$, contradicting Fact~1. Hence, $|X| \ge 2$.

Suppose $2 \le |X| \le 4$. In this case we let $H'$ be obtained from
$H$ by deleting the vertices $s_1$ and $s_2$ and the edges
$f_1,f_2,\{s_1,s_2\}$ and adding the new edge $f=X$.  By
Subclaim~3(d), $b(H-f_1-f_2) = 1$ and therefore $B$ is the only
component of $H - f_1 - f_2$ in $\cB$. This implies that if $b(H')>0$
or $b^1(H')>0$, then the new edge $f$ belongs to some subhypergraph
$R$ which contributes to $b(H')$ or $b^1(H')$, and this $R$ is the
only subhypergraph that contributes to $2b(H')+b^1(H')$. Therefore,
$2b(H')+b^1(H') \le 2$. We now show that $\tau(H) \le \tau(H')+1$.
Assume that $S'$ is a $\tau(H')$-set and note that some vertex in $X$
belongs to $S'$. Without loss of generality we may assume that there
is a vertex in $S' \cap X$ belonging to $f_1$. This implies that $S'
\cup \{s_2\}$ is a transversal of $H$, and so $\tau(H) \le |S'| + 1 =
\tau(H') +1$. We now consider the following possibilities.

Suppose that $|X|=2$. Suppose that $|\Vx{f_1}| = 2$. As observed
earlier, $2b(H')+b^1(H') \le 2$. Since $|\Vx{f_2}| \le |X|+1 = 3$, we
have that $\omega(f_2) \ge 6$. Thus,

\[
\begin{array}{lcl} \1
\phi(H)-\phi(H') & = & 6|\{s_1,s_2\}| +  \omega(f_1)  + \omega(f_2) + \omega(\{s_1,s_2\}) \\
& & \hspace*{0.5cm} - \omega(f) - (2b(H')+b^1(H')) \\ \1
& \ge & (6 \times 2) + 10 + 6 + 10 - 10 - 2 > 24  \\ \1
& \ge & 24 (\tau(H) - \tau(H')),
\end{array}
\]

\noindent
contradicting Fact~1. Hence, $|\Vx{f_1}| = 3$. Thus, $3 = |\Vx{f_1}|
\le |\Vx{f_2}| \le |X|+1 = 3$, implying that $|\Vx{f_2}| = 3$.
Assume that $X=\{x,y\}$, which implies that $f_1=\{s_1,x,y\}$ and
$f_2=\{s_1,x,y\}$. If $b(H')=b^1(H')=0$, then $\phi(H)-\phi(H') \ge
(2 \times 6) + (2 \times 6) + 10 - 10 = 24 = 24(\tau(H) - \tau(H'))$,
contradicting Fact~1. Hence, $2b(H')+b(H')>0$. This implies that the
new edge $f$ belongs to some subhypergraph $R$ which contributes to
$b(H')$ or $b^1(H')$, and this $R$ is the only subhypergraph that
contributes to $2b(H')+b^1(H')$. Since $\{x,y\}$ is a $2$-edge in
$R$, using Step~(B) in  Definition~\ref{create_B} we can extend $R$
to a subhypergraph  $R' \in \cB$, by adding the vertices
$\{s_1,s_2\}$ and the edges $f_1$, $f_2$ and $\{s_1,s_2\}$ and
deleting the edge $\{x,y\}$. However this implies that $R'$ is a
subhypergraph in $H$ contributing to $b(H)$ or $b^1(H)$, a
contradiction. Hence, $|X| \ge 3$.

Suppose that $|X| = 3$. Then, $\omega(f) = 6$. Suppose that
$|\Vx{f_1}| \le 3$. Then, $\omega(f_1) \ge 6$, while $\omega(f_2) \ge
4$. As observed earlier, $2b(H')+b^1(H') \le 2$. Thus,
$\phi(H)-\phi(H') \ge (2 \times 6) + 6 + 4 + 10 - 6 - 2 = 24 =
24(\tau(H) - \tau(H'))$, contradicting Fact~1. Hence, $|\Vx{f_1}| \ge
4$, implying that $|\Vx{f_1}| = |\Vx{f_2}| = 4$.
If $b(H')=b^1(H')=0$, then $\phi(H)-\phi(H') \ge (2 \times 6) + (2
\times 4) + 10 - 6 = 24 = 24(\tau(H) - \tau(H'))$, contradicting
Fact~1. Hence, $2b(H')+b(H')>0$.
This implies that the new edge $f$ belongs to some subhypergraph $R$
which contributes to $b(H')$ or $b^1(H')$, and this $R$ is the only
subhypergraph that contributes to $2b(H')+b^1(H')$. Since $f$ is a
$3$-edge in $R$, using Step~(C) in  Definition~\ref{create_B} we can
extend $R$ to a subhypergraph  $R' \in \cB$, by adding the vertices
$\{s_1,s_2\}$ and the edges $f_1$, $f_2$ and $\{s_1,s_2\}$ and
deleting the edge $f$, a contradiction.

Hence, $|X| = 4$, and so $\omega(f) = 4$. As observed earlier,
$2b(H')+b^1(H') \le 2$. Thus, $\phi(H)-\phi(H') \ge (2 \times 6) + (2
\times 4) + 10 - 4 - 2 = 24 = 24(\tau(H) - \tau(H'))$, contradicting
Fact~1. This completes the case when $2 \le |X| \le 4$.

It remains for us to consider the case when $5 \le |X| \le 6$. In
this case we note that $|\Vx{f_1} \cap \Vx{f_2}| \le 1$. Further,
$|\Vx{f_1}| \ge 3$, and so neither $f_1$ nor $f_2$ is a $2$-edge. Let
$X'$ be the set of vertices from $X$ which belong to some $2$-edge
in~$H$. We note that by Subclaim~3(c), every $2$-edge in $H$ contains
at most one vertex of $X$.

Suppose that $|X'| \le 3$. Let $f \subseteq X$ be chosen such that
$|V(f)|=4$, $X' \subseteq \Vx{f}$ and if any vertex belongs to
$\Vx{e_1} \cap \Vx{e_2}$, then it also belongs to $f$. In particular,
we note that $\omega(f) = 4$. Let $H'$ be obtained from $H$ by
deleting the vertices $s_1$ and $s_2$ and the edges
$f_1,f_2,\{s_1,s_2\}$ and adding the new edge~$f$. Analogously to the
case when $2 \le |X| \le 4$, we have that $\tau(H) \le \tau(H')+1$.
By Subclaim~3(d), $b(H-f_1-f_2) = 1$ and therefore $B$ is the only
component of $H - f_1 - f_2$ in $\cB$. This implies that if
$2b(H')+b^1(H') \ge 3$, then there must exists a subhypergraph $R \in
\cB$ which does not contain the edge $f$ but contributes to
$2b(H')+b^1(H')$. But then $R$ contributed to $b^2(H)$, which by
Subclaim~3(b) implies that $R=H_2$, a contradiction to the definition
of $X'$. Therefore, $2b(H')+b^1(H') \le 2$. Hence, $\phi(H)-\phi(H')
\ge (2 \times 6) + (2 \times 4) + 10 - 4 - 2 = 24 = 24(\tau(H) -
\tau(H'))$, contradicting Fact~1. Hence, $|X'| \ge 4$.

Let $f \subseteq X'$ be chosen such that $|V(f)|=4$. Let $H''$ be obtained from $H$ by deleting the vertices $s_1$ and $s_2$ and the edges $f_1,f_2,\{s_1,s_2\}$ and adding the new edge $f$. Analogously to the case when $2 \le |X| \le 4$, we have that $\tau(H) \le \tau(H'')+1$. By Subclaim~3(d), $b(H-f_1-f_2) = 1$ and therefore $B$ is the only component of $H - f_1 - f_2$ in $\cB$. For the sake of contradiction suppose that there exists a subhypergraph $R \in \cB$ which contains the edge $f$ and contributes to $2b(H'')+b^1(H'')$. By
Lemma~\ref{B_properties}(ix) and Subclaim~3(c) we note that at most two of the four $2$-edges intersecting $f$ can belong to $R$. As observed earlier, neither $f_1$ nor $f_2$ is a $2$-edge. But this implies that the subhypergraph $R \in \cB$ is intersected by at least two $2$-edges in $H''$ that do not belong to $R$, contradicting the fact that $R$ contributes to $2b(H'')+b^1(H'')$.  Therefore, $2b(H'')+b^1(H'') \le |X \setminus \Vx{f}| \le 2$. Hence, $\phi(H)-\phi(H'') \ge (2 \times 6) + (2 \times 4) + 10 - 4 - 2 = 24 = 24(\tau(H) - \tau(H''))$, contradicting Fact~1. This completes the proof of Claim~\ref{claim2}.~\smallqed

\newpage
\begin{claim} \label{claim3}
No $2$-edges in $H$ intersect.
\end{claim}
\proof Suppose to the contrary that there are two $2$-edges, $e$ and $e'$, that intersect in $H$ and let $x$ be the vertex common to both edges. Let $H'=H(\{x\},\emptyset)$ and let $X = \{x\}$. If $d(x)=2$, then Claim~\ref{claim1} and~\ref{claim2} imply that $b(H') = 0$ and $b^1(H') \le 1$, and so $2b(H')+b^1(H') \le 1$. This implies that, $\phi(H)-\phi(H') \ge 6 + (2 \times 10) - (2b(H')+b^1(H')) > 24 = 24|X|$, contradicting Fact~1.  Therefore, $d(x)=3$, which by Claim~\ref{claim1} and~\ref{claim2} implies that $2b(H')+b^1(H') \le 2$ and therefore that $\phi(H)-\phi(H') \ge 6 + (2 \times 10) + 4 - (2b(H')+b^1(H')) > 24 = 24|X|$, contradicting Fact~1.~\smallqed

\begin{claim} \label{claim4}
If $e=\{x,y\}$ is a $2$-edge in $H$ and $d_H(x)=3$, then $x$ is contained in two distinct $4$-edges.
\end{claim}
\proof Assume that $e=\{x,y\}$ is a $2$-edge in $H$ and $d_H(x)=3$. Let $e$, $e'$ and $e''$ be the three distinct edges in $H$ containing $x$. By Claim~\ref{claim3}, neither $e'$ nor $e''$ is a $2$-edge.
Suppose to the contrary that $e'$ is a $3$-edge.  Let $H'=H(\{x\},\emptyset)$. Then, $\tau(H) \le \tau(H') + 1$.

Suppose that $e''$ is a $4$-edge. If $b(H')>0$, then by Claim~\ref{claim1} and~\ref{claim2} we note that any component $R \in \cB$ in $H$ must intersect $e$, $e'$ and $e''$ and therefore contain $y$. This implies that $b(H') \le 1$. Since $|\Vx{e} \setminus \{x\}| + |\Vx{e'} \setminus \{x\}| + |\Vx{e''} \setminus \{x\}| = 6$, we note that by Claim~\ref{claim1} and~\ref{claim2} either $b(H')=1$ and $b^1(H') \le 1$ or $b(H')=0$ and $b^1(H') \le 3$. Thus, $2b(H')+b^1(H') \le 3$. Furthermore if $2b(H')+b^1(H')=3$, then $b^1(H')$ is odd.
By the minimality of $H$ we have $24\tau(H') \le \phi(H')$ when $2b(H')+b^1(H') \le 2$ and
$24\tau(H') \le \phi(H')-1$ when $2b(H')+b^1(H') = 3$.  On the one hand if $2b(H')+b^1(H')=3$, then

\[
\begin{array}{lcl} \1
24\tau(H) & \le & 24(\tau(H')+1) \\
& \le & (\phi(H')-1) + 24 \\ \1
& = & [\phi(H) - 6|\{x\}| - \omega(e) - \omega(e') - \omega(e'') + 2b(H')+b^1(H') - 1] + 24 \\ \1
& = & [\phi(H) - 6 - 10 - 6 - 4 + 3 - 1] + 24 \\ \1
& = & \phi(H),
\end{array}
\]

\noindent a contradiction. On the other hand if $2b(H')+b^1(H') \le 2$, then $24\tau(H) \le 24(\tau(H')+1) \le \phi(H') + 24 = [\phi(H) - 6 - 10 - 6 - 4 + 2] + 24 = \phi(H)$, once again a contradiction. Hence, $e''$ is not a $4$-edge, implying that $e''$ is a $3$-edge. Since $|\Vx{e} \setminus \{x\}| + |\Vx{e'} \setminus
\{x\}| + |\Vx{e''} \setminus \{x\}| = 5$, we note that by Claim~\ref{claim1} and~\ref{claim2} $2b(H')+b^1(H')
\le 3$. Therefore,

\[
\begin{array}{lcl} \1
\phi(H)-\phi(H') & \ge & 6|\{x\}| + \omega(e) + \omega(e') + \omega(e'') - (2b(H')+b^1(H'))  \\
& \ge & 6 + 10 + (2 \times 6) - 3 \\ \1
& > & 24 \\ \1
& = & 24|\{x\}|,
\end{array}
\]

\noindent a contradiction. This completes the proof of Claim~\ref{claim4}.~\smallqed

\begin{claim} \label{claim5}
If $R \in \cB$ is a subhypergraph in $H$ and $e$ is a $2$-edge in $E(H)\setminus E(R)$, then $\Vx{e} \cap V(R) = \emptyset$.
\end{claim}
\proof Assume that $R \in \cB$ is a subhypergraph in $H$ and $e$ is a $2$-edge in $E(H) \setminus E(R)$. Suppose to the contrary that $\Vx{e} \cap V(R) \ne \emptyset$ and let $x \in \Vx{e} \cap V(R)$. If $d_R(x) = 1$, then by Lemma~\ref{B_properties}(x) we have that $R = H_2$ and so $x$ belongs to a $2$-edge in $R$, a contradiction to Claim~\ref{claim3}. Hence, $d_R(x) \ge 2$. However since $\Delta(H) \le 3$ and the edge $e \notin E(R)$ contains the vertex $x$, we have that $d_R(x) \le 2$. Consequently, $d_R(x) = 2$. By Lemma~\ref{B_properties}(xi), $x$ is therefore contained in a $3$-edge or a $2$-edge in $R$, a contradiction to Claim~\ref{claim4}.~\smallqed

\begin{claim} \label{claim6}
If $B \in \cB$ contributes to $b^3(H)$, then $B=H_2$.
\end{claim}
\proof  Assume that $b^3(H)>0$ and that $B \in \cB$ is a subhypergraph in $H$ that contributes to $b^3(H)$. Suppose to the contrary that $B \ne H_2$. Let $f_1,f_2,f_3 \in E(H) \setminus E(B)$ be the three edges in $H$ that intersect $B$.

Suppose that $|V(f_i) \cap V(B)| \ge 2$ for all $i=1,2,3$. Then by Lemma~\ref{B_properties}(viii) there exists a $\tau(B)$-set, $S$, intersecting $f_1$, $f_2$ and $f_3$. Let $H' = H(S,V(B) \setminus S)$. By Claim~\ref{claim1} and~\ref{claim2} we note that any component $R \in \cB$ in $H'$ must intersect all of $f_1$, $f_2$ and $f_3$, while any subhypergraph in $H'$ that contributes to $b^1(H')$ must intersect at least two of $f_1$, $f_2$ and $f_3$. Since $|(\Vx{f_1} \cup \Vx{f_2} \cup \Vx{f_3}) \setminus V(B)| \le 6$, this implies that $2b(H')+b^1(H') \le 4$. Therefore by Lemma \ref{B_properties}(iii), we have that

\[
\begin{array}{lcl} \1
\phi(H)-\phi(H') & \ge & 6 n(B) + 4 e_4(B) + 6 e_3(B) + 10 e_2 (B)  \\
& & \hspace*{0.5cm} + \omega(f_1) + \omega(f_2) + \omega(f_3) - (2b(H')+b^1(H')) \\
& \ge & (24|S| - 2) + 12 - 4 \\ \1
& > & 24|S|,
\end{array}
\]

\noindent
contradicting Fact~1. Hence we may assume without loss of generality that $|\Vx{f_1} \cap B|=1$.

If there is no $\tau(B)$-set intersecting both $f_2$ and $f_3$, then by Lemma~\ref{B_properties}(vii) we must have $\Vx{f_2} \cap B=\{b_2\}$ and $\Vx{f_3} \cap B=\{b_3\}$ and $\{b_2,b_3\}$ is an $(A)$-pair in $B$. However in this case by Lemma~\ref{B_properties}(iv) there exists a $\tau(B)$-set intersecting both $f_1$ and $f_2$. Hence in both cases there exists a $\tau(B)$-set intersecting two of $f_1,f_2,f_3$ such that the edge not covered intersects $B$ in exactly one vertex. Without loss of generality we may assume that $\Vx{f_1} \cap B = \{b_1\}$ and that $S_B$ is a $\tau(B)$-set intersecting both $f_2$ and $f_3$.

Let $H_1^* = H(V(B),\emptyset)$. If $b^1(H_1^*)>0$, then let $B_1 \in \cB$ be a subhypergraph in $H_1^*$ and let $e_1 \in E(H_1^*)$ be the only edge intersecting $B_1$ in $H_1^*$. In this case let $H_2^*=
H_1^*(V(B_1),\emptyset)$. If $b^1(H_2^*)>0$, then let $B_2 \in \cB$ be a subhypergraph in $H_2^*$ and let $e_2 \in E(H_2^*)$ be the only edge intersecting $B_2$ in $H_2^*$. In this case let $H_3^* = H_2^*(V(B_2),\emptyset)$. Continue the above process until $b^1(H_\ell^*) =0$, for some $\ell \ge 1$. This defines $H_1^*,H_2^*,\ldots,H_\ell^*$ and $B_1,B_2,\ldots, B_{\ell-1}$ and $e_1,e_2,\ldots,e_{\ell-1}$.

We first consider the case when $b(H_\ell^*)=0$. Recall that $S_B$ is a $\tau(B)$-set intersecting both $f_2$ and $f_3$. Let $S'=S_B$. We now construct a hypergraph $H'$ where initially we let $H'=H(S_B,V(B)\setminus(S_B \cup \{b_1\}))$. If $b^1(H') > 0$, let $R \in \cB$ be a subgraph in $H'$ intersected by exactly one edge $e \in E(H') \setminus E(R)$ and do the following. Let $S_R$ be a $\tau(R)$-set intersecting $e$ (which exists by Lemma~\ref{B_properties}(vi)) and add $S_R$ to $S'$ and let $H'$ be $H'(S_R,V(R) \setminus S_R)$.  We continue this process until $b^1(H') = 0$. When the above process stops assume that $b^1(H')>0$ was true $r$ times. Let $S'$ consist of the set $S_B$ and the $r$ $\tau(R)$-sets $S_R$ resulting from constructing $H'$.

We show first that $b(H')=0$. Suppose to the contrary that $b(H')>0$ and let $R^* \in \cB$ be a component in $H'$.
This implies that $R$ must contain the edge $f_1$, for if this were not the case, then such a component would also
be a component in $H_\ell^*$, but $b^1(H_\ell^*)=b(H_\ell^*)=0$. However, $f_1$ is not a $2$-edge by Claim~\ref{claim5},
but it does contain a vertex of degree one in $H'$ (namely $b_1$).
However this is a contradiction to Lemma~\ref{B_properties}(x). Therefore, $b(H')=0$ and $2b(H')+b^1(H') = 0$.

Let $V'$ denote all vertices removed from $H$ to obtain $H'$ and let $E'$ be all edges removed. We note that $H' = H(S',V' \setminus (S' \cup \{b_1\}))$. Furthermore the vertex $b_1$ was not removed from $H$ when we initialized $H'$ for the first time. By applying Lemma~\ref{B_properties}(iii) $r+1$ times, we note that $24|S'| = 6(|V'|+1) + 4 e_4(E') + 6 e_3(E') + 10 e_2(E') + 2(r+1)$. Note that apart from the vertices and edges in subhypergraphs from $\cB$ that were deleted when constructing $H'$ from $H$, a further $r+2$ edges have been removed, namely the two edges $f_2$ and $f_3$ and the $r$ edges from subhypergraphs contributing to $b^1(H')$ when constructing $H'$. Therefore since we have removed in total $r+1$ subhypergraphs in $H$ belonging to $\cB$, we have that

\[
\begin{array}{lcl} \1
\phi(H)-\phi(H') & \ge & 6 |V'| + 4 e_4(E') + 6 e_3(E') + 10 e_2 (E') + 4(r+2)  \\ \1
& = &  6(|V'|+1) + 4 e_4(E') + 6 e_3(E') + 10 e_2 (E') + 2(r + 1) + 2r \\ \1
& \ge & 24|S'| + 2r \\ \1
& \ge & 24|S'|,
\end{array}
\]

\noindent
contradicting Fact~1. Hence, $b(H_\ell^*) > 0$.

Since $B \ne H_2$, we have by Lemma~\ref{B_properties}(x) that $\delta(B) \ge 2$. Since $\Delta(H) = 3$,
each vertex in $V(B)$ is intersected by at most one of the three edges $f_1$, $f_2$ and $f_3$,
implying that $\Vx{f_1} \cap V(B)$, $\Vx{f_2} \cap V(B)$ and $\Vx{f_3} \cap V(B)$  are distinct sets.
By Lemma~\ref{B_properties}(iv) and~\ref{B_properties}(vii), we may assume that there exists a $\tau(B)$-set
intersecting both $f_1$ and $f_2$ and a $\tau(B)$-set intersecting both $f_2$ and $f_3$ (by renaming $f_1$, $f_2$ and $f_3$ if necessary).

Let $R \in \cB$ be a component in $H_\ell^*$. Recall by Claim~\ref{claim1} and~\ref{claim2} that we have $b(H)=b^1(H)=b^2(H)=0$. This implies that there is an edge in
$\{f_1,f_3,e_1,e_2,\ldots,e_{\ell-1}\}$ that intersects $R$. Assume it is $e_{j_1}$. However now there is an edge in
$\{f_1,f_3,e_1,e_2,\ldots,e_{j_1-1}\}$ that intersects $B_{j_1}$.  Assume it is $e_{j_2}$. However now there is an edge in
$\{f_1,f_3,e_1,e_2,\ldots,e_{j_2-1}\}$ that intersects $B_{j_2}$.  Assume it is $e_{j_3}$. Continuing the above process we
note that $j_1 > j_2 > j_3 > \cdots > j_s$ and the edge that intersects $B_{j_s}$ is without loss
of generality $f_1$ (otherwise it is $f_3$). By Lemma~\ref{B_properties}(vi) we can find a minimum transversal in $B_{j_i}$ that
covers the edge $e_{j_{i+1}}$ for each $i \in \{1,2,\ldots\}$. Furthermore we can find a $\tau(B_{j_s})$-set that covers $f_1$ and a
$\tau(R)$-set covering $e_{j_1}$. Taking the union of all of these transversals
we obtain a minimum transversal in each of $R,B_{j_1},B_{j_2},\ldots,B_{j_s}$ that together cover all the edges $e_{j_1},e_{j_2},\ldots,e_{j_s},f_1$.
Similarly by Lemma~\ref{B_properties}(vi) we can readily find a minimum transversal in each hypergraph in
$\{B_1,B_2, \ldots, B_{\ell-1}\} \setminus \{B_{j_1},B_{j_2},\ldots,B_{j_s}\}$ that cover
all edges in $\{e_1,e_2,\ldots,e_{\ell-1}\} \setminus \{e_{j_1},e_{j_2},\ldots,e_{j_s}\}$. Let $S_B$ be a $\tau(B)$-set covering $f_2$ and $f_3$
(if $f_3$ would have intersected $B_{j_s}$ instead of $f_1$, then we would have let $S_B$ cover $f_1$ and $f_2$). Let $S^*$ denote the union of all of these
transversals together with $S_B$. Then, $S^*$ covers every edge in $E^* \cup E^{**}$, where $E^*=\{f_1,f_2,f_3,e_1,e_2,\ldots,e_{\ell-1}\}$
and $E^{**} = E( R \cup B \cup B_1 \cup B_2 \cup \cdots \cup B_{\ell-1})$.

Let $H'$ be obtained from $H$ be removing $S^*$ and all edges incident with $S^*$ and all resulting isolated vertices.
Since $b(H)=b^1(H)=b^2(H)=0$, we note that every component in $H_\ell^*$ which belong to $\cB$ is incident with at least three edges
from $E^*$. Further every edge in $E^*$ intersects at most three such components, implying that $b(H_\ell^*) \le |E^*|= \ell + 2$. Recall that $b^1(H_\ell^*) =0$.
Since $H'$ is obtained from $H_\ell^*$ by removing vertices from the component $R$, we have that $b(H') \le \ell + 1$ and $b^1(H') =0$,
and so $2b(H') + b^1(H') \le 2(\ell + 1)$. Let $V^*=V( R \cup B \cup B_1 \cup B_2 \cup \cdots \cup B_{\ell-1})$ and note that
$H' = H(S^*,V^* \setminus S^*)$. Applying Lemma~\ref{B_properties}(iii) to the $\ell + 1$ hypergraphs $R, B, B_1, B_2, \ldots, B_{\ell-1}$, we therefore have that

\[
\begin{array}{lcl} \1
\phi(H)-\phi(H') & \ge & 6 |V^*| + 4 e_4(E^{**}) + 6 e_3(E^{**}) + 10 e_2 (E^{**}) + 4|E^*| - 2(\ell+1)  \\ \1
& = &  (24|S^*| - 2(\ell + 1)) + 4(\ell + 2) - 2(\ell+1)  \\ \1
& = & 24|S^*| + 4 \\ \1
& > & 24|S^*|,
\end{array}
\]

\noindent
contradicting Fact~1. This completes the proof of Claim~\ref{claim6}.~\smallqed

\begin{claim} \label{claim7}
$e_2(H)=0$, which by Claim~\ref{claim6} also implies that $b^3(H)=0$.
\end{claim}
\proof Suppose to the contrary that $e=\{x,y\}$ is a $2$-edge in $H$. Recall by Claim~\ref{claim1} and~\ref{claim2}
that $b(H)=b^1(H)=b^2(H)=0$. Hence since $H_2 \in \cB$, we have that $d_H(x)=3$ or $d_H(y)=3$ (or both).
Renaming vertices if necessary, we may assume that $d_H(x)=3$. Let $e$, $e_1$ and $e_2$ be the edges in $H$
containing $x$. By Claim~\ref{claim4}, the edges $e'$ and $e''$ are both $4$-edges.
Let $e'=\{x,u_1,v_1,w_1\}$ and $e''=\{x,u_2,v_2,w_2\}$. Let $H'=H(\{x\},\emptyset)$ and let $X = \{x\}$.

If $b(H')>0$, then since $b(H)=b^1(H)=b^2(H)=0$ the component contributing to $b(H')$ must intersect $e$, $e'$ and $e''$ and therefore contains the vertex $y$, contradicting Claim~\ref{claim5}. Therefore, $b(H')=0$.
If $b^1(H') =0$, then

\[
\begin{array}{lcl} \1
\phi(H)-\phi(H') & \ge & 6|X| + \omega(e) + \omega(e') + \omega(e'') - (2b(H')+b^1(H'))  \\
& = & 6 + 10 +  (2 \times 4) - 0 \\ \1
& = & 24 \\ \1
& = & 24|X|,
\end{array}
\]

\noindent
contradicting Fact~1. Hence, $b^1(H') \ge 1$. If $b^1(H')=1$, then $24\tau(H') < \phi(H')$ since
$b^1(H')$ is odd, and so

\[
\begin{array}{lcl} \1
24\tau(H) & \le & 24(\tau(H')+1) \\
& \le & (\phi(H')-1) + 24 \\ \1
& = & [\phi(H) - 6|X| - \omega(e) - \omega(e') - \omega(e'') + 2b(H')+b^1(H') - 1] + 24 \\ \1
& = & [\phi(H) - 6 - 10 - (2 \times 4) + 1 - 1] + 24 \\ \1
& = & \phi(H),
\end{array}
\]

\noindent a contradiction. Hence, $b^1(H') \ge 2$. Let $R \in \cB$
contribute to $b^1(H')$. By Claim~\ref{claim5}, the vertex $y \notin
V(R)$ and therefore $R$ contributes to $b^3(H)$ and is intersected by
both $e'$ and $e''$. By Claim~\ref{claim6}, we have that $R = H_2$.
Let $e_2=\{z,w\}$ denote the edge in $R$, and so $y \notin \{z,w\}$.

Suppose that the edges $e'$ and $e''$ intersect the edge $e_2$ in the same vertex, say $z \in \Vx{e_2} \cap \Vx{e'} \cap \Vx{e''}$.
Now let $H^*$ be obtained from $H$ by deleting the vertices $x$ and $z$ and edges $e$, $e'$, $e''$ and $e_2$ and adding a $2$-edge $\{y,w\}$.
Let $S^*$ be a $\tau(H^*)$-set. In order to cover the $2$-edge $\{y,w\}$, we note that $|S^* \cap \{y,w\}| \ge 1$.
If $y \in S^*$, then $S^* \cup \{z\}$ is a transversal in $H$. If $w \in S^*$, then $S^* \cup \{x\}$ is a transversal in $H$.
In both cases, there exists a transversal in $H$ of size~$|S^*| + 1$, implying that $\tau(H) \le \tau(H^*)+1$.
Furthermore since $|\Vx{e} \setminus \{x\}| + |\Vx{e'} \setminus \{x,z\}| + |\Vx{e''} \setminus \{x,z\}| + |\Vx{e_2} \setminus \{z\}| = 6$ and since we added the edge $\{y,w\}$,
we note that $2b(H^*) + b^1(H^*) \le 6$ (in fact one can show that it is at most $3$). Therefore,

\[
\begin{array}{lcl} \1
\phi(H)-\phi(H^*) & \ge & 6|\{x,z\}| + \omega(e') + \omega(e'') + \omega(e) + \omega(e_2) \\ \1
& & \hspace*{0.5cm} -  \omega(\{y,w\}) - (2b(H^*) + b^1(H^*))  \\ \1
& = & (2 \times 6) +  (2 \times 4) +  (2 \times 10) - 10 - 6 \\ \1
& = & 24 \\ \1
& \ge & 24(\tau(H) - \tau(H^*),
\end{array}
\]

\noindent contradicting Fact~1. Hence, $e'$ and $e''$ do not
intersect $R$ in the same vertex. Renaming vertices in $e'$ and
$e''$, if necessary, we may assume that $e_2 = \{u_1,u_2\}$, where we
recall that $e'=\{x,u_1,v_1,w_1\}$ and $e''=\{x,u_2,v_2,w_2\}$. Since
$b^1(H') \ge 2$ there is also another subhypergraph $R' \in \cB$
which contributes to $b^1(H')$. Analogously to the above arguments
for $R$, we have that $R'$ contributes to $b^3(H)$, $R'$ is
isomorphic to $H_2$ and we may assume that the edge, $e_3$, in $R'$
is $\{v_1,v_2\}$. Since $R$ contributes to $b^3(H)$, there is an edge
$f$ in $E(H) \setminus \{e_2\}$ that intersects $R$ distinct from
$e'$ and $e''$. By Claim~\ref{trivial}, the edge $f$ contains exactly
one of $u_1$ and $u_2$. Therefore exactly one vertex in $\{u_1,u_2\}$
has degree~$2$ in $H$ and the other vertex has degree~$3$ in $H$.
Analogously, there is an edge $f'$ in $E(H) \setminus \{e_3\}$ that
intersects $R'$ distinct from $e'$ and $e''$. Further, exactly one
vertex in $\{v_1,v_2\}$ has degree~$2$ in $H$ and the other vertex
has degree~$3$ in $H$. Without loss of generality we may assume that
$d_H(u_1) = 3$ (and so, $d_H(u_2) = 2$). By Claim~\ref{claim5}, we
note that $f$ and $f'$ are $4$-edges.

Suppose that $d_H(v_2)=3$. In this case, we let $H'' = H(Y,Y')$,
where $Y = \{u_1,v_2\}$ and $Y' = \{u_2,v_1\}$. It is not difficult
to see that $2b(H'') + b^1(H'') \le 6$. Therefore the following holds (even if $f=f'$).

\[
\begin{array}{lcl} \1
\phi(H)-\phi(H'') & \ge & 6|Y| + 6|Y'| + \omega(e') + \omega(e'') + \omega(e_2) + \omega(e_3) \\ \1
& & \hspace*{0.5cm}  + \omega(f) - (2b(H'')+b^1(H''))  \\  \1
& = & (4 \times 6) +  (4 \times 3)  +  (2 \times 10)  - 6 \\  \1
& > & 48  \\ \1
& = & 24|Y|,
\end{array}
\]

\noindent contradicting Fact~1. Therefore, $d_H(v_1)=3$. Let $H^x$ be obtained from $H$ by deleting the vertices $u_2$ and $v_2$ and edges $e'$, $e''$, $e_2$ and $e_3$ and adding the $2$-edge $e^x = \{u_1,v_1\}$. Let $S^x$ be a $\tau(H^x)$-set. In order to cover the $2$-edge $e^x$, we note that $|S^x \cap \{u_1,v_1\}| \ge 1$. If $u_1 \in S^x$, then $S^x \cup \{v_2\}$ is a transversal in $H$. If $v_1 \in S^x$, then $S^x \cup \{u_2\}$ is a transversal in $H$. In both cases, there exists a transversal in $H^x$ of size~$|S^x| + 1$, implying that $\tau(H) \le \tau(H^x)+1$.

If $H^x$ contains a component, $R^x$, that belongs to $\cB$, then since $b(H)=b^1(H)=b^2(H)=0$ the component $R^x$ must intersect at least three of the edges
$e'$, $e''$, $e_2$ and $e_3$ and therefore contains both vertices $u_1$ and $v_1$ (recall that $\{u_1,v_1\}$ is an edge in $H^x$).
Hence, $b(H^x) \le 1$. Suppose $b^1(H^x) \ge 1$. In this case, let $R^x \in \cB$ contribute to $b^1(H^x)$. Since $b(H)=b^1(H)=b^2(H)=0$,
the subhypergraph $R^x$ must intersect at least two of the edges $e'$, $e''$, $e_2$ and $e_3$. In particular, if $w_1 \in V(R^x)$,
then $R^x$ must contains at least one of the vertices $u_1$, $x$ and $w_2$. An analogous argument holds if $w_2 \in V(R^x)$.
Further since $\{u_1,v_1\}$ is an edge of $H^x$, this implies that $b^1(H^x) \le 3$. Moreover, if $b^1(H^x) = 3$, then $b(H^x) = 0$.
Thus if $b(H^x) = 1$, then $b^1(H^x) \le 2$. Hence, $2b(H^x) + b^1(H^x) \le 4$.
Therefore,

\[
\begin{array}{lcl} \1
\phi(H)-\phi(H^x) & \ge & 6|\{u_2,v_2\}| + \omega(e') + \omega(e'') + \omega(e_2) + \omega(e_3) \\ \1
& & \hspace*{0.5cm} -  \omega(\{u_1,v_1\}) - (2b(H^x) + b^1(H^x))  \\ \1
& = & (2 \times 6) +  (2 \times 4) +  (2 \times 10) - 10 - 4 \\ \1
& > & 24 \\ \1
& \ge & 24(\tau(H) - \tau(H^x),
\end{array}
\]

\noindent contradicting Fact~1. This completes the proof of Claim~\ref{claim7}.~\smallqed

\begin{claim} \label{claim8}
There are no $3$-edges $e_1,e_2 \in E(H)$ with $|\Vx{e_1} \cap \Vx{e_2}| = 2$.
\end{claim}
\proof Suppose to the contrary that $e_1,e_2 \in E(H)$ are $3$-edges and $|\Vx{e_1} \cap \Vx{e_2}| = 2$. Let $H'$ be obtained from $H$ by removing $e_1$ and $e_2$ and adding the edge $f=\Vx{e_1} \cap \Vx{e_2}$. Every transversal in $H'$ is also a transversal in $H$, and so $\tau(H) \le \tau(H')$. By Claims~\ref{claim1},~\ref{claim2} and~\ref{claim7} we have that $b(H)=b^1(H)=b^2(H)=b^3(H)=0$. This implies that $b(H-e_1-e_2) = b^1(H-e_1-e_2)=0$, which in turn implies that $2b(H')+b^1(H') \le 2$. Therefore,

\[
\begin{array}{lcl} \1
\phi(H)-\phi(H^x) & \ge & \omega(e_1) + \omega(e_2) - \omega(f) - (2b(H')+b^1(H')) \\ \1
& \ge & 6 + 6  - 10 - 2 \\ \1
& = & 0 \\ \1
& \ge &24(\tau(H) - \tau(H')),
\end{array}
\]

\noindent contradicting Fact~1.~\smallqed

\begin{claim} \label{claim9}
There are no $4$-edges  $e_1,e_2 \in E(H)$ with $|\Vx{e_1} \cap \Vx{e_2}| = 3$.
\end{claim}
\proof This is proved analogously to Claim~\ref{claim8}. Suppose to the contrary that $e_1,e_2 \in E(H)$ are $4$-edges and $|\Vx{e_1} \cap \Vx{e_2}| = 3$. Let $H'$ be obtained from $H$ by removing $e_1$ and $e_2$ and adding the edge $f=\Vx{e_1} \cap \Vx{e_2}$. Then, $\tau(H) \le \tau(H')$ and $2b(H')+b^1(H') \le 2$. Therefore,

\[
\begin{array}{lcl} \1
\phi(H)-\phi(H^x) & \ge & \omega(e_1) + \omega(e_2) - \omega(f) - (2b(H')+b^1(H')) \\ \1
& \ge & 4 + 4  - 6 - 2 \\ \1
& = & 0 \\ \1
& \ge & 24(\tau(H) - \tau(H')),
\end{array}
\]

\noindent contradicting Fact~1.~\smallqed

\begin{claim} \label{claim10}
There is no $3$-edge $e_1$ and $4$-edge $e_2$ in $H$ with $|\Vx{e_1} \cap \Vx{e_2}| = 2$.
\end{claim}
\proof Suppose to the contrary that $e_1=\{u,v,x\}$ is a $3$-edges and  $e_2=\{u,v,s,t\}$ is a $4$-edge with $\Vx{e_1} \cap \Vx{e_2} = \{u,v\}$. Suppose that $d_H(u)=3$ and let $e_u$ be the third edge that contains $u$. If $d_{H(\{u\},\emptyset)}(v)=0$, then let $H'=H(\{u\},\{v\})$. In this case, we note that since $b(H)=b^1(H)=b^2(H)=b^3(H)=0$, we have $2b(H')+b^1(H') \le 1$. Therefore,

\[
\begin{array}{lcl} \1
\phi(H)-\phi(H') & \ge & 2|\{u,v\}| + \omega(e_1) + \omega(e_2) + \omega(e_u) - (2b(H')+b^1(H')) \\ \1
& \ge & 12 + 6 + 4 + 4 - 1 \\ \1
& > & 24 \\ \1
& = & 24|\{u\}|,
\end{array}
\]

\noindent contradicting Fact~1. Hence, $d_{H(\{u\},\emptyset)}(v) > 0$. In this case let $H'$ be obtained from
$H$ by removing $e_1$ and $e_2$ and adding the edge $f=\{u,v\}$. Since $b(H-e_1-e_2)=b^1(H-e_1-e_2)=0$, we note that if $2b(H')+b^1(H')>0$, then the edge $f$ must belong to a subhypergraph $R \in \cB$ which contributes
to $b(H')$ or $b^1(H')$. Since $d_H(u)=3$ and $d_{H(\{u\},\emptyset)}(v)>0$,  we note that $R \ne H_2$.
By Lemma\ref{B_properties}(xii) we note that $R$ contains two $3$-edges overlapping in two vertices or two $4$-edges overlapping in three vertices,
a contradiction against Claim~\ref{claim8} and~\ref{claim9}.
Therefore $2b(H')+b^1(H')=0$ and $\phi(H)-\phi(H') = 6+4-10 = 0$, a contradiction to Fact~1.
Therefore, $d_H(u)=2$. Analogously, $d_H(v)=2$.

Let $H^*=H(\emptyset, \{u\})$. Hence, $H^*$ is obtained from $H$ by deleting the vertex $u$ and the two edges $e_1$ and $e_2$ and adding the $2$-edge $e_1' = \{x,v\}$ and the $3$-edge $e_2' = \{v,s,t\}$. Since every transversal in $H^*$ is a transversal in $H$, we have that $\tau(H) \le \tau(H^*)$. If $b(H^*)=b^1(H^*)=0$, then we have that

\[
\begin{array}{lcl} \1
\phi(H)-\phi(H^*) & \ge & 6|\{u\}| + \omega(e_1) + \omega(e_2) - \omega(e_1') - \omega(e_2') - (2b(H^*)+b^1(H^*)) \\ \1
& \ge & 6 + 6 + 4 - 10 - 6 \\ \1
& = & 0 \\ \1
& \ge & 24(\phi(H)-\phi(H^*)),
\end{array}
\]

\noindent contradicting Fact~1. Hence, $2b(H^*)+b^1(H^*)>0$. Let $R \in \cB$ be a subhypergraph in $H^*$ contributing to $b(H^*)$ or $b^1(H^*)$. Since $b(H-e_1-e_2) = b^1(H-e_1-e_2) = 0$, the edge $e_1'$ or
$e_2'$ must belong to $R$, implying that $v \in V(R)$. However we note that $d_{H^*}(v)=2$ and that $v$ is incident to the $2$-edge $e_1' = \{x,v\}$ and the $3$-edge $e_2' = \{v,s,t\}$.

Suppose $d_R(v) = 1$. Then by Lemma~\ref{B_properties}(x) we have that $R = H_2$. But since the edge $e_2'$ intersects $R$, we have that $R$ contributes to $b^1(H^*)$ and that $e_2'$ is the only edge intersecting $R$. This in turn implies that $d_H(x) = 1$. But then letting $H^x = H(\{u,v,x\}, \emptyset)$, we have that every transversal in $H^x$ can be extended to a transversal in $H$ by adding to it the vertex~$u$, and so $\tau(H) \le \tau(H^x) + 1$. Further, $b(H^x)=b^1(H^x)=0$, and so $\phi(H)-\phi(H^x) = 6|\{u,v,x\}| + \omega(e_1) + \omega(e_2) - (2b(H^x)+b^1(H^x)) = 18 + 6 + 4 > 24 \le 24(\tau(H) - \tau(H^x))$, contradicting Fact~1. Hence, $d_R(v) = 2$.

Since $d_R(v) = 2$, both edges $e_1'$ and $e_2'$ belong to $R$.
By Lemma\ref{B_properties}(xii) we note that $R$ contains two $3$-edges overlapping in two vertices or two $4$-edges overlapping in three vertices.
By Claim~\ref{claim8} and~\ref{claim9} we note that $R$ contains two $3$-edges overlapping in two vertices and $e_2'=\{v,s,t\}$ is one of these $3$-edges.
By Lemma\ref{B_properties}(xiii) and Claim~\ref{claim9} we note that $\{x,s,t\}$ is an edge in $R$ and therefore also in $H^*$ and $H$.
Considering the edges $\{x,s,t\}$ and $\{u,v,s,t\}$ instead of $e_1$ and $e_2$, we have that $d_H(s) = d_H(t) = 2$
(analogously to the arguments showing that $d_H(u)=2$ and $d_H(v)=2$).

Let $F$ be the hypergraph with $V(F) = \{u,v,x,s,t\}$ and with $E(F) = \{e_1,e_2,e_3\}$. We note that $F$ is obtained by using Step~(D)
in Definition~\ref{create_B} on two disjoint copies of $H_2$, and so $F \in \cB$.
On the one hand, if $d_H(x) = 2$, then $H = F$ since recall that, by Claim~\ref{claim1}, $H$ is connected. But this implies that $b(H)=1$.
On the other hand, if $d_H(x)=3$, then $F$ is a component of $H - e'$, where $e'$ denote the edge of $H$ containing $x$ different from
$e_1$ and $e_2$. But this implies that the subhypergraph $F \in \cB$ contributes to $b^1(H)$, and so $b^1(H) \ge 1$. In both cases,
we contradict Claim~\ref{claim1}. This completes the proof of Claim~\ref{claim10}.~\smallqed

\begin{claim} \label{claim11}
No $B \in \cB$ is a subhypergraph of $H$.
\end{claim}
\proof Suppose to the contrary that $R \in \cB$ is a subhypergraph of $H$. By Claim~\ref{claim7}, we have that $e_2(H)=0$, implying that in order to create $R$ in Definition~\ref{create_B} we must have used Step~(D) last. However this implies that a $3$-edge and a $4$-edge overlap in two vertices, a
contradiction to Claim~\ref{claim10}.~\smallqed

\newpage
\begin{claim} \label{claim12}
There are no overlapping edges in $H$.
\end{claim}
\proof  Suppose to the contrary that $e_1,e_2 \in E(H)$ have $|\Vx{e_1} \cap \Vx{e_2}| \ge 2$. By Claims~\ref{claim7},~\ref{claim8},~\ref{claim9} and~\ref{claim10} we note that $e_1$ and $e_2$ are both $4$-edges and $|\Vx{e_1} \cap \Vx{e_2}| = 2$. Let $e_1=\{u,v,x_1,y_1\}$ and $e_2=\{u,v,x_2,y_2\}$.
Suppose that $d_H(u)=2$. Let $H'=H(\emptyset,\{u\})$. Hence, $H'$ is obtained from $H$ by deleting the vertex $u$ and the two edges $e_1$ and $e_2$ and adding the edges $e_1' = \{v,x_1,y_1\}$ and $e_2' = \{v,x_2,y_2\}$. Since every transversal in $H'$ is a transversal in $H$, we have that $\tau(H) \le \tau(H')$.  Every $R \in \cB$ contributing to $2b(H')+b^1(H')$ must contain the vertex~$v$, implying that $2b(H')+b^1(H') \le 2$. Therefore,

\[
\begin{array}{lcl} \1
\phi(H)-\phi(H') & \ge & 6|\{u\}| + \omega(e_1) + \omega(e_2) - \omega(e_1') - \omega(e_2') - (2b(H')+b^1(H')) \\ \1
& \ge & 6 + 4 + 4 - 6 - 6 - 2\\ \1
& = & 0 \\ \1
& \ge & 24(\phi(H)-\phi(H')),
\end{array}
\]

\noindent contradicting Fact~1. Therefore, $d_H(u)=3$. Analogously, $d_H(v)=3$. Let $f_u$ be the edge in $E(H) \setminus \{e_1,e_2\}$ containing $u$ and let $f_v$ be the edge in $E(H) \setminus \{e_1,e_2\}$ containing
$v$. Without loss of generality, we may assume that $|\Vx{f_u}| \ge |\Vx{f_v}|$.
Suppose that $f_u = f_v$. In this case, let $H' = H(\{u\},\{v\})$. By
Claim~\ref{claim11}, no $B \in \cB$ is a subhypergraph of $H$, and so
$b(H') = b^1(H') = 0$. Therefore,

\[
\begin{array}{lcl} \1
\phi(H)-\phi(H') & \ge & 6|\{u,v\}| + \omega(e_1) + \omega(e_2) + \omega(f_u) - (2b(H')+b^1(H')) \\ \1
& \ge & 12 + 4 + 4 + 4 \\ \1
& = & 24 \\ \1
& = & 24|\{u\}|,
\end{array}
\]

\noindent contradicting Fact~1. Hence, $f_u \ne f_v$, implying that
$v \notin V(f_u)$ and $u \notin V(F_v)$. By
Claims~\ref{claim8},~\ref{claim9} and~\ref{claim10}, there is a
vertex $w \in \Vx{f_v} \setminus (\Vx{f_u} \cup \{v\})$. Let $f^* =
(\Vx{f_u} \setminus \{u\}) \cup \{w\}$. Then, $|\Vx{f^*}| =
|\Vx{f_u}| \ge |\Vx{f_v}|$. Let $H^*$ be obtained from $H$ be
deleting the edges $e_1,e_2,f_u,f_v$ and the vertices $u$ and $v$,
but adding the edge $f^*$. Let $S^*$ be a $\tau(H^*)$-set and note
that $|S^* \cap V(f^*)| \ge 1$. If $w \in S^*$, then let $S = S^*
\cup \{u\}$, while if $w \notin S^*$, let $S = S^* \cup \{v\}$. In
both cases, $S$ is a transversal in $H$ and $|S| = |S^*| + 1 =
\tau(H^*) + 1$, implying that $\tau(H) \le \tau(H^*) + 1$. Recalling
that $\omega(f^*) = \omega(f_u)$, we have

\[
\begin{array}{lcl} \1
\phi(H)-\phi(H^*) & \ge & 6|\{u,v\}| + \omega(e_1) + \omega(e_2) + \omega(f_u) + \omega(f_v) \\ \1
& & \hspace*{0.5cm} - \omega(f^*) - 2b(H^*) - b^1(H^*) \\ \1
& \ge & 12 + 4 + 4 + \omega(f_v) - 2b(H^*) - b^1(H^*) \\ \1
& = & 20 + w(f_v) - 2b(H^*) - b^1(H^*).
\end{array}
\]

By Claim~\ref{claim11}, no $B \in \cB$ is a subhypergraph of $H$.
Hence any subgraph $R \in \cB$ contributing to $2b(H^*) + b^1(H^*)$
must contain the edge $f^*$, implying that $2b(H^*) + b^1(H^*) \le
2$.
If $f_v$ is a $3$-edge, then $w(f_v)=6$ and $w(f_v) - 2b(H^*) -
b^1(H^*) \ge 4$. But then $\phi(H)-\phi(H^*) \ge 20 + w(f_v) -
2b(H^*) - b^1(H^*) \ge 24 \ge 24(\tau(H) - \tau(H^*))$, contradicting
Fact~1. Hence, $f_v$ is a $4$-edge, implying that $f^*$ and $f_u$ are
$4$-edges. In particular, $\omega(f_v) = 4$.
Furthermore if $2b(H^*) + b^1(H^*) =0$, then $w(f_v) - 2b(H^*) -
b^1(H^*) \ge 4$ and therefore that $\phi(H)-\phi(H^*) \ge 24(\tau(H)
- \tau(H^*))$, contradicting Fact~1. Hence, $2b(H^*) + b^1(H^*) \ge
1$. Let $R \in \cB$ be a subhypergraph in $H^*$ contributing to
$2b(H^*) + b^1(H^*)$.

By Claim~\ref{claim7}, we have $e_2(H)=0$. Since no $2$-edges are
added when constructing $H^*$, we therefore have that $e_2(H^*)=0$.
This implies that the last step performed in the creation of $R$ in
Definition~\ref{create_B} is Step~(D). This in turn implies that in
$H^*$ there is a $4$-edge intersected by two $3$-edges. Moreover,
such a $4$-edge intersects each of these $3$-edges in two vertices.
By Claim~\ref{claim10}, this $4$-edge must therefore be the new edge
$f^*$ added when constructing $H^*$. By Step (D) in
Definition~\ref{create_B} we furthermore note that the two $3$-edges
that intersect the $4$-edge $f^*$ intersect it in disjoint sets.
Hence there is a $3$-edge, not containing the vertex $w \in V(F^*)$,
that intersects $f^*$ in two vertices. But this implies it
intersected $f_u$ in two vertices, a contradiction to
Claim~\ref{claim10}. This completes the proof of
Claim~\ref{claim12}.~\smallqed

\begin{claim} \label{claim13}
$\delta(H) \ge 2$.
\end{claim}
\proof Suppose to the contrary that a vertex $x \in V(H)$ has
$d_H(x)=1$. Let $e$ be the edge containing~$x$ and let $e' = V(e_x)
\setminus \{x\}$. Let $H' = H(\emptyset,\{x\})$ and note that
$\tau(H)=\tau(H')$. By Claim~\ref{claim11}, no $B \in \cB$ is a
subhypergraph of $H$, implying that $2b(H') + b^1(H') \le 2$. If $e$
is a $3$-edge, then $\phi(H)-\phi(H^*) \ge 6|\{x\}| + \omega(e) -
\omega(e') - 2b(H') - b^1(H') \ge 6 + 6 - 10 - 2 \ge 0 \ge 24(\tau(H)
- \tau(H'))$, contradicting Fact~1. Hence, $e$ is a $4$-edge. But
then $\phi(H)-\phi(H^*) \ge 6|\{x\}| + \omega(e) - \omega(e') -
2b(H') - b^1(H') \ge 6 + 4 - 6 - 2 > 0 \ge 24(\tau(H) - \tau(H'))$,
once again contradicting Fact~1.~\smallqed

\begin{claim} \label{claim14}
Every vertex of degree~$2$ in $H$ is incident with two $3$-edges.
\end{claim}
\proof Assume that $d_H(x)=2$. By Claim~\ref{claim7}, we have
$e_2(H)=0$. Suppose to the contrary that $x$ is incident with at
least one $4$-edge, $e$. Let $f$ denote the remaining edge that
contains~$x$. Let $H'=H(\emptyset,\{x\})$ and note that $\tau(H) \le
\tau(H')$. We first show that $2b(H')+b^1(H') =0$. If this is not the
case, let $R \in \cB$ be a subhypergraph of $H'$ that contributes to
$b(H')$ or $b^1(H')$. Since all hypergraphs in $\cB \setminus H_2$
have overlapping edges while there are no overlapping edges in $H'$,
by Claim~\ref{claim12}, we must have that $R = H_2$. By
Claim~\ref{claim11}, no $B \in \cB$ is a subhypergraph of $H$,
implying that $f$ is a $3$-edge and $R$ necessarily contains the
$2$-edge $V(f) \setminus \{x\}$. However since $\delta(H) \ge 2$ by
Claim~2, both vertices in $R$ are incident with at least one edge in
$E(H) \setminus \{e'\}$. Further since there are no overlapping edges
in $H$, these edges are distinct. But this implies that there are
least two edges in $E(H') \setminus E(R)$ intersecting $V(R)$, and so
$R$ does not contribute to $b(H')$ or $b^1(H')$, a contradiction.
Therefore, $2b(H')+b^1(H') = 0$. Letting $e' = V(e) \setminus \{x\}$
and $f' = V(f) \setminus \{x\}$, we note that $\omega(e') - \omega(e)
=2$ and $\omega(f') - \omega(f) \le 4$. Therefore,

\[
\begin{array}{lcl} \1
\phi(H)-\phi(H') & \ge & 6|\{x\}| + \omega(e) + \omega(f) - \omega(e') - \omega(f') - (2b(H')+b^1(H')) \\ \1
& \ge & 6 - (\omega(e') - \omega(e)) - (\omega(f') - \omega(f)) - 0 \\ \1
& \ge & 6 - 2 - 4 \\ \1
& = & 0 \\ \1
& \ge & 24(\tau(H) - \tau(H'),
\end{array}
\]

\noindent contradicting Fact~1. This completes the proof of
Claim~\ref{claim14}.~\smallqed

\begin{claim} \label{claim15}
Every vertex of degree~$3$ in $H$ is incident with a $3$-edge and a
$4$-edge.
\end{claim}
\proof Assume that $d_H(x)=3$ and suppose to the contrary that $x$ is
contained in only $3$-edges or only $4$-edges. Suppose first that $x$
is contained in only $3$-edges and let $H'=H(\{x\},\emptyset)$. By
Claim~\ref{claim11}, we have that $b(H')=b^1(H')=0$. Therefore since
each of the three edges incident with $x$ has weight~$6$, we have
that $\phi(H)-\phi(H') = 6 + (3 \times 6) = 24 = 24|\{x\}|$,
contradicting Fact~1. Hence, $x$ is contained in only $4$-edges.

We now let $H^*=H(\emptyset,\{x\})$ and note that $\tau(H) \le
\tau(H')$. Since there are no $2$-edges in $H^*$ and no overlapping
edges in $H^*$ by Claim~\ref{claim12}, we note that
$b(H^*)=b^1(H^*)=0$. Therefore since each of the three deleted edges
has weight~$4$ and each of the three added edges has weight~$6$, we
have that $\phi(H)-\phi(H') = 6 - (3 \times 2) = 0 \ge 24(\tau(H) -
\tau(H')$, contradicting Fact~1. This completes the proof of
Claim~\ref{claim15}.~\smallqed

\begin{claim} \label{claim16}
Every $3$-edge in $H$ contains a vertex of degree~$3$.
\end{claim}
\proof Assume that $e=\{u_1,u_2,u_3\} \in E(H)$ and suppose to the
contrary that $d_H(u_1)=d_H(u_2)=d_H(u_3)=2$. For $i = 1,2,3$, let
$e_i$ be the edge in $E(H) \setminus \{e\}$ containing $u_i$. By
Claim~\ref{claim12}, $e_1$, $e_2$ and $e_3$ are distinct edges and by
Claim~\ref{claim14} they are all $3$-edges.

Suppose first that $|\Vx{e_i} \cup \Vx{e_j}| \le 5$ for all $1 \le i
< j \le 3$. In this case, by Claim~\ref{claim12}, every pair of edges
in $\{e_1,e_2,e_3\}$ intersect in exactly one vertex. So let
$\Vx{e_i} \cap \Vx{e_j} = \{v_{i,j}\}$ for $1 \le i < j \le 3$. If
$v_{1,2}$, $v_{1,3}$ and $v_{2,3}$ are not distinct vertices, then we
must have $v_{1,2}=v_{1,3}=v_{2,3}$, which contradicts
Claim~\ref{claim15}. Hence, $v_{1,2}$, $v_{1,3}$ and $v_{2,3}$ are
distinct vertices. Hence, $e_1 = \{u_1,v_{12},v_{13}\}$, $e_2 =
\{u_2,v_{12},v_{23}\}$, and $e_3 = \{u_3,v_{13},v_{23}\}$. Let $H'$
be obtained from $H$ by deleting the edges $e,e_1,e_2,e_2$ and
vertices $u_1,u_2,u_3$ and adding the edge $f =
\{v_{1,2},v_{1,3},v_{2,3}\}$. We will first show that $\tau(H) \le
\tau(H')+1$. Let $S'$ be a $\tau(H')$-set. Since $S'$ intersects the
edge $f$, we may assume, renaming vertices if necessary, that
$v_{1,2} \in S'$. But then $S' \cup \{u_3\}$ is a transversal of $H$,
and so $\tau(H) \le |S'| + 1 = \tau(H')+1$, as desired. Clearly
$2b(H')+b^1(H') \le 2$ as any subgraph contributing to
$2b(H')+b^1(H')$ must contain the added edge $f$ since by
Claim~\ref{claim11} no subhypergraph of $H$ belongs to~$\cB$.
Therefore,

\[
\begin{array}{lcl} \1
\phi(H)-\phi(H') & \ge & 6|\{u_1,u_2,u_3\}| + \omega(e) + \omega(e_1) + \omega(e_2) +  \omega(e_3) \\ \1
& & \hspace*{0.5cm} - \omega(f) - (2b(H') + b^1(H')) \\ \1
& \ge & (3 \times 6) + (4 \times 6) - 6 - 2 \\ \1
& > & 24 \\ \1
& \ge & 24(\tau(H) - \tau(H'),
\end{array}
\]

\noindent contradicting Fact~1. We may therefore assume, renaming
vertices if necessary, that $|\Vx{e_1} \cup \Vx{e_2}| = 6$; that is,
the $3$-edges $e_1$ and $e_2$ do not intersect. Let $f_{1,2} =
(\Vx{e_1} \cup \Vx{e_2}) \setminus \{u_1,u_2\}$ and let $f_3 =
\Vx{e_3} \setminus \{u_3\}$. Let $H^*$ be obtained from $H$ by
deleting the edges $e,e_1,e_2,e_3$ and vertices $u_1,u_2,u_3$ and
adding the edges $f_{1,2}$ and $f_3$. We will first show that
$\tau(H) \le \tau(H^*)+1$. Let $S^*$ be a $\tau(H^*)$-set. Since
$S^*$ intersects the edge $f_{12}$, we may assume, renaming vertices
if necessary, that $S^* \cap \Vx{e_1} \ne \emptyset$. But then $S^*
\cup \{u_2\}$ is a transversal of $H$, and so $\tau(H) \le |S^*| + 1
= \tau(H^*)+1$, as desired. Clearly, $2b(H^*)+b^1(H^*) \le 4$ as any
subhypergraph contributing to $2b(H^*)+b^1(H^*)$ must contain the
edge $f_{1,2}$ or the edge $f_3$. Therefore, since $f_{1,2}$ is a
$4$-edge and $f_3$ a $2$-edge, we have that

\[
\begin{array}{lcl} \1
\phi(H)-\phi(H^*) & \ge & 6|\{u_1,u_2,u_3\}| + \omega(e) + \omega(e_1) + \omega(e_2) +  \omega(e_3) \\ \1
& & \hspace*{0.5cm} - \omega(f_{1,2}) - \omega(f_3) - (2b(H^*) + b^1(H^*)) \\ \1
& \ge & (3 \times 6) + (4 \times 6) - 4 - 10 - 4 \\ \1
& = & 24 \\ \1
& \ge & 24(\tau(H) - \tau(H^*),
\end{array}
\]

\noindent contradicting Fact~1. This completes the proof of
Claim~\ref{claim16}.~\smallqed

\begin{claim} \label{claim17}
Every $3$-edge in $H$ contains at least two vertices of degree~$3$.
\end{claim}
\proof Assume that $e=\{u_1,u_2,u_3\} \in E(H)$ and suppose to the
contrary that $d_H(u_2)=d_H(u_3)=2$. By Claim~\ref{claim16} we have
$d_H(u_1)=3$. Let $e_1'$ and $e_2'$ be the two edges in $E(H)
\setminus \{e\}$ containing $u_1$. For $i = 2,3$, let $e_i$ be the
edge in $E(H) \setminus \{e\}$ containing $u_i$ and let $f_i = V(e_i)
\setminus \{u_i\}$. By Claim~\ref{claim14}, the edges $e_2$ and $e_3$
are both $3$-edges, and so $f_2$ and $f_3$ are both $2$-edges. Let
$H' = H(\{u_1\},\{u_2,u_3\})$. Let $V(f_2) = \{v_2,w_2\}$.

We will first show that $b(H')=0$. If this is not the case, then let
$R \in \cB$ be a component in $H'$. By Claim~\ref{claim12}, we have
that $R = H_2$, which by Claim~\ref{claim7}, implies that $f_2$ or
$f_3$ is the edge in $R$. Renaming vertices if necessary, we may
assume that $E(R) = \{f_2\}$, and so $V(R) = \{v_2,w_2\}$. Since
there is no edge in $E(H') \setminus \{f_2\}$ that intersects $V(R)$,
we note that the edges $f_2$ and $f_3$ do not intersect. By
Claim~\ref{claim14}, each vertex in $V(R)$ is either incident to
three edges in $H$ or two $3$-edges in $H$. Suppose both $v_2$ and
$w_2$ are incident to two $3$-edges in $H$. This implies that there
are two distinct $3$-edges that contain (exactly) one of $v_2$ and
$w_2$ and these two $3$-edges are different from the edge $e_2$ (and
from the edge $e$). Since the vertex $u_1$, which has degree~$3$ in
$H$, cannot be incident to three $3$-edges by Claim~\ref{claim15}, at
least one of these $3$-edges that contain $v_2$ or $w_2$ is different
from both $e_1'$ and $e_2'$. This $3$-edge belongs to $E(H')
\setminus \{f_2\}$ and intersects $V(R)$, a contradiction to the fact
that $R$ is a component in $H'$. Hence at least one of $v_2$ and
$w_2$ is incident to three edges in $H$ and the other to at least two
edges in $H$. But once again this implies that there exists an edge
that contain $v_2$ or $w_2$ and is different from the deleted edges
$e, e_1', e_1'',e_2,e_3$ and the edge $f_2$, a contradiction again to
the fact that $R$ is a component in $H'$. Therefore, $b(H')=0$.

If $b^1(H')=0$, then

\[
\begin{array}{lcl} \1
\phi(H)-\phi(H') & \ge & 6|\{u_1,u_2,u_3\}| + \omega(e) + \omega(e_1') + \omega(e_1'') + \omega(e_2) +  \omega(e_3) \\ \1
& & \hspace*{0.5cm} - \omega(f_2) - \omega(f_3) - (2b(H') + b^1(H')) \\ \1
& \ge & (3 \times 6) + (3 \times 6)+ (2 \times 4) - (2 \times 10) - 0 \\ \1
& = & 24 \\ \1
& = & 24|\{u_1\}|,
\end{array}
\]

\noindent contradicting Fact~1. Hence, $b^1(H') \ge 1$. Let $R \in
\cB$ be a subhypergraph in $H'$ contributing to $b^1(H')$. By
Claim~\ref{claim12}, there are no overlapping edges in $H$ and
therefore in $H'$, implying that $R = H_2$. This in turn implies by
Claim~\ref{claim7} that $f_2$ or $f_3$ is the edge in $R$. Renaming
vertices if necessary, we may assume that $E(R) = \{f_2\}$. Let $e'$
be the edge in $E(H') \setminus \{f_2\}$ that intersects $R$. Since
there are no overlapping edges in $H$, we note that $|\Vx{f_2} \cap
\Vx{e'}| = 1$. Renaming vertices in $f_2$ if necessary, we may assume
that $\Vx{f_2} \cap \Vx{e'} = \{v_2\}$.

We now consider the hypergraph $H^* = H'(\{v_2\},\{w_2\})$ obtained
from $H'$ by deleting the vertices $v_2$ and $w_2$ and deleting the
edge $e'$. We note that $H^* = H(\{\{u_1,v_2\},\{u_2,u_3,w_2\})$. By
Claims~\ref{claim7},~\ref{claim11} and~\ref{claim12} the only
possibly subhypergraph in $H^*$ in $\cB$ is the hypergraph isomorphic
to $H_2$ that consists of the $2$-edge $f_3$, implying that
$2b(H^*)+b^1(H^*) \le 2$. Therefore,

\[
\begin{array}{lcl} \1
\phi(H)-\phi(H^*) & \ge & 6|\{u_1,u_2,u_3,v_2,w_2\}| + \omega(e) + \omega(e_1') + \omega(e_1'') +
\omega(e_2) +  \omega(e_3) \\ \1
& & \hspace*{0.5cm} + \omega(e') - \omega(f_3) - (2b(H^*) + b^1(H^*)) \\ \1
& \ge & (5 \times 6) + (3 \times 6)+ (3 \times 4) - 10 - 2 \\ \1
& = & 48 \\ \1
& \ge & 24|\{u_1,v_2\}|,
\end{array}
\]

\noindent contradicting Fact~1. This completes the proof of
Claim~\ref{claim17}.~\smallqed

\begin{claim} \label{claim18}
No vertex is contained in two $3$-edges and one $4$-edge, such that
one of the $3$-edges contains a degree-$2$ vertex.
\end{claim}
\proof Assume that $e_1=\{x,u_1,v_1\}$, $e_2=\{x,u_2,v_2\}$ and
$e_3=\{x,u_3,v_3,w_3\}$ are edges in $H$ and suppose to the contrary
that $d_H(u_1)=2$. By Claim~\ref{claim14}, $u_1$ is incident with two
$3$-edges, say $e_1$ and $f_1 = \{u_1,x_1,y_1\}$. Let
$H'=H(\{x\},\{u_1\})$. If $2b(H')+b^1(H')>0$, then let $R \in \cB$ be
a subhypergraph in $H'$ contributing to $2b(H')+b^1(H')$. By
Claim~\ref{claim12}, there are no overlapping edges in $H$ and
therefore in $H'$, implying that $R = H_2$. This in turn implies by
Claim~\ref{claim7} that the edge in $E(R)$ is $g = \{x_1,y_1\}$. By
supposition, $d_H(u_1) = 2$. Hence by Claim~\ref{claim17} the two
vertices, namely $x_1$ and $y_1$, in the $3$-edge $f_1$ both have
degree~$3$ in $H$. Since there are no overlapping edges in $H$, there
are therefore four distinct edges in $H$ excluding the edge $f_1$
that intersect $V(R)$. Further we note that the vertex $u_1$ is the
only vertex common to both edges $e_1$ and $f_1$, implying that the
edge $e_1$ does not intersect $V(R)$. Hence removing the three edges
$e_1$, $e_2$ and $e_3$ from $H$ can remove at most two edges
intersecting $V(R)$, implying that at least two edges in $H$ that
intersect $V(R)$ remain in $H'$. But then $R$ does not contribute to
$2b(H')+b^1(H')$, a contradiction. Therefore, $2b(H')+b^1(H')=0$.
This implies that

\[
\begin{array}{lcl} \1
\phi(H)-\phi(H') & \ge & 6|\{x,u_1\}| + \omega(e_1) + \omega(e_2) + \omega(e_3) + \omega(f_1) \\ \1
& & \hspace*{0.5cm} - \omega(g) - (2b(H') + b^1(H')) \\ \1
& \ge & (2 \times 6) + (3 \times 6) + 4 - 10 - 0 \\ \1
& = & 24 \\ \1
& = & 24|\{x\}|,
\end{array}
\]

\noindent contradicting Fact~1.~\smallqed

\begin{claim} \label{claim19}
$H$ is $3$-regular.
\end{claim}
\proof Suppose to the contrary that $\delta(H) = 2$. Let $x$ be a
vertex of degree~$2$ in $H$. By Claim~\ref{claim14}, $x$ is incident
with two $3$-edges in $H$, say $e_1=\{x,u_1,v_1\}$ and
$e_2=\{x,u_2,v_2\}$. By Claims~\ref{claim15},~\ref{claim17}
and~\ref{claim18} each vertex in $\{u_1,v_1,u_2,v_2\}$ is contained
in one $3$-edge and two $4$-edges. Let $f_1$ and $f_2$ be the two
$4$-edges containing $u_1$ and let $h_1$ and $h_2$ be the two
$4$-edges containing $v_1$. For $i \in \{1,2\}$, let $h_i' = V(h_i)
\setminus \{v_1\}$. Let $e_2' = \{u_2,v_2\}$. We note that $h_1'$ and
$h_2'$ are both $3$-edges. We now consider the hypergraph
$H'=H(\{u_1\},\{x,v_1\})$.

We will first show that $2b(H')+b^1(H')=0$. If this is not the case,
then let $R \in \cB$ be a subhypergraph in $H'$ contributing to
$2b(H')+b^1(H')$. By Claim~\ref{claim12}, there are no overlapping
edges in $H$ and therefore in $H'$, implying that $R = H_2$. This in
turn implies by Claim~\ref{claim7} that $V(R) = \{u_2,v_2\}$ as
$e_2'$ is the only $2$-edge in $H'$. By Claim~\ref{claim17}, both
vertices $u_2$ and $v_2$ have degree~$3$ in $H$. Since removing all
edges containing $u_1$ can remove at most two edges intersecting
$V(R)$ in $H$, at least two edges in $H$ that intersect $V(R)$ remain
in $H'$. But then $R$ does not contribute to $2b(H')+b^1(H')$, a
contradiction. Therefore, $2b(H')+b^1(H')=0$. This implies that

\[
\begin{array}{lcl} \1
\phi(H)-\phi(H') & \ge & 6|\{x,u_1,v_1\}| + \omega(e_1) + \omega(e_2) + \omega(f_1) + \omega(f_2) + \omega(h_1) + \omega(h_2) \\ \1
& & \hspace*{0.5cm} - \omega(h_1') - \omega(h_1') - \omega(e_2') - (2b(H') + b^1(H')) \\ \1
& \ge & (3 \times 6) + (2 \times 6) + (4 \times 4) - (2 \times 6) - 10 - 0 \\ \1
& = & 24 \\ \1
& = & 24|\{u_1\}|,
\end{array}
\]

\noindent contradicting Fact~1.~\smallqed

\begin{claim} \label{claim20}
All vertices are contained in two $3$-edges and one $4$-edge.
\end{claim}
\proof Suppose to the contrary that there is a vertex $x$ in $H$ that
is not adjacent with two $3$-edges and one $4$-edge. By
Claim~\ref{claim20}, $d_H(x) = 3$. By Claim~\ref{claim15}, the vertex
$x$ is incident with a $3$-edge and a $4$-edge. By our supposition,
the remaining edge incident with $x$ is a $4$-edge. Let
$e_1=\{x,u_1,v_1\}$, $e_2=\{x,u_2,v_2,w_2\}$ and
$e_3=\{x,u_3,v_3,w_3\}$ be the three edges incident with $x$. For $i
\in \{1,2,\}$, let $e_i' = V(e_i) \setminus \{x\}$.

By Claim~\ref{claim15} and~\ref{claim20}, we have that $d_H(u_1) = 3$
and $u_1$ is incident with either two $3$-edges and one $4$-edge or
with one $3$-edge and two $4$-edges. Suppose that $u_1$ is incident
with two $3$-edges, say $e_1$ and $f_1$. In this case, let $f_2$ be
the $4$-edge that contains~$u_1$. Let $H'=H(\{u_1\},\{x\})$. Since
$e_2(H')=0$ and there are no overlapping edges in $H'$, we note that
$b(H')=b^1(H')=0$. Therefore,

\[
\begin{array}{lcl} \1
\phi(H)-\phi(H') & \ge & 6|\{x,u_1\}| + \omega(e_1) + \omega(e_2) + \omega(e_3) + \omega(f_1) + \omega(f_2)  \\ \1
& & \hspace*{0.5cm} - \omega(e_2') - \omega(e_3') - (2b(H') + b^1(H')) \\ \1
& \ge & (2 \times 6) + (2 \times 6) + (3 \times 4) - (2 \times 6) - 0 \\ \1
& = & 24 \\ \1
& = & 24|\{u_1\}|,
\end{array}
\]

\noindent contradicting Fact~1. Hence, $u_1$ is incident with one
$3$-edge and two $4$-edges. Analogously, $v_1$ is incident with one
$3$-edge and two $4$-edges. Let $h_1$ and $h_2$ be the two $4$-edges
containing $u_1$ and let $g_1$ and $g_2$ be the two $4$-edges
containing $v_1$. For $i \in \{1,2\}$, let $h_i' = V(h_i) \setminus
\{u_1\}$ and let $g_i' = V(g_i) \setminus \{v_1\}$. We now consider
the hypergraph $H^*=H(\{x\},\{u_1,v_1\})$ and note that $e_2(H^*)=0$.
Further since there are no overlapping edges in $H^*$, we note that
$b(H^*)=b^1(H^*)=0$. Therefore,

\[
\begin{array}{lcl} \1
\phi(H)-\phi(H^*) & \ge & 6|\{x,u_1,v_1\}| + \omega(e_1) + \omega(e_2) + \omega(e_3) + \omega(h_1) + \omega(h_2) + \omega(g_1)  \\ \1
& & \hspace*{0.5cm} + \omega(g_2) - \omega(h_1') - \omega(h_2') - \omega(g_1') - \omega(g_1')
- (2b(H^*) + b^1(H^*)) \\ \1
& \ge & (3 \times 6) + 6 + (6 \times 4) - (4 \times 6) - 0 \\ \1 & =
& 24 \\ \1
& = & 24|\{x\}|,
\end{array}
\]

\noindent contradicting Fact~1.~\smallqed

\medskip
We now return to the proof of Theorem~\ref{main_thm} to obtain a
final contradiction implying the non-existence of our counterexample,
$H$, to the theorem.
Let $e=\{u_1,u_2,u_3\}$ be an arbitrary $3$-edge in $H$. By
Claim~\ref{claim20}, each vertex of $H$ is contained in two $3$-edges
and one $4$-edge. For $i \in \{1,2,3\}$, let $e_i$ be the $3$-edge
and $f_i$ the $4$-edge in $E(H) \setminus \{e\}$ that contains the
vertex $u_i$. By Claim~\ref{claim12}, the  edges $e_1$, $e_2$ and
$e_3$ are all distinct. For $i \in \{1,2,3\}$, let $f_i' = V(f_i)
\setminus \{u_i\}$ and note that $f_i'$ is a $3$-edge.

Suppose that $\Vx{e_i} \cap \Vx{e_j} \ne \emptyset$ for all $1 \le i
< j \le 3$. Let $\Vx{e_i} \cap \Vx{e_j} = \{v_{i,j}\}$ for $1 \le i <
j \le 3$. If $v_{1,2}$, $v_{1,3}$ and $v_{2,3}$ are not distinct
vertices, then we must have $v_{1,2}=v_{1,3}=v_{2,3}$, which implies
that a vertex is incident with three $3$-edges, contradicting
Claim~\ref{claim20}. Hence, $v_{1,2}$, $v_{1,3}$ and $v_{2,3}$ are
distinct vertices. Thus, $e_1 = \{u_1,v_{12},v_{13}\}$, $e_2 =
\{u_2,v_{12},v_{23}\}$, and $e_3 = \{u_3,v_{13},v_{23}\}$. Let $h =
\{v_{1,2}, v_{1,3}, v_{2,3}\}$. Let $H'$ be obtained by deleting the
edges $e,e_1,e_2,e_3,f_1,f_2,f_3$ and vertices $u_1,u_2,u_3$ and
adding the $3$-edges $f_1'$, $f_2'$, $f_3'$ and $h$. By
Claim~\ref{claim11} and by construction, we note that if $R$ is a
subhypergraph contributing to $2b(H')+b^1(H')$, then $R$ must contain
the added $3$-edge $h$, implying that $2b(H')+b^(H') \le 2$. Suppose
that $S'$ is a $\tau(H')$-set. Since $|S' \cap V(h)| \ge 1$, we may
assume renaming vertices if necessary that $v_{1,2} \in S'$. But then
$S' \cup \{u_3\}$ is a transversal of $H$, and so $\tau(H) \le |S'| +
1 = \tau(H')+1$. Therefore,

\[
\begin{array}{lcl} \1
\phi(H)-\phi(H') & \ge & 6|\{u_1,u_2,u_3\}| + \omega(e) + \omega(e_1) + \omega(e_2) + \omega(e_3) +
 \omega(f_1) + \omega(f_2)  \\ \1
& & \hspace*{0.5cm} + \omega(f_3) - \omega(f_1') - \omega(f_2') - \omega(f_3') - \omega(h) - (2b(H') + b^1(H')) \\ \1
& \ge & (3 \times 6) + (4 \times 6) + (3 \times 4) - (4 \times 6) - 2 \\ \1
& > & 24 \\ \1
& \ge & 24(\tau(H) - \tau(H')),
\end{array}
\]

\noindent contradicting Fact~1. Hence, $\Vx{e_i} \cap \Vx{e_j} =
\emptyset$ for some $i$ and $j$ where $1 \le i < j \le 3$. Renaming
vertices if necessary, we may assume that $\Vx{e_1} \cap \Vx{e_2} =
\emptyset$. For $i \in \{1,2,3\}$, let $e_i=\{u_i,x_i,y_i\}$. Since
$H$ has no overlapping edges by Claim~\ref{claim12}, we know that
$|V(f_2) \cap V(e_1)| \le 1$. Renaming the vertices $x_1$ and $y_1$
if necessary, we may assume that $x_1 \notin \Vx{f_2}$. This implies
that there is no common edge containing both $u_2$ and $x_1$. We now
consider the hypergraphs $H^*=H(\{x_1,u_2\},\{u_1\})$. Then,
$e_2(H^*)=0$ and $H^*$ has no overlapping edges, implying that
$b(H^*)=b^1(H^*)=0$. By Claim~\ref{claim20}, the vertex $x$ is
contained in two $3$-edges, say $e_1$ and $e_x$, and in one $4$-edge,
say $f_x$. We now have that

\[
\begin{array}{lcl} \1
\phi(H)-\phi(H^*) & \ge & 6|\{u_1,u_2,x_1\}| + \omega(e) + \omega(e_1) + \omega(e_2) +
 \omega(f_1) + \omega(f_2)  \\ \1
& & \hspace*{0.5cm} + \omega(e_x) + \omega(f_x) - \omega(f_1') - (2b(H^*) + b^1(H^*)) \\ \1
& \ge & (3 \times 6) + (4 \times 6) + (3 \times 4) - 6 - 0 \\ \1
& = & 48 \\ \1
& = & 24|\{\{x_1,u_2\}|,
\end{array}
\]

\noindent contradicting Fact~1. This completes the proof of
Theorem~\ref{main_thm}.~\qed

\section{Proof of Theorem~\ref{main_thmC}} \label{3n7}

Before giving a proof of Theorem~\ref{main_thmC} we first present short proofs of Theorem~\ref{main_thmC_h1} and Theorem~\ref{main_thmC_h2}.
Recall the statement of Theorem~\ref{main_thmC_h1}, first proved by Chv\'{a}tal and McDiarmid.

\noindent \textbf{Theorem~\ref{main_thmC_h1}.} \emph{ {\rm (\cite{ChMc})}
If $H$ is a $4$-uniform hypergraph, then
$\tau(H) \le \frac{n(H)}{6} + \frac{m(H)}{3}$.
}

\noindent \textbf{Proof of Theorem~\ref{main_thmC_h1}.}
We will prove the theorem by induction on $n(H)$. Clearly the theorem holds when $n(H) \le 4$, so assume that $H$ is a $4$-uniform hypergraph with $n(H)>4$. Further since
$\tau(H)$ is additive with respect to the components of $H$,
we may assume that $H$ is connected, and so $\delta(H) \ge 1$. If $\Delta(H)>2$, then let $x$ be any
vertex with $d_H(x)>2$ and let $H'=H-x$. By induction, $6\tau(H') \leq n(H') + 2m(H') \leq (n(H)-1) + 2(m(H)-3)$.
By adding $x$ to any transversal in $H'$ we obtain a transversal in $H$, implying that $6\tau(H) \leq 6(\tau(H')+1) \leq (n(H)+2m(H) -7)+6 < n(H)+2m(H)$.
We may therefore assume that $\Delta(H) \leq 2$.

If some $x \in V(H)$ has $d_H(x)=1$, then let $e$ denote the edge containing $x$ and let $y$ be a vertex in $e$ of maximum degree in $H$.
Note that $d_H(y) = 2$ as $H$ is connected, $n(H) \ge 5$ and $\Delta(H) \leq 2$.
Consider the hypergraph $H'=H-y$ obtained from $H$ by deleting $y$ and all
edges containing $y$ and then removing any resulting isolated vertices.
Then, $n(H') \leq n(H) - 2$ (since both $x$ and $y$ get removed) and $m(H') = m(H) - 2$. By induction, $6\tau(H') \leq n(H') + 2m(H') \leq n(H) + 2m(H) -6$.
By adding $y$ to any transversal in $H'$ we obtain a transversal in $H$, implying that $6\tau(H) \leq 6(\tau(H')+1) \leq n(H')+2m(H')+ 6 \leq n(H)+2m(H)$, by induction.
We may therefore assume that $\delta(H) \geq 2$, implying that $H$ is $2$-regular.

Since $H$ is $2$-regular and $4$-uniform we have $2n(H)=4m(H)$, which, by Theorem~\ref{main_thmA} implies that
$24\tau(H) \leq 6n(H) + 4m(H) = 4n(H) + 8m(H)$, completing the proof of Theorem~\ref{main_thmC_h1}.~\qed

\medskip
\noindent \textbf{Theorem~\ref{main_thmC_h2}.} \emph{ {\rm (\cite{ThYe07})}
If $H$ is a $4$-uniform hypergraph, then
$\tau(H) \le \frac{5n(H)}{21} + \frac{4m(H)}{21}$.
}

\noindent \textbf{Proof of Theorem~\ref{main_thmC_h2}.}
The following holds for $4$-uniform hypergraph, by Theorem~\ref{main_thmA} and Theorem~\ref{main_thmC_h1}.

\[
\tau(H)  \leq   \frac{6}{7} \left( \frac{n(H)}{4} + \frac{m(H)}{6} \right) + \frac{1}{7} \left( \frac{n(H)}{6} + \frac{m(H)}{3} \right)
         =  \frac{5n(H)}{21} + \frac{4m(H)}{21}. \hspace*{1cm} \Box
\]

\medskip
Recall that for a graph $G$, the \emph{open neighborhood hypergraph}, abbreviated
ONH, of $G$ is the hypergraph $H_G$ with vertex set $V(H_G) = V(G)$
and with edge set $E(H_G) = \{ N_G(x) \mid x \in V \}$ consisting of
the open neighborhoods of vertices in $G$.  The transversal number of
the ONH of a graph is precisely the total domination number of the
graph; that is, for a graph $G$, we have $\gt(G) = \tau(H_G)$.
We are now in a position to prove Theorem~\ref{main_thmC}. Recall the
statement of the theorem.

\noindent \textbf{Theorem~\ref{main_thmC}.} \emph{ {\rm (\cite{ThYe07})} If $G$ is a graph
of order~$n$ with $\delta(G) \ge 4$, then $\gt(G) \le 3n/7$. }

\noindent \textbf{Proof of Theorem~\ref{main_thmC}.} Let $G$ be a
graph of order~$n$ with $\delta(G) \ge 4$ and let $H_G$ be the ONH of
$G$. Then, each edge of $H_G$ has size at least~$4$. Let $H$ be
obtained from $H_G$ by shrinking all edges of $H_G$, if necessary, to
edges of size~$4$. Then, $H$ is a $4$-uniform hypergraph with $n$
vertices and $n$ edges; that is, $n(H) = m(H) = n(G) = n$.
By Theorem~\ref{main_thmC_h2} we note that $21\tau(H) \leq 5n(H)+4m(H) = 9n$.
This completes the proof of the theorem since $\gt(G) = \tau(H_G)
\le \tau(H)$.~\qed

\section{Proof of Theorem~\ref{main_thmD}} \label{char:3n7} 

Chv\'{a}tal and McDiarmid proved the following bound in~\cite{ChMc}.

\begin{thm} {\rm (\cite{ChMc})}
If $H$ is a $4$-uniform hypergraph, then $6 \tau(H) \le n(H) +
2m(H)$.
 \label{thm_n2m6}
\end{thm}

In order to present a proof of Theorem~\ref{main_thmD}, we
shall need a characterization of the hypergraphs that achieve
equality in Theorem~\ref{thm_n2m6}. For this purpose, let $H_4$ be
the hypergraph on four vertices with only one hyperedge containing
all four of these vertices. Let $H_6$ be the hypergraph with vertex
set $\{a_1,a_2,b_1,b_2,c_1,c_2\}$ and edge set $E(H_6) = \{
\{a_1,a_2,b_1,b_2\}$, $\{a_1,a_2,c_1,c_2\}$, $\{b_1,b_2,c_1,c_2\}
\}$. The following result is given in~\cite{HeLo12}.

\begin{thm} {\rm (\cite{HeLo12})}
Let $H$ be a $4$-uniform hypergraph. If $6 \tau(H) = n(H) + 2m(H)$,
then every component of $H$ is isomorphic to $H_4$ or $H_6$.
 \label{thm_n2m6_char}
\end{thm}

We shall also need the following result in~\cite{HeYe08}.

\begin{thm} {\rm (\cite{HeYe08})}
The ONH of a connected bipartite graph consists of two components,
while the ONH of a connected graph that is not bipartite is
connected.
 \label{thm_ONH}
\end{thm}

We are now in a position to prove Theorem~\ref{main_thmD}. Recall the
statement of the theorem.

\noindent   \textbf{Theorem~\ref{main_thmD}.} \emph{If $G$ is a
connected graph  of order~$n$ with $\delta(G) \ge 4$, then $\gt(G)
\le 3n/7$. Furthermore we have equality if and only if $G$ is the
bipartite complement of the Heawood Graph.}

\noindent \textbf{Proof of Theorem~\ref{main_thmD}.} Let $G$ be a
connected graph  of order~$n$ with $\delta(G) \ge 4$ and let $H_G$ be
the ONH of $G$. If $G$ is not $4$-regular, then let $x$ be an
arbitrary vertex in $G$ with $d_G(x) \ge 5$. Now let $H$ be obtained
by shrinking all edges of size greater than four to size four in such
a way that we never remove $x$ from any edge. We note that the
resulting hypergraph $H$ is $4$-uniform with $n(H) = m(H) = n(G) =
n$, but $H$ is not $4$-regular. Alternatively if $G$ is $4$-regular,
then let $H=H_G$ in which case again $n(H) = m(H) = n(G) = n$, but in
this case $H$ is $4$-regular.

Let $x_1$ be a vertex of maximum degree in $H$. Let $x_2$ be a vertex
of maximum degree in $H-\{x_1\}$. Let $x_3$  be a vertex of maximum
degree in $H-\{x_1,x_2\}$. Continue this process as long as the
maximum degree in the resulting hypergraph is at least four and let
$X=\{x_1,x_2,\ldots,x_\ell\}$ be the resulting set of chosen
vertices. Let $H' = H-X$ and note that the following holds.
\begin{description}
\item[(a):] $\Delta(H') \le 3$.
\item[(b):] $n(H') \le n(H) - |X|$ and $m(H') \le m(H) - 4|X|$.
    Furthermore if $H$ is not $4$-regular, then since $x_1$
    removes at least five edges from $H$ we have that $m(H') <
    m(H) - 4|X|$.
\end{description}

If $|X|<n/7$, then by Theorem~\ref{main_thmA} we have
\[
\begin{array}{lcl} \1
\tau(H) & \le & \tau(H') + |X|  \\ \1
& \le & \frac{n(H')}{4}+ \frac{m(H')}{6} + |X| \\ \1
& \le & \frac{n(H)-|X|}{4} +\frac{m(H)-4|X|}{6} + |X| \\ \1
& = & \left( \frac{1}{4}+\frac{1}{6} \right)n + \left( 1 - \frac{1}{4} - \frac{4}{6} \right)|X|  \\ \1
& = & \frac{5n}{12} + \frac{|X|}{12}   \\ \1
& < & \left( \frac{5}{12} + \frac{1}{7 \times 12} \right) n \\    \1
& = & 3n,
\end{array}
\]

\noindent and the desired result follows from the observation that
$\gt(G) = \tau(H_G) \le \tau(H)$. Hence in what follows we may assume
that $|X| \ge n/7$. By Theorem~\ref{thm_n2m6}, we now have that

\[
\begin{array}{lcl} \1
\tau(H) & \le & \tau(H') + |X|  \\ \1
& \le & \left( \frac{n(H') + 2m(H')}{6} \right)  + |X| \\ \1
& \le & \left( \frac{n(H)-|X| + 2(m(H) - 4|X|)}{6} \right) +|X| \\ \1
& = & \frac{n}{2} - \frac{9|X|}{6} + |X|  \\ \1
& = & \frac{n - |X|}{2}   \\ \1
& \le & \frac{n - n/7}{2} \\    \1
& = & 3n.
\end{array}
\]

Hence, $\gt(G) = \tau(H_G) \le \tau(H) \le 3n/7$, proving the desired
upper bound. Suppose that $\gt(G) = 3n/7$. Then we must have equality
throughout the above inequality chains. In particular, this implies
that the following holds.
\begin{description}
\item[(c):] $6\tau(H') = n(H') + 2m(H')$.
\item[(d):] $n(H') = n(H) - |X|$.
\item[(e):] $m(H') = m(H) - 4|X|$.
\item[(f):] $|X| = n/7$.
\item[(g):] $H$ is $4$-regular (by (b) and (e)).
\end{description}

Since (c) holds, Theorem~\ref{thm_n2m6_char} implies that every
component of $H'$ is isomorphic to $H_4$ or $H_6$. By (d), (e) and
(f), and noting that $n(H) = m(H) = n$, we have that $n(H') = 6n/7$
and $m(H') = 3n/7$. Hence if $\dbar$ denotes the average degree in
$H'$, we have that

\[
4 m(H') = \sum_{v \in V(H')} d_{H'}(v) =  n(H') \cdot \dbar,
\]

and so $\dbar = 4m(H')/n(H') = 2$. We show that every component of
$H'$ is an $H_6$-component. Suppose to the contrary that there is an
$H_4$-component in $H'$. Each vertex in such a component has
degree~$1$ in $H'$. Since the average degree in $H'$ is~$2$, this
implies that there must also be a vertex of degree~$3$ in $H'$.
However such a vertex does not belong to an $H_4$- or an
$H_6$-component, a contradiction. Therefore the following holds.
\begin{description}
\item[(h):] Every component of $H'$ is an $H_6$-component.
\end{description}

Suppose that $N_H(u_1) \cap N_H(u_2) \ne \emptyset$ for some $u_1,u_2
\in V(H)$. We show that $u_1$ and $u_2$ are contained in a common
edge of $H$. Suppose to the contrary that no edge in $H$ contains
both $u_1$ and $u_2$ and let $w \in N(u_1) \cap N(u_2)$ be arbitrary.
Let $f_1,f_2 \in E(H)$ be chosen so that $\{u_1,w\} \subset V(f_1)$
and $\{u_2,w\} \subset V(f_2)$. By the $4$-regularity of $H$, we can
choose the set $X$ by starting with $x_1=u_1$ and $x_2=u_2$. We note
that with this choice of the set $X$, the vertex $w \in V(H')$. By
(h), the vertex $w$ belongs to some $H_6$-component in $H'$, implying
that there is a vertex $w' \in V(H')$ such that $w$ and $w'$ both
belong to two overlapping edges, say $e_1$ and $e_2$, in $E(H')$.
However, if we had created $X$ starting with $x_1=w'$, then $w$ would
belong to an $H_6$-component, $R$, of $H'$. Since $d_H(w) = 4$ and
$d_{H'}(w) = 2$, and since $e_1,e_2 \notin E(H')$, we have that
$f_1,f_2 \in E(H')$ and $f_1,f_2 \in E(R)$. In particular, $u_1$ and
$u_2$ are contained in a common edge of $R$ and therefore of $H$, a
contradiction. Therefore, the following holds.
\begin{description}
\item[(i):] If $N_H(u_1) \cap N_H(u_2) \ne \emptyset$ for some
    $u_1,u_2 \in V(H)$, then there exists an edge $e \in E(H)$,
    such that $u_1,u_2 \in V(e)$.
\end{description}

Let $u_1 \in V(H)$ be arbitrary. If some edge $e$ contains vertices
from $N[u_1]$ and from $V(H) \setminus N[u_1]$, then let $u_2 \in
V(e) \cap (V(H) \setminus N[u_1])$ and let $w \in V(e) \cap N[u_1]$
be arbitrary. Since $u_1$ and $u_2$ are not adjacent, $u_1 \ne w$ and
no edge contains both $u_1$ and $u_2$. However this is a
contradiction by (i). Therefore, the following holds.
\begin{description}
\item[(j):] If $x \in V(H)$, then $N[x]$ is the vertex set of
    some component in $H$.
\end{description}

Let $x \in V(H)$ be arbitrary and let $R_x$ be the component of $H$
containing~$x$. By (j), $V(R_x) = N[x]$. By the $4$-regularity of
$H$, we can choose the set $X$ by starting with $x_1=x$. Thus by (h),
$R_x - \{x\}$ only contains components isomorphic to $H_6$. However
since $H$ is a $4$-regular $4$-uniform hypergraph, and since $H_6$ is
$2$-regular, we note that $R_x-\{x\}$ must contain only one
component, which is isomorphic to $H_6$. Therefore, $|R_x|=7$ and
$R_x - \{x\} = H_6$. This is true for every vertex $x$ of $H$,
implying that $R_x$ must be isomorphic to the complement of the Fano
plane. Hence, the following holds.
\begin{description}
\item[(k):] Every component of $H$ is isomorphic to the
    complement of the Fano plane, which we will denote by $\overline{F_7}$.
\end{description}

By (k), every component of $H$ is isomorphic to  $\overline{F_7}$ (the complement of the
Fano plane). If $H \ne H_G$, then by construction $H$ is not $4$-regular, a contradiction to (g).
Hence, $H = H_G$. Since $\overline{F_7}$ is not the ONH of any graph, applying
the result of Theorem~\ref{thm_ONH} we have that $H$ consists of
precisely two components since $G$ is by assumption connected.
 Let $G'$ be constructed such that $V(G')=V(\overline{F_7}) \cup E(\overline{F_7})$ and
let $xy$ be an edge in $G'$ if and only if $x$ belongs to $y$ in $\overline{F_7}$ ($x$ is a vertex and $y$ is an edge in $\overline{F_7}$).
Now it is not difficult to see that $G'$ is the incidence bipartite graph of the complement of the Fano plane and that the ONH of
$G'$ is $H$. Therefore $G'=G$.~\qed

\section{Closing Comment}

Let $H$ be a $4$-uniform hypergraph of order $n = n(H)$ and size~$m = m(H)$. In this paper we have shown that
if $\Delta(H) \le 3$, then $\tau(H) \le n/4 + m/6$. It is known
that $\tau(H) \le n/4 + m/6$ is not always true when $\Delta(H) \ge 4$. We close with the following conjectures. Recall that a hypergraph is \emph{linear} if every two edges intersect in at
most one vertex.

\begin{conj}
If $H$ is a $4$-uniform linear hypergraph, then $\tau(H) \le \frac{n}{4} + \frac{m}{6}$.
 \label{conj2}
\end{conj}

\begin{conj}
If $H$ is a $4$-uniform linear hypergraph, then $\tau(H) \le \frac{n + m}{5}$.
 \label{conj3}
\end{conj}

We remark that Conjecture~\ref{conj3} implies Conjecture~\ref{conj2}. If there is a vertex of degree at least~$5$, then we may remove it and use induction in order to prove Conjecture~\ref{conj2} and if there is no such vertex we note that Conjecture~\ref{conj2} follows from Conjecture~\ref{conj3} as in this case $n/5 + m/5 < n/4 + m/6$.
Conjecture~\ref{conj3}, if true, would be best possible due to the $4$-uniform hypergraph $H_{10}$, illustrated in Figure~\ref{f:H10}, of order~$n = 10$, size~$m = 5$, and $\tau = 3$.

\bigskip
\tikzstyle{vertexX}=[circle,draw, fill=black!100, minimum size=8pt, scale=0.5, inner sep=0.1pt]
\newcommand{\SsD}[1]{}
\begin{figure}[htb]
\tikzstyle{vertexX}=[circle,draw, fill=black!100, minimum size=8pt, scale=0.5, inner sep=0.1pt]
\begin{center}
\begin{tikzpicture}[scale=0.3]
 \node (a1) at (3.0,1.0) [vertexX] {\SsD{v_1}};
\node (a2) at (1.0,4.0) [vertexX] {};
\node (a3) at (8.0,4.0) [vertexX] {};
\node (a4) at (1.0,7.0) [vertexX] {};
\node (a5) at (1.0,10.0) [vertexX] {};
\node (a6) at (8.0,10.0) [vertexX] {};
\node (a7) at (5.0,7.0) [vertexX] {};
\node (a8) at (5.0,10.0) [vertexX] {};
\node (a9) at (11.0,7.0) [vertexX] {};
\node (a10) at (5.0,4.0) [vertexX] {};
\draw[color=black!90, thick,rounded corners=4pt] (2.0499807070266725,9.99363487064001) arc (-0.3473307448246601:179.2952450519341:1.05); 
\draw[color=black!90, thick,rounded corners=4pt] (2.1306542800660906,0.4111553522085877) arc (-145.88855253567698:33.910731462853335:1.05); 
\draw[color=black!90, thick,rounded corners=4pt] (-0.049920569951011684,3.98708501669599) arc (180.7047549480659:213.2686875876366:1.05); 
\draw[color=black!90, thick,rounded corners=4pt] (2.0499807070266725,9.99363487064001) -- (2.0040741850589106,4.3071400835096005) -- (3.8714031954967276,1.5857955879640029);
\draw[color=black!90, thick,rounded corners=4pt] (2.1306542800660906,0.4111553522085877) -- (0.12208735591026942,3.424005738442319);
\draw[color=black!90, thick,rounded corners=4pt] (-0.049920569951011684,3.98708501669599) -- (-0.050000000000000044,7.0) -- (-0.049920569951011684,10.012914983304011);
\draw[color=black!90, thick,rounded corners=4pt] (5.749999847417262,10.000478407863328) arc (0.03654767108452628:180.2805186063686:0.75); 
\draw[color=black!90, thick,rounded corners=4pt] (5.625214161602007,3.5857449431421515) arc (326.4722994678067:359.9634523289155:0.75); 
\draw[color=black!90, thick,rounded corners=4pt] (2.378477767908561,1.419773885581363) arc (145.9650494508976:326.14756548023377:0.75); 
\draw[color=black!90, thick,rounded corners=4pt] (5.749999847417262,10.000478407863328) -- (5.749999847417262,3.9995215921366722);
\draw[color=black!90, thick,rounded corners=4pt] (5.625214161602007,3.5857449431421515) -- (3.6228562684849397,0.5822081034665504);
\draw[color=black!90, thick,rounded corners=4pt] (2.378477767908561,1.419773885581363) -- (4.282421395542916,4.2181305719641875) -- (4.250008988939981,9.996328034699204);
\draw[color=black!90, thick,rounded corners=4pt] (4.992814850727984,10.899971318226275) arc (90.45742566858571:89.79652738238184:0.9); 
\draw[color=black!90, thick,rounded corners=4pt] (8.634484929395986,10.6383015544156) arc (45.171808746367944:89.54254959264932:0.9); 
\draw[color=black!90, thick,rounded corners=4pt] (10.35499147420594,6.3723344826637565) arc (-135.78073383493228:44.8281912536321:0.9); 
\draw[color=black!90, thick,rounded corners=4pt] (0.99680382681625,10.899994324691539) arc (90.20347513216664:270.54736293902323:0.9); 
\draw[color=black!90, thick,rounded corners=4pt] (4.992814850727984,10.899971318226275) -- (8.007185537855243,10.899971315123839);
\draw[color=black!90, thick,rounded corners=4pt] (8.634484929395986,10.6383015544156) -- (11.6383015544156,7.634484929395986);
\draw[color=black!90, thick,rounded corners=4pt] (10.35499147420594,6.3723344826637565) -- (7.654620431594964,9.168908576792932) -- (1.008597826158542,9.100041069056289);
\draw[color=black!90, thick,rounded corners=4pt] (0.99680382681625,10.899994324691539) -- (5.003196133685564,10.899994324831809);
\draw[color=black!90, thick,rounded corners=4pt] (1.0068640423311093,4.799970552534827) arc (89.50839314811405:269.9719735745519:0.8); 
\draw[color=black!90, thick,rounded corners=4pt] (11.56689056583425,6.435522288865075) arc (-44.877806206814796:135.7304985509228:0.8); 
\draw[color=black!90, thick,rounded corners=4pt] (8.0003913227057,3.2000000957084183) arc (270.0280264254481:314.87780620681474:0.8); 
\draw[color=black!90, thick,rounded corners=4pt] (1.0068640423311093,4.799970552534827) -- (7.693244035398324,4.738850985098684) -- (10.42714848092762,7.558427378534099);
\draw[color=black!90, thick,rounded corners=4pt] (11.56689056583425,6.435522288865075) -- (8.564477711134925,3.43310943416575);
\draw[color=black!90, thick,rounded corners=4pt] (8.0003913227057,3.2000000957084183) -- (0.9996086772942984,3.2000000957084183);
\draw[color=black!90, thick,rounded corners=4pt] (1.0074859813489707,7.649956890942194) arc (89.34011597049562:270.6598840295044:0.65); 
\draw[color=black!90, thick,rounded corners=4pt] (8.649999943388577,10.000271283699341) arc (0.02391294072948469:135.63085566770644:0.65); 
\draw[color=black!90, thick,rounded corners=4pt] (7.535347910927709,3.5454689932240533) arc (224.36914433229356:359.9760870592705:0.65); 
\draw[color=black!90, thick,rounded corners=4pt] (1.0074859813489707,7.649956890942194) -- (4.7513218109650674,7.600549047371075) -- (7.535347910927709,10.454531006775946);
\draw[color=black!90, thick,rounded corners=4pt] (8.649999943388577,10.000271283699341) -- (8.649999943388577,3.999728716300659);
\draw[color=black!90, thick,rounded corners=4pt] (7.535347910927709,3.5454689932240533) -- (4.751321810965067,6.399450952628925) -- (1.0074859813489707,6.350043109057806);
 \end{tikzpicture}
\end{center}
\vskip -0.75 cm
\caption{The hypergraph $H_{10}$.} \label{f:H10}
\end{figure}
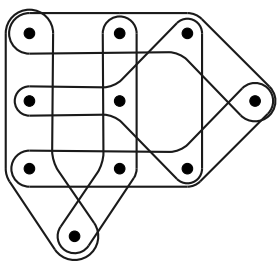


\end{document}